\documentclass[11pt]{amsart}
\usepackage{amsmath,amsfonts,amssymb,amsxtra,latexsym,amscd,enumerate,amsthm}

\usepackage{latexsym}
\usepackage{epsfig}
\usepackage{verbatim}
\def\z{\zeta}
\def\t{\theta}
\def\r{\mathcal{R}}

\def\g{\gamma}
\def\G{\Gamma}
\def\a{\alpha}

\def\d{\delta}
\def\b{\beta}
\def\e{\varepsilon}
\def\p{\Phi}
\def\v{\varphi}
\def\ep{\epsilon}
\def\es{\emptyset}
\def\l{\Lambda}
\def\ll{\Lambda}
\def\o{\omega}
\def\O{\Omega}
\def\les{\lesssim}
\def\n{\nu}
\def\R{\mathbb{R}}
\def\Q{\mathcal{Q}}
\def\c{\mathcal{C}}
\def\C{\mathbb{C}}
\def\Z{\mathbb{Z}}
\def\K{\mathcal{K}}
\def\D{\mathcal{D}}
\def\P{\mathbb{P}}
\def\p{\mathcal{P}}
\def\U{\mathcal{U}}
\def\V{\mathcal{V}}

\def\J{\mathcal{J}}
\def\I{\mathcal{I}}
\def\s{\mathbf{s}}
\def\M{\mathcal{M}}
\def\S{\mathcal{S}}
\def\H{\mathcal{H}}

\def\A{\mathcal{A}}
\def\B{\mathcal{B}}

\def\TT{\mathbb{T}}
\def\N{\mathbb{N}}

\def\n{\mathcal{N}}
\def\m{\mathcal{M}}
\def\f{\mathcal{F}}

\def\Z{\mathbb{Z}}
\def\beq{\begin{equation}}
\def\eeq{\end{equation}}
\def\eq{\Leftrightarrow}
\def\beq{\begin{equation}}
\def\eeq{\end{equation}}
\def\ra{\rightarrow}
\def\Leq{\Longleftrightarrow}

\def\z{\zeta}
\def\t{\theta}
\def\r{\rho}
\def\g{\gamma}
\def\G{\Gamma}
\def\a{\alpha}

\def\d{\delta}
\def\b{\beta}
\def\e{\varepsilon}
\def\p{\Phi}
\def\v{\varphi}
\def\va{\vartheta}
\def\ep{\epsilon}
\def\u{\upsilon}
\def\es{\emptyset}
\def\ll{\Lambda}
\def\l{\Lambda}
\def\o{\omega}
\def\O{\Omega}
\def\les{\lesssim}
\def\n{\nu}
\def\R{\mathbb{R}}
\def\Q{\mathcal{Q}}
\def\C{\mathbb{C}}
\def\c{\mathcal{C}}
\def\m{\mathcal{M}}
\def\Z{\mathbb{Z}}
\def\s{\mathfrak{s}}
\def\K{\mathcal{K}}
\def\D{\mathcal{D}}
\def\P{\mathbb{P}}
\def\p{\mathcal{P}}
\def\t{\mathcal{T}}
\def\U{\mathcal{U}}
\def\L{\mathcal{L}}
\def\V{\mathcal{V}}
\def\J{\mathcal{J}}
\def\M{\mathcal{M}}
\def\H{\mathcal{H}}
\def\A{\mathcal{A}}
\def\B{\mathcal{B}}
\def\TT{\mathbb{T}}
\def\N{\mathbb{N}}
\def\n{\mathcal{N}}
\def\f{\mathcal{F}}
\def\k{\kappa}
\def\Z{\mathbb{Z}}
\def\I{\mathcal{I}}
\def\beq{\begin{equation}}
\def\eeq{\end{equation}}
\def\eq{\Leftrightarrow}
\def\beq{\begin{equation}}
\def\eeq{\end{equation}}
\def\ra{\rightarrow}
\def\Leq{\Longleftrightarrow}

\newtheorem{thm}{Theorem}

\newtheorem{conj}[thm]{Conjecture}

\newtheorem{d0}[thm]{Definition}
\newtheorem{o0}[thm]{Observation}
\newtheorem{c0}[thm]{Corollary}
\newtheorem{t1}[thm]{Theorem}
\newtheorem{l1}[thm]{Lemma}
\newtheorem{p1}[thm]{Proposition}
\newtheorem{r1}[thm]{Remark}

\begin{document}
\title[A unified approach to three themes in harmonic analysis (I$\,\&\,$II)]{\Large{A unified approach to three themes in harmonic analysis (I$\,\&\,$II)}
\newline\\
\small{\noindent (I) The Linear Hilbert Transform and Maximal Operator along variable curves\\
\noindent (II) Carleson Type operators in the presence of curvature\\
\noindent (III) The bilinear Hilbert transform and maximal operator along variable curves}}

\author{Victor Lie}

\date{\today}
\address{Department of Mathematics, Purdue, IN 46907 USA}

\email{vlie@math.purdue.edu}

\address{Institute of Mathematics of the
Romanian Academy, Bucharest, RO 70700, P.O. Box 1-764, Romania.}

\thanks{The author was supported by the National Science Foundation under Grant No. DMS-1500958. The most recent revision of the paper was performed while the author was supported by the National Science Foundation under Grant No. DMS-1900801.}


\dedicatory{In loving memory of Elias Stein,\\whose deep mathematical breadth and vision for unified theories and fundamental\\ concepts have shaped the field of harmonic analysis for more than half a century.}

\maketitle

\begin{abstract}
In the present paper and its sequel \cite{lvUA3}, we address three rich historical themes in harmonic analysis that
rely fundamentally on the concept of \textit{non-zero curvature}. Namely, we focus on the boundedness properties of (I)
the \emph{linear} Hilbert transform and maximal operator along variable curves, (II) Carleson-type operators in the presence of curvature, and (III) the \emph{bilinear} Hilbert transform and maximal operator along variable curves.

Our Main Theorem states that, given a general variable curve $\g(x,t)$ in the plane that is assumed only to be \textit{measurable} in $x$ and to satisfy suitable non-zero curvature (in $t$) and non-degeneracy conditions, all of the above itemized operators defined along the curve $\g$ are $L^p$-bounded for $1<p<\infty$.

Our result provides a new and unified treatment of these three themes. Moreover, it establishes a unitary approach for both the singular integral and the maximal operator versions within themes (I) and (III).

At the heart of our approach stays a methodology encompassing three key ingredients: 1) discretization on the multiplier side that confines the phase of the multiplier to oscillate at the \emph{linear} level, 2) \emph{Gabor}-frame discretization of the input function(s) and 3) extraction of the \emph{cancelation} hidden in the non-zero curvature of $\g$ via $TT^{*}-$orthogonality methods and time-frequency \emph{correlation}.
\end{abstract}

\section{Introduction}\label{Intro}

This paper constitutes the first part of a study that is meant to present a new and unified approach to three distinct themes in harmonic analysis that focus on the boundedness properties of
\begin{enumerate}
\item[(I)] The \emph{linear} Hilbert transform and maximal operator along variable curves;

\item[(II)] Carleson-type operators in the presence of curvature;

\item[(III)] The \emph{bilinear} Hilbert transform and maximal operator along variable curves.
\end{enumerate}

The key underlying concept that governs all of the above topics is that of \emph{curvature}. This will of course be fundamental in terms of the methods that we develop in order to establish the relevant connections among the above three themes and put them under the same umbrella.

The generic formulation of the subject that we intend to address is given by the following:

$\newline$
\textbf{Main Problem.}(\textsf{Informal}) \textit{For each point $x\in\R$ we associate a curve $\Gamma_{x}=(t,\,-\g_x(t))$ in the plane, where here $t\in\R$ and
\beq\label{gam}
\g_{x}(\cdot):=\g(x,\cdot):\:\R\,\rightarrow\,\R\,,
\eeq
is a real function obeying some ``suitable" smoothness and non-zero curvature conditions in the $t$-parameter. Define now the variable family of curves in the plane $\G\equiv \{\Gamma_{x}\}_{\{x\in\R\}}$.
$\newline$
\indent\textsf{Task:} Under minimal\footnote{The main target is to achieve minimal regularity in the $x$-variable.} conditions on the curve family $\Gamma$, study the $L^p$-boundedness, $1\leq p\leq \infty$, of the following operators:}

\begin{itemize}
\item \textbf{the linear Hilbert transform along $\G$} \textit{defined as}
\beq\label{defHVT}
\eeq
$$H_{\G}\,:\: S(\R^2)\longrightarrow S'(\R^2)\,,$$
$$H_{\G}(f)(x,y):= \textrm{p.v.}\int_{\R}f(x-t,\,y+\g(x,t))\,\frac{dt}{t}\:;$$

\item \textbf{the (sub)linear maximal operator along $\G$} \textit{defined as}
\beq\label{defMVT}
\eeq
$$M_{\G}\,:\: S(\R^2)\longrightarrow L^{\infty}(\R^2)\,,$$
$$M_{\G}(f)(x,y):= \sup_{h>0}\frac{1}{2 h}\int_{-h}^{h}|f(x-t,\,y+\g(x,t))|\,dt\:;$$

\item \textbf{the $\g$ - Carleson operator} \textit{given by}
\beq\label{SWgeng}
\eeq
$$C_{\g}\,:\: S(\R)\longrightarrow L^{\infty}(\R)\,,$$
$$C_{\g}f(x):=\textit{p.v.}\,\int_{\R} f(x-t)\, e^{i\,\g(x,t)}\,\frac{dt}{t}\,;$$

\item \textbf{the $\g$ - maximal operator}\footnote{This operator is trivially dominated from above by the standard Hardy-Littlewood operator; hence the $L^p$-boundedness of this operator within the range $1<p\leq\infty$ is superfluous. However, we choose to mention it in this enumeration for two reasons: 1) an aesthetic one - it preserves a sense of symmetry, corresponding naturally to the maximal operators considered for themes (I) and (III); and 2) a questioning one - its structure includes an oscillatory behavior which makes $M_{\g}$ susceptible to better $L^1$-behavior (\textit{i.e.} end-point behavior) than its standard counterpart, though we will not analyze this last aspect in our present paper. For more on this, please see the Final Remarks section.} \textit{given by}
\beq\label{HLgeng}
\eeq
$$M_{\g}\,:\: S(\R)\longrightarrow L^{\infty}(\R)\,,$$
$$M_{\g}f(x):=\sup_{a>0}\,\left|\frac{1}{2a}\int_{-a}^{a} f(x-t)\, e^{i\,\g(x,t)}\,dt\right|\,;$$

\item \textbf{the bilinear Hilbert transform along $\G$} \textit{defined as}
\beq\label{defbHVT}
\eeq
$$H_{\G}^{\B}\,:\: S(\R)\times S(\R)\longrightarrow S'(\R)\,,$$
$$H_{\G}^{\B}(f,g)(x):= \textrm{p.v.}\int_{\R}f(x-t)\:g(x+\g(x,t))\:\frac{dt}{t}\:;$$

\item \textbf{the (sub)bilinear maximal operator along $\G$} \textit{defined as}
\beq\label{defbMVT}
\eeq
$$M_{\G}^{\B}\,:\: S(\R)\times S(\R)\longrightarrow L^{\infty}(\R)\,,$$
$$M_{\G}^{\B}(f,g)(x):= \sup_{h>0}\frac{1}{2 h}\int_{-h}^{h}|f(x-t)\:g(x+\g(x,t))|\:dt\:.$$
\end{itemize}
$\newline$

The problem above lies at the interface of several relevant and interconnected topics in the harmonic analysis of the plane:
\begin{itemize}
\item the study of singular linear/bilinear/maximal integral operators;

\item the boundedness of Carleson-type operators;

\item the problem of differentiability of functions along (smooth) variable vector fields.
\end{itemize}

Before describing the historical evolution of the above themes and their deep connections with our main problem, we make a detour and introduce our main results.

\subsection{Main results}\label{Mres}

In order to visualize the global picture and key message of our main theorem we purposely choose a first more informal presentation of it - in that we will not describe here the precise properties of the class of curves $\g$ for which it holds but defer this more technical aspect of the presentation until Section \ref{Curv}. Our Main Theorem will address the (sub)linear cases \eqref{defHVT} - \eqref{HLgeng} mentioned in our Main Problem, leaving the treatment of the bilinear Hilbert transform and its maximal analogue for the second part of our study in \cite{lvUA3}.

With these being said, we have:
$\newline$

\noindent \textbf{Main Theorem.}\label{Mainth}
\textit{Let  $\G\equiv\{\G_x\}_{x\in\R}$ be a family of twisted\footnote{We use the term “twisted” throughout to describe such a family of variable curves for which $\g\,:\R^2\,\rightarrow\,\R$ is a generic measurable function; the most salient point we wish to emphasize is that $\g$ need not necessarily split as an elementary tensor of the form $\g(x,t)=u(x)\,\tilde{\g}(t)$.} variable curves defined by $\G_x=(t,\gamma(x,t))$ with
$\g:\:\R^2\,\rightarrow\,\R$ measurable.
$\newline$
\indent Assume now that $$\g\in \mathbf{M}_x\mathbf{NF}_{t}\,,$$ that is, at an informal level\footnote{For the precise definition of the class $\mathbf{M}_x\mathbf{NF}_{t}$ please see Section \ref{Curv}.}, one has\footnote{Below, the class $C^{2+}$ simply means the standard $C^{2+\d}$ where here  $\d$ can be any number strictly greater than zero.}
\beq\label{mppst}
\eeq
\begin{itemize}
\item $\g_t(\cdot):=\g(\cdot,\,t)$ is \emph{$x-$measurable} for every $t\in\R\setminus\{0\}$;
\item $\g_x(\cdot):=\g(x,\,\cdot)$ is finitely piecewise \emph{$t-$smooth} within the class $C^{2+}(\R\setminus\{0\})$ for almost every $x\in\R$;
\item $\g$ is ``non-flat" near the origin and infinity; [In particular, outside of a controlled region, $\g$ can be decomposed into a finite number of pieces on which it has $x-$uniform non-vanishing curvature in the $t-$variable.]
\item $\g$ obeys a suitable non-degeneracy condition.
\end{itemize}
\indent Then, for any $1<p<\infty$, we have that\footnote{Of course, in the case of the maximal operators $M_{\G}$ and $M_{\g}$ one trivially gets the desired bounds for the limiting case $p=\infty$}.
\begin{enumerate}
\item[(I)] $H_{\G}$ and $M_{\G}$ are bounded operators from $L^p(\R^2)$ to $L^p(\R^2)$;
\item[(II)] $C_{\g}$ and $M_{\g}$ are bounded operators from $L^p(\R)$ to $L^p(\R)$.
\end{enumerate}}
$\newline$

Next, we clarify the extent of generality of our Main Theorem which constitutes in fact one of the main motivations for our program - to unify and extend the existing results in the literature treating operators along curves $\g$ that are polynomial in $t$ with measurable coefficients in $x$:

\begin{t1}\label{Gencurvpolyn} Let $d\in\N$ and \footnote{Throughout the paper, for convenience, we allow a notational abuse and introduce the  following convention: given $\a,\,t\in \R$ we let $t^{\a}$ stand for either $|t|^{\a}$ or $\textrm{sgn}\, (t)\, |t|^{\a}$.}
\beq\label{rich0}
\g(x,t):=\sum_{j=1}^{d} a_j(x)\,t^{\a_j}\,,
\eeq
where here $\{\a_j\}_{j=1}^d\subset \mathbb{R}\setminus\{0,\,1\}$ and $\{a_j\}_{j=1}^d$ measurable functions.

Then, one has that \footnote{This result remains true for more general classes of curves - see for this Observation \ref{W} at the end of Section \ref{Curv}.}
\beq\label{rich}
\g\in \mathbf{M}_x\mathbf{NF}_{t}\,.
\eeq
\end{t1}

From the theorems above and their corresponding proofs we deduce the following consequences:

\begin{c0}\label{Polyncasegen} Let $\g$ be as in \eqref{rich0} and assume wlog that  $\{\a_j\}_{j=1}^d\subset \mathbb{R}\setminus\{0,\,1\}$ strictly increasing.  Then, letting $\G=(t,-\gamma(x,t))$ and $1<p<\infty$, there exists $C(p,\,d,\, \{\a_j\}_{j=1}^d)>0$ such that
\beq\label{lpbd}
\|C_{\g}\|_{L^p(\R)\rightarrow L^p(\R)},\:\|M_{\g}\|_{L^p(\R)\rightarrow L^p(\R)}\leq C(p,\,d,\, \{\a_j\}_{j=1}^d)\:,
\eeq
and
\beq\label{lpbd1}
\|H_{\G}\|_{L^p(\R^2)\rightarrow L^p(\R^2)},\:\|M_{\G}\|_{L^p(\R^2)\rightarrow L^p(\R^2)}\leq C(p,\,d,\, \{\a_j\}_{j=1}^d)\:.
\eeq
Moreover, the constant $C(p,\,d,\, \{\a_j\}_{j=1}^d)$ depends only on the quantities
\beq\label{lpbddep}
C(p,\,d,\, \{\a_j\}_{j=1}^d)=C\left(p,\,d,\,\max_{1\leq j\leq d}|\a_j|,\,\max_{1\leq j\leq d}\frac{1}{|\a_j|},\,\max_{1\leq j\leq d}\frac{1}{|\a_j-1|},\,\max_{1\leq j<k\leq d}\frac{1}{|\a_k-\a_j|}\right)\:.
\eeq

It follows from the above that if
\beq\label{condunif}
 \min_j |\a_j|,\,\min_j |\a_j-1|,\,\min_{j\not=k} |\a_j-\a_k|\gtrsim 1\,,
\eeq
then the $L^p$ bounds appearing in \eqref{lpbd} and \eqref{lpbd1} are all \underline{uniform} in $\{a_j\}_{j=1}^d$ and depend only on $\sum_{j=1}^d |\a_j|$ and $p$.
\end{c0}

\begin{c0}\label{Polyncaseap}
Let $1<p<\infty$, $d\in\N$ and $\{\a_j\}_{j=1}^d\subset \mathbb{R}\setminus\{1\}$.
Then, the \emph{generalized Polynomial Carleson-type operator} defined as
\beq\label{SWgen}
C_{\vec{\a},d}f(x):=\sup_{\{a_j\}_{j=1}^d\subset\R}\,\left|\,p.v.\,\int_{\R} f(x-t)\, e^{i\,\sum_{j=1}^d a_j\,t^{\a_j}}\,\frac{dt}{t}\right|\,,\:\:\:\:\:\:\:f\in \mathcal{S}(\R)\,,
\eeq
is a bounded operator from $L^p(\R)$ to $L^p(\R)$.
\end{c0}

\begin{o0}
Notice that if one takes in Corollary \ref{Polyncaseap} the particular values $\a_j=j+1$ for $j\in\{1,\ldots, d\}$ one recovers the by now classical  Stein-Wainger result, \cite{sw}, on the $L^p$ boundedness of the Polynomial Carleson operator with no linear term in one dimension. In a different direction, taking $d=1$ and thus placing ourselves in the simplified tensor-product case, we get the main result in \cite{Guoosc}.
\end{o0}

\begin{c0}\label{Crossprod} Assume $\g$ has a \underline{tensor-product} structure, that is
\beq\label{crp}
\g(x,t)=u(x)\,\tilde{\g}(t)\,,
\eeq
with
\begin{itemize}
\item $u$ real measurable;

\item $\tilde{\g}\in\mathbf{NF}$. Informally, membership in the class $\mathbf{NF}$ means that $\tilde{\g}\in C(\R\setminus{0})$ is doubling and ``non-flat" near the origin and infinity.

     Formally, the following two conditions hold \footnote{The class of curves $\mathbf{NF}$ serves as an extension of the class $\n\f$ introduced by the author in \cite{lv4}.}:
\begin{itemize}
\item There exists $j_0\in\N$ such that letting $\D:=\R\setminus \{ 2^{-j_0}<|t|<2^{j_0}\}$ one has $\tilde{\g}\in C^{2+}(\D)$ with $|\tilde{\g}''|>0$ on $\D$ and
\beq\label{crp1}
\sup_{c>0}\#\{j\in\Z,\,|j|\geq j_0\,|\,|2^{-j}\,\tilde{\g}'(2^{-j})|\in[c,2c]\}<\infty\:;
\eeq
\item Let $I:=\{s\,|\,\frac{1}{10}\leq|s|\leq 10\}$ and $j\in \Z$ with $|j|\geq j_0$ and set\footnote{From the first item above we notice that by possibly shrinking the domain $\D$ one can assume wlog that in fact also $|\tilde{\g}'|>0$ on $\D$.}
\beq\label{crasymptotic0}
\tilde{Q}_{j}(t):=\frac{\tilde{\g}(2^{-j}\,t)}{2^{-j}\,\tilde{\g}'(2^{-j})}\in C_t^{2+}(I)\,.
\eeq
Then, uniformly in $|j|\geq j_0$, one has that
\beq\label{crfstterma0}
 \sup_{t\in I} |\tilde{Q}_{j}(t)|<c_1(\g)\:,
\eeq
and
\beq\label{crfstterm0}
 \inf_{t\in I} |\tilde{Q}''_{j}(t)|>c_0(\g)\:,
\eeq
where  $0<c_0(\g)\leq c_1(\g)$ are constants depending only on $\g$.
\end{itemize}
\end{itemize}
Then
\beq\label{gamg}
\g\in \mathbf{M}_x\mathbf{NF}_{t}\:,
\eeq
if one of the below conditions is satisfied:

\noindent i) $u(x)=u_0$ for almost every $x\in\R$ (i.e. $u$ constant), or

\noindent ii) more generally, if $u(x)$ a generic measurable function, then the curve $\g$ satisfies the non-degeneracy condition

\beq\label{nondeg}
\exists\:M>0\:\:\textrm{s.t.}\:\:\forall\:|j|\geq |j_0|\:\:\:\left|\frac{\tilde{Q}'_{j}}{\tilde{Q}''_{j}}(t)- \frac{\tilde{Q}'_{j}}{\tilde{Q}''_{j}}(s)\right|\geq M\,|t-s|\,,\:\:\:\:\:\:\forall\:t,s\in I\:.
\eeq
\end{c0}

The next corollary was treated as a model case for our Main Theorem, point (I), and was proved\footnote{This result was shown in an unpublished note in May 2016 and will be made available soon - in the form of \cite{lvBLC} - in order to provide a better understanding of the contrast with our current approach. Strictly speaking the result proved in \cite{lvBLC} covers only the case $p=2$. However, this was intentionally meant in order to maintain the simplicity and reduce the size of the presentation in \cite{lvBLC}, leaving the general $L^p$ discussion for the present, more elaborate study.} by the author in \cite{lvBLC}.

\begin{c0}\label{Crossprod1}  Consider the class of curves $\n\f$ introduced in \cite{lv4}. Then,
taking  $\tilde{\g}\in \n\f$ and $\g$ as in \eqref{crp}, one has that the corresponding $H_{\G}$ and $M_{\G}$ are bounded operators from $L^p(\R^2)$ to $L^p(\R^2)$ for $1<p<\infty$.
\end{c0}

\begin{o0} i) Remark that  $\n\f$ is included in the set of curves that belong to $\mathbf{NF}$ and simultaneously satisfy \eqref{nondeg}; thus, Corollary \ref{Crossprod1} is an immediate consequence of Corollary \ref{Crossprod}.

ii) Recall from \cite{lv4} that the above set $\n\f$ contains, as particular instances, both the class of curves introduced in \cite{LY1} as well as any generalized polynomial of the form $\tilde{\g}(t)=\sum_{j=1}^{d} a_j\,t^{\a_j}$ with $\{a_j\}_{j=1}^d\subset \mathbb{R}$, $\{\a_j\}_{j=1}^d\subset \mathbb{R}\setminus\{0,\,1\}$, and $d\in\N$. Consequently, one deduces that the results in \cite{LY1} as well as the corresponding ones addressing the monomial case $\g(x,t)=u(x)\,t^{\a}$ with $\a\in (0,\infty)\setminus \{1\}$ in \cite{GHLR} follow immediately from Corollary \ref{Crossprod1}.
\end{o0}

We end our section with the following

\begin{o0}\label{bilinearhilbert}
In the second part of our study, \cite{lvUA3}, - under suitable conditions imposed on $\g$ - following and further developing some of the key ideas introduced here we will prove the $L^p$ boundedness of the bilinear Hilbert transform and maximal operator corresponding to \eqref{defbHVT} and \eqref{defbMVT}. Thus, we are able to provide a unified method for all the operators defined in the statement of our Main Problem, and, in particular, to identify and highlight as natural a common approach to both the singular and the maximal operators within the themes (I) and (III). Moreover, as a consequence of these methods, we are able to immediately encompass and generalize the previous results appearing in \cite{li}, \cite{LX}, \cite{lv4} and \cite{lv10}.
\end{o0}

\subsection{Main ideas and relevance of the results}

By the very nature of this problem that involves the study of a singular/maximal operator with rough $x-$dependence and highly oscillatory phase multiplier, the key stepping stone is given by the \textit{discretization} of the operator. This is the central element that dictates both the proof's mechanism and, finally, the strength of the output.

In the present context the discretization of our operator follows three levels\footnote{The description below reflects the treatment of the linear Hilbert and maximal operators along curves defined by \eqref{defHVT} and \eqref{defMVT}. However, at the conceptual level, the same philosophy applies to the other classes of operators introduced in our Main Problem.}:
\begin{itemize}
\item The first level focuses on the global analysis of the operator's multiplier, which is performed in Section \ref{anmult}. Inspired\footnote{For more on this please see Section \ref{conBHT}.} by the author's work in \cite{lv4} - though in our present context the extra rough $x-$dependence of the multiplier makes this analysis significantly more subtle - we decompose our multiplier into three components:
 \begin{itemize}
   \item \emph{a low-frequency part} - this addresses the situation in which there is virtually no oscillation of the phase of the multiplier and can be dealt with via Taylor series arguments;

   \item \emph{a high-frequency far from diagonal part} - this focuses on the region in the time-frequency plane where we have no stationary points of the phase of the multiplier and hence one uses the oscillation to integrate by parts and obtain supplementary decay. A powerful new tool in this context is given by Lemma \ref{translk}, which appeals to multiple ingredients such as shifted square functions and vector-valued Calder\'on-Zygmund inequalities and that results in a common treatment of the pieces corresponding to both the singular and the maximal operator. More importantly, Lemma \ref{translk} provides an alternative and unified approach to \emph{global} $L^p$-bounds, $p\not=2$, for $H_{\Gamma}$ and $M_{\Gamma}$ - this latter aspect being discussed in Section \ref{LPS}.

   \item \emph{a high-frequency close to diagonal part} - this is of course the most difficult component to treat, as it refers precisely to the region in the time-frequency plane where the multiplier's phase has stationary points.
 \end{itemize}

\item The second level focuses on the fine analysis of the high-frequency piece close to diagonal of the multiplier which is performed in two stages:
\begin{itemize}
\item first one needs to address the rough $x-$dependence of the phase of the multiplier - this is the content of \eqref{xlocfor} via the decomposition offered by \eqref{xloc}.
\item the second is to target the $(\xi, \eta)$ - frequency localization of the phase, which has as a key output a \emph{linearizing} effect on the oscillation of the phase; this operation is performed in \eqref{xiloc} and \eqref{etaloc}, and its specific choice is of fundamental importance for the $L^2$-decay result provided by Theorem \ref{l2dec}.
 \end{itemize}

\item Finally the third level focuses on the operator's input - that is, on the function $f$ itself - and appeals to an adapted Gabor frame decomposition of $f$ which of course needs to be synchronized with the discretization of the multiplier (for this, one is invited to see Step 2 in the proof of Theorem \ref{l2dec}).  This final element brings into the picture genuine manifestations of wave-packet analysis and provides the required smoothness for the kernel of our operator and hence for a successful application of the $TT^{*}-$method.
\end{itemize}

With the stage thus set, our approach is as follows:
\begin{itemize}
\item For the $L^2$ case one proves via orthogonality methods and time-frequency correlation a suitable \emph{exponential decay} in terms of the height of the multiplier's phase - this is the content of Theorem \ref{l2dec}, which is the central result of the paper.

\item For the $L^p$ case with $p\not=2$ one first proves \emph{tame polynomial growth} relative to the hight of multiplier's phase followed in a second stage by a standard interpolation argument. The tame $L^p-$bounds are the content of Theorem \ref{hmdiagp}. A noteworthy aspect here is that we dedicate an entire section - see Section \ref{otherlp} - to various other approaches to Theorem \ref{hmdiagp} in which one explores ideas that involve further discretization techniques, shifted maximal and square function estimates, etc.
\end{itemize}

The strategy discussed above may be summarized by what we would like to refer from now on as the \emph{$LGC-$methodology} encompassing the following three key ingredients:\footnote{For a more details on the specifics of this methodology one is invited to consult Section \ref{ldecsec}.}
\begin{itemize}
\item L - \emph{phase Linearization}: a first stage of dicretization is performed on the multiplier side with the aim of forcing the multiplier's phase to oscillate at the linear level;

\item G - \emph{Gabor frame discretization}: a second stage decomposition involving wave-packet analysis is applied at the level of the input functions;

\item C - \emph{Cancellation/Correlation}: the third and final stage relies on extracting the cancellation encoded in the non-zero curvature of the phase based on $TT^{*}$-methods and the time-frequency correlation of the variables involved in the decompositions employed at the first two stages.
\end{itemize}

Beyond these, one may notice the new non-degeneracy condition part of the definition of the class of curves $\mathbf{M}_x\mathbf{NF}_{t}$ - see \eqref{ndeg0} - which seems to be a very malleable notion that can be verified for a very wide classes of curves as envinced by Theorem \ref{Gencurvpolyn}, Corollaries \ref{Crossprod} and \ref{Crossprod1} and Observation \ref{W}.

We end this section by briefly discussing the relevance of our result in the context of the present literature, which manifests in three directions:

\begin{itemize}
\item This is a first study for rough\footnote{\textit{I.e.} no smoothness in the $x$-variable is assumed.} $x-$dependent curves that presents a novel and unified approach to both the Hilbert transform (singular integral) and the maximal operator cases.

\item While treating a situation of non-zero curvature in $t$, we introduce \emph{elements of wave-packet analysis (Gabor frames)} and then blend them together with orthogonality methods, with the aim of developing a theory that unifies/fills the existing gap between the two different approaches corresponding to the zero/non-zero curvature cases.\footnote{For more on the proof dichotomy between the zero and non-zero curvature cases the reader is invited to consult Section \ref{Dich}.} Although our treatment does not yet cover the zero-curvature case (\textit{e.g} $\g(x,t):=\sum_{j=1}^{d} a_j(x)\,t^{j}$ with $d\in\N$ and $\{a_j\}_{j=1}^d$ measurable functions - thus allowing the linear term $a_1(x)\,t$), the methods developed here provide some intuition for the general situation; we hope to return to this topic in the near future.

\item The proof of the present result provides a unified perspective on several important directions within Harmonic Analysis, gathering under the same umbrella themes involving maximal and singular oscillatory integrals of Stein-Wainger type, (see \cite{sw}); more generally, (Polynomial) Carleson-like operators (with no linear term); boundedness of the bilinear Hilbert and maximal operator along ``non-flat" curves (see \cite{li}, \cite{LX}, \cite{lv4}, \cite{lv10} and \cite{GL}); and the boundedness of Hilbert transforms and maximal operators along variable ``non-flat" curves (see \cite{GHLR}, \cite{LY1}).
\end{itemize}

\subsection{A fundamental dichotomy: curvature versus modulation invariance symmetry}\label{Dich}

Before passing to the historical evolution of our three distinct but inter-related themes, we make a brief digression in order to evoke a fundamental dichotomy that serves as a cornerstone in the field of harmonic analysis: \emph{non-zero versus zero curvature} problems.\footnote{While many of the elements in this section are part of the harmonic analysis folklore, we choose to present them here in order to provide a perspective on our main themes of study.} This dichotomy marks each of our chosen themes and is quintessential to identifying the method of proof in a given problem\footnote{Below we are taking as reference point the singular integral variants; however, all of the discussion below has a direct analogue for the maximal case.}:
$\newline$

\begin{itemize}
\item (A) \underline{\emph{Hilbert transform along curves}}:  We consider in \eqref{defHVT} a generic class of curves with $\g(x,t)=\sum_{j=1}^{n} a_j(x)\,t^j$ and $\{a_j(\cdot)\}_j$ arbitrary real measurable functions. Then, one has:
\begin{itemize}
\item the zero-curvature case (prototype: $n=1$, with $\g(x,t)=a_1(x) t$).

In this situation, letting $M_{1,a} f(x,y):= e^{i a x}\, f(x,y)$, one has that
\beq\label{hilbcsym}
\|H_{\Gamma}M_{1,a} f\|_{L^2(\R^2)}=\|H_{\Gamma}f\|_{L^2(\R^2)}\:.
\eeq
\item the nonzero-curvature case (prototype: $n>1$, with $\g(x,t)=\sum_{j=2}^{n} \tilde{a}_j(x)\,t^j$ - no linear term allowed).

In this situation $H_{\Gamma}$ has \emph{no} modulation invariance symmetry.
\end{itemize}
$\newline$

\item (B) \underline{\emph{Carleson-type operators}}: We consider here polynomial Carleson-type operators, which - following Kolmogorov's linearization - can be written in the form
\beq\label{polcarlm}
C_{\g}f(x):=\int_{\R}e^{i\,\g(x,t)}\,f(x-t)\,\frac{dt}{t}\:,
\eeq
with $n\in\N$, $\g(x,t):=\sum_{j=1}^{n} a_j(x)\,t^j$ and $\{a_j(\cdot)\}_j$ arbitrary real measurable functions. Now, analogously with the example above, we have:
\begin{itemize}
\item the zero-curvature case (prototype: $n=1$, \emph{i.e.} $\g(x,t)=a_1(x) t$).

In this situation, in addition to the standard commutation relations with translation and dilation symmetries the operator $C_{\g}$ is invariant under the modulation symmetry $M_a f(x):= e^{i a x}\, f(x)$ with $a\in\R$, \textit{i.e.}:
\beq\label{polcarlmsym}
C_{\g}M_a f=C_{\g}f\:.
\eeq

\item the nonzero-curvature case (prototype: $n>1$, with $\g(x,t)=\sum_{j=2}^{n} a_j(x)\,t^j$ - no linear term allowed).

In this situation $C_{\g}$ has \emph{no} modulation invariance symmetry.
 \end{itemize}
$\newline$

\item (C) \underline{\emph{Bilinear Hilbert transform along curves}}:  Taking the generic case $\g(x,t)=\g(t)=\sum_{j=1}^{n} a_j\,t^j$ with $\{a_j\}_j$ real and $\Gamma=(t,-\g(t))$, we define
\beq\label{biHilb}
H_{\Gamma}^{\B}(f,g)(x):=\int_{\R}f(x-t)\,g(x+\g(t))\,\frac{dt}{t}\:.
\eeq

\begin{itemize}
\item the zero-curvature case (prototype: $n=1$, with $\g(t)=a_1 t$ with $a_1\in\R\setminus\{0,\,1\}$);

In this situation we have that
\beq\label{bihilbcsym}
H_{\Gamma}^{\B}(M_{a_1}f,\,M_1 g)=M_{1+a_1} H_{\Gamma}^{\B}(f,g)\:.
\eeq
\item the nonzero-curvature case (prototype: $n>1$, with $\g(t)=\sum_{j=2}^{n} a_j\,t^j$ - no linear term allowed);

In this situation $H_{\Gamma}^{\B}$ has \emph{no} modulation invariance symmetry.
\end{itemize}

\end{itemize}

Once we have seen this dichotomy it is important to stress the following:

\begin{itemize}
\item In the \emph{zero-curvature (flat) case} all the above operators obey suitable invariance under modulation symmetry. As a consequence, any method of proof requires an approach based on wave-packet analysis and thus in particular a time-frequency discretization of the corresponding operator; moreover, the proof should involve concepts like mass and/or energy of wave-packets in the spirit of the known proofs of Carleson's Theorem (see \cite{c1}, \cite{f} and \cite{lt3}).

\item In the \emph{nonzero-curvature (non-flat) case} there is no modulation-invariance symmetry, and thus one expects that more standard analysis can be performed on the object under study, involving $T T^{*}$ and more general orthogonality methods, the (non)stationary phase principle including Van der Corput estimates, Littlewood-Paley techniques, square-function arguments, etc. While discretization techniques in physical and frequency space are still relevant, the zero frequency plays a favorite role in this discretization, and, usually, one is able to obtain a suitable \textit{scale type decay} where here the concept of ``scale" should be properly adapted to the context.
\end{itemize}

While both situations are interesting and historically motivated, generically speaking the zero-curvature situation tends to be more difficult and accordingly most of the celebrated problems in this area - some of which remain open - regard precisely this case. The situation of nonzero curvature can also prove challenging, but to a lesser extent. In this context, while often regarded as model problems for the flat case, the corresponding non-flat case problems usually can only provide limited intuition, since, they require yet distinct methods of proof.

This last fact motivates a very interesting further direction of study - that of striving to \emph{unify} the two approaches corresponding to the zero/non-zero curvature cases, and thus to provide a method of proof for the situation in which $\g$ is given by a polynomial in $t$ with the linear term included.

At this point, it is worth mentioning that, with the notable exception of the Polynomial Carleson operator proved in \cite{lv3}, no unified treatment is known for the other two fundamental objects: the Hilbert and bilinear Hilbert transform - and their maximal analogues - along curves.

Finally, this paper and its companion \cite{lvUA3} can be regarded as a first step into this program by providing a unifying treatment for the non-zero curvature cases of all three themes - and their maximal variants - enumerated at the beginning of our Introduction. In fact, our approach here goes further, by also partly incorporating elements of time-frequency analysis via Gabor frame decompositions. Isolating now the first theme represented by $H_{\Gamma}$ and $M_{\Gamma}$, which is also the main focus in our present paper, we hope to return with a unifying approach for the zero/non-zero curvature cases in the near future.

\subsection{Historical background; motivation (I)}

The first theme of our paper has a long and rich history, and, as is the case with many others in harmonic analysis, it originates in the field of partial differential equations. Thus, in what follows, we will start by describing the original PDE motivation for our theme (I) and its initial development within the field of harmonic analysis.

\subsubsection{\underline{The original motivation - a PDE quest}}\label{pdeorig}  As already noted, the study of the boundedness  of the Hilbert transform and maximal operators along curves is part of a larger class of deep and fundamental topics - see \emph{e.g.} the pointwise convergence of Fourier Series, the Calder\'on-Zygmund theory, and even the Restriction Problem visualized via Strichartz estimates - that traces back to the area of PDE. More precisely, the relevant starting point in our story pertains to the study of constant coefficient differential operators. Here, in order to provide more context to our description, we choose a parallel and contrasting presentation with the corresponding moment of birth for the Calder\'on-Zygmund theory.

\medskip
\noindent A. \emph{Constant coefficient elliptic differential operators}.
\medskip

In this section we take as a main prototype the following:
\medskip

\noindent \underline{\textsf{Model:}} The Laplace/Poisson equation in $\R^d$, $d\geq 2$, given by
\beq\label{lap}
\triangle u=f\,,
\eeq
where $u,\, f$ are suitable (smooth) functions.
\medskip

\noindent \underline{\textsf{Aim:}} Understand the $L^p$-boundedness, $1<p<\infty$, of the second derivatives of our solution $u$ in terms of the $L^p$ bounds of the input function $f$.
\medskip

Applying now standard PDE techniques, one obtains that the fundamental solution $U^0$ associated to \eqref{lap}, that is, the solution to
\beq\label{lapf}
\triangle U^0=\d_0\,,
\eeq
with $\d_0$ the Dirac mass at the origin, is given by
\beq\label{lapf1}
\eeq
$$U^0(x):=-\frac{1}{2\pi}\,\log\frac{1}{|x|}\:\:\textrm{if}\:\:d=2\,,$$
$$U^0(x):=\frac{1}{(d-2)\o_d}\,|x|^{2-d}\:\:\textrm{if}\:\:d>2,\:\:(\o_d=Area(S^d))\,.$$
Thus, for suitable $f$, the solution to \eqref{lap} becomes
$$u(x):=\int_{\R^d}U^0(x-y)\,f(y)\,dy\,.$$

Indeed, letting $k_{ij}:= U^{0}_{x_ix_j}$, we have that a.e.
$$ u_{x_i x_j}(x)=\frac{1}{d}\d_{ij} f(x)+\int_{\R^d} k_{ij}(x-y)\,f(y)\,dy\,,$$
and thus one can verify \eqref{lap} pointwise or in the sense of distributions.

For $d>2$ one has that the kernel $K:=k_{ij}$ obeys the following key properties:
\beq\label{CZ}
\eeq
\begin{itemize}
\item $K$ is homogeneous of degree $-d$, \textit{i.e.} if $\d_{\a}(x) = (\a x_1,\ldots,\, \a x_d)$ then
$K(\d_\a(x))=\a^{-d}\,K(x)\,,\:\:\:\a>0\:.$
\item $K$ is $C^{\infty}$ away from the origin \textit{or} alternatively, one can relax this assumption by requesting for example that
 $$\sup_{x\in\R^n\setminus\{0\}}\int_{|y|>2|x|} |K(y-x)-K(y)|\,dy<\infty\,;$$
\item $\int_{|x|=1} K(x)\,d\sigma(x)=0$.
\end{itemize}

Now, as it turns out, these particular conditions on the kernel $K$ associated with an operator $Tf:=K*f$ provide the right framework for developing a seminal new theory of such integral operators more generally, beyond the study of the Laplace and Poisson equations. This fundamental theory developed by Calder\'on and Zygmund in \cite{CZ1} and \cite{CZ2} had as its main result the fact that any such $Tf:=K*f$ acts boundedly from $L^p(\R^d)$ to $L^p(\R^d)$ for any $1<p<\infty$.  As an immediate consequence, we thus deduce that
$$\|u_{x_ix_j}\|_{L^p(\R^d)}\lesssim_{p} \|f\|_{L^p(\R^d)},\:\:\:\:\:1<p<\infty\;,$$
as desired.

\medskip
\noindent B. \emph{Constant coefficient parabolic differential operators}.
\medskip

Proceeding in the mirror with point A above, we now take as a prototype for our discussion the following:
\medskip

\noindent \underline{\textsf{Model:}} The heat equation in $\R^{d+1}_{+}=\R^d\times R_{+}$, $d\geq 2$; this is represented by
\beq\label{heat}
\partial_t u-\triangle u=f\,,
\eeq
where, as before, $u,\, f$ are suitable (smooth) functions.
\medskip

\noindent \underline{\textsf{Aim:}} Control the $L^p$-bounds, $1<p<\infty$, of the first time derivative and second spatial derivatives of our solution $u$ in terms of the $L^p$ bounds of the input function $f$.
\medskip

Again applying standard PDE techniques, we obtain that for $t>0$ and $x\in\R^d$, the fundamental solution of \eqref{heat} is described by
\beq\label{heatfs}
U^0(x,t):=\frac{1}{4\pi}\,|t|^{-\frac{d}{2}}\,e^{-\frac{|x|^2}{t}}\,.
\eeq
From this we deduce that for suitable $f\in L^p(\R^{d+1}_{+})$, $1<p<\infty$, one has that the solution to \eqref{heat} is given by
$$u(x,t):=\int_{0}^{t}\int_{\R^d}U^0(x-y,\,t-s)\,f(y,s)\,dy\,ds\,.$$

Indeed, one can check this since for $k_{ij}=U^{0}_{x_i x_j}$ and $k'=U^{0}_{t}$
\beq\label{uder}
\eeq
 $u_{x_i x_j}(x,t)=\int_{0}^{t}\int_{\R^d}k_{ij}(x-y,\,t-s)\,f(y,s)\,dy\,ds\,,$
$$u_{t}(x,t)=f(x,t)+\int_{0}^{t}\int_{\R^d}k'(x-y,\,t-s)\,f(y,s)\,dy\,ds\,$$
hold both pointwise and in terms of distributions.

Extracting now the quintessence from the kernels $K:=k_{ij}$ or $K=k'$, we have the following properties:
\beq\label{QCZ}
\eeq
\begin{itemize}
\item $K$ obeys an anisotropic dilation symmetry, \textit{i.e.} if $\a>0$, $\d_\a(x,t)=(\a x_1,\ldots,\, \a x_d,\a^2 t)\:\Rightarrow\:K(\d_\a(x,t))=\a^{-d-2}\,K(x,t)\,;$
\item $K(x,t)=0$ for $t<0$;
\item $K$ is $C^{\infty}$ away from the origin or one can ask for less - for example a condition of the type
$$\sup_{|x|,\,|t|>0}\int_{\{(y,s)\,|\, |s|>2 |t|,\,|s|> |t|+|x|^2\}} |K(y-x, s-t)-K(y,s)|\,dy ds<\infty\,;$$
\item $\int_{\R^d} K(x,1)\,d\,x=0$;
\item $\sup_{y\in\R^d}\int_{\R^d} (1+|x|+\log \frac{|y|}{|<x,y>|})\,|K(x,1)|\,d\,x<\infty$.
\end{itemize}
One should notice the similarities between \eqref{CZ} and \eqref{QCZ}; exploiting these similarities, it turns out that one can combine the method of rotations with Calder\'on-Zygmund operator techniques in order to show that the operator defined as
\beq\label{QCZ1}
Tf(x,t):=K*f(x,t)\,,
\eeq
is bounded from $L^p(\R^{d+1})$ to $L^p(\R^{d+1})$ for any $1<p<\infty$.

Finally, combining the second item above with the trivial extension assumption that $f(x,t)=0$ for any $(x,t)\in\R^d\times\R_{-}$, one immediately notices that both $u_{x_i x_j}$ and $u_{t}$ in \eqref{uder} can be realized as convolution operators of the type \eqref{QCZ1}, and hence we conclude that
\beq\label{QCZ2}
\|u_t\|_{L^p(\R^{d+1}_{+})}+\|u_{x_ix_j}\|_{L^p(\R^{d+1}_{+})}\lesssim_{p} \|f\|_{L^p(\R^{d+1}_{+})},\:\:\:\:\:1<p<\infty\;.
\eeq

We end this subsection by mentioning that the systematic study of constant coefficient parabolic differential operators was initiated by F. Jones, \cite{Jon}, E. Fabes, \cite{Fabs}, and E. Fabes and M. Rivi\`ere, \cite{fr}.

\medskip
\noindent C. \emph{Connections between the theme of constant coefficient parabolic differential operators and that of the Hilbert transform along curves}.
\medskip

As we have already expressed earlier, the area of PDE was the source for many interesting problems that later developed into main themes of study within the field of harmonic analysis. To exemplify this, we turn our attention to the two directions described above:
\begin{itemize}
\item The study of constant coefficient elliptic differential operators was the starting point for the development of the Calder\'on-Zygmund theory, which in turn became the central pillar in the theory of singular integral operators - a classical, important branch of harmonic analysis.

\item In a parallel setting, the study of constant coefficient parabolic differential operators was the starting point for the theory of Hilbert transform and maximal operators along curves, which later connected naturally with the theory of Radon transforms and that of singular integral operators with anisotropic symmetries - another representative branch within harmonic analysis area.
\end{itemize}

While the first connection has been already clarified in our point $A$ above, to make transparent the second connection above we proceed as follows:

Assume we want to prove the $L^2$-boundedness of \eqref{QCZ1} with $K$ obeying \eqref{QCZ}. Then, in order to avoid the singularity at the origin, we first define a family of truncated kernels $\{K_{\ep, R}\}_{\{0<\ep<R<\infty\}}$ with
$K_{\ep, R}(x,t)=K(x,t)$ if $\ep<t<R$ and $0$ otherwise. With this, applying Plancherel and some standard reasoning, we see
that the $L^2$-boundedness of $Tf(x,t):=K*f(x,t)$ is in fact equivalent to the $L^{\infty}$-boundedness, uniformly in $0<\ep<R$, of the expression
\beq\label{khat}
\hat{K}_{\ep,R}(\xi,\eta)=\int_{\R^d}K(x,1)\,\int_{\ep}^{R} \frac{e^{i \xi s}\,e^{i x\cdot \eta s^{\frac{1}{2}} }}{s}\,ds\,dx\:.
\eeq

Based on our hypothesis on $K$, the uniform boundedness of $\hat{K}_{\ep,R}(\xi,\eta)$ in dimension $d=1$ is now essentially equivalent to the $L^2$-boundedness of $H_{\Gamma}$ along a parabola, $\g(x,t)=t^2$, since the corresponding multiplier for $H_{\Gamma}$ is given by
\beq\label{mpar}
m_{H_{\Gamma}}(\xi,\eta)=\int_{\R} \frac{e^{-i \xi t}\,e^{i \eta t^{2} }}{t}\,dt\:.
\eeq

This moral equivalence between the (uniform) $L^{\infty}$-boundedness of the expressions in \eqref{khat} and \eqref{mpar} is the fundamental observation that effects the successful transition from the realm of PDE toward that of harmonic analysis.

\subsubsection{\underline{Harmonic analysis takes over - singular oscillatory integral operators}}\label{sgo}

Now that we have gained a better perspective on the origin of theme I (which happens to be our main theme of interest in the present study), it is of no surprise that the first results concerning the Hilbert transform along curves were direct consequences of the work done in the study of constant coefficient parabolic differential operators.  Indeed, from the work in \cite{Fabs}, one can immediately deduce the $L^2(\R^2)$-boundedness of the Hilbert transform along $\Gamma=(t,t^\a)$ with $\a>0$ and $\a\not=1$ ($\g(x,y,t)=t^{\a}$). The proof of Fabes relies on complex integration methods.

In parallel and related with this result, E. Stein and S. Wainger initiated a systematic study of singular oscillatory integrals and associated operators. One of their first results, obtained in 1970 in \cite{swmul}, was the following:

If $\{a_j\}_{j=1}^{n}\subset\R_{+}$ and $\{b_j\}_{j=1}^{n}\subset\R$, $n\in\N$, one has
\beq\label{stwosc}
\left|\int_{\R} e^{i\,\sum_{j=1}^n b_j\,t^{a_j}}\,\frac{dt}{t}\right|<C(a_1,\ldots, a_n)\,,
\eeq
with $C$ \textit{independent} of $\{b_j\}_{j=1}^{n}$.

This result is based on Van der Corput estimates and has as a direct consequence the $L^2(\R^2)-$boundedness of $H_{\Gamma}$ for $\Gamma=(t,\gamma(t))$ with $\gamma(t)=\sum_{j=1}^n b_j\,t^{a_j}$, and it can be extended to more general classes called ``homogeneous curves".

An important advance was the passage from the $L^2$ case to the general $L^p$ case with $1<p<\infty$. Essentially relying on complex interpolation methods, this was achieved by Nagel, Rivi\`ere, and Wainger in \cite{NRW1} and \cite{NRW2}. The passage from the Hilbert transform $H_{\Gamma}$ toward the maximal operator $M_{\Gamma}$ was first realized in the special case $\Gamma=(t, t^2)$ by the same authors in \cite{NRW3}. The useful observation there was the fact that Fourier transform methods can be effective even when dealing with maximal (positive kernel) operators. In the same year, Stein, \cite{stemax1}, \cite{stemax2}, introduced a method relying on the so-called $g-$function technique, where here the $g-$function is related to the square-function introduced by Littlewood and Paley and with its modified continuous version known as the Luzin area integral. Relying on this approach, Stein proved general $L^p$ bounds for both the Hilbert transform and the maximal operator along homogeneous curves. A bit later, Stein and Wainger, \cite{sw1}, extended these results to more general classes of curves.

\subsubsection{\underline{Further connections}}

As we have already seen, by comparing \eqref{CZ} with \eqref{QCZ} and integrating this into our discussion from Section \ref{pdeorig}, point C above, the anisotropic dilation symmetry plays a fundamental role in the behavior of our convolution operators\footnote{Nowadays there is a well-established theory of so-called anisotropic H\"ormander classes of symbols both in homogeneous and inhomogenous forms - for more on this see \textit{e.g.} \cite{BB1} and \cite{BB2} and the bibliography therein.} that is connected to our main theme concerning the Hilbert transform (and maximal operators) along curves. Now, as explained in \cite{sw1}, anisotropic dilations appear naturally in several other related harmonic analysis problems, among which we mention:

- (a) \textsf{Estimates for suitable subelliptic partial differential operators}: anisotropic dilation structures feature prominently in the study of second order hypoelliptic operators of H\"ormander type, \cite{Hor67}, or in the related context given by the study of the inhomogeneous Cauchy-Riemann equation for domains in several complex variables, \cite{FoSt}. To exemplify, in the latter case a key role is played by the identification $\R^{2n+1}=\C^n\times\R=\{(z,t)\,|\,z\in \C^n,\,t\in\R\}$ with the latter visualized as the Heisenberg group with the standard group multiplication $(z,t)\cdot (z', t')=(z+z',t+t'+2 Im z\cdot \bar{z}')$. Notice that indeed in this situation the natural (\textit{i.e.} compatible with the structure of the Heisenberg group) class of dilations are anisotropic and are given by $\d_s(z,t)=(sz, s^2 t),$ $s>0$.

- (b) \textsf{Analysis on symmetric spaces}: here we can mention a class of problems that aim to extend Fatou's theorem - \textit{i.e.} the almost everywhere existence of boundary values of harmonic functions - to the setting of Lie groups/symmetric spaces (see \cite{st71}, \cite{Kor72}, \cite{st76}).

- (c) \textsf{Radon transforms}: the literature here is quite rich; the interested reader may wish to consult \cite{CNSW} for a detailed account and further bibliography. Another, more recent paper, with a different direction of investigation combining Radon singular integral expressions with maximal oscillatory behavior, is given by \cite{PP}.

\subsubsection{\underline{Zygmund's differentiation conjecture; other curved models}}
$\newline$

\noindent\textsf{The zero-curvature case}
$\newline$

This topic originates in Lebesgue's theory of integration.  In \cite{Leb}, he showed that for any (locally) integrable function over the real line and for almost every point, the value of the integrable function is the limit of infinitesimal averages taken about the point. Given Lebesgue's result, it is natural to ask about similar differentiability results in higher dimensions, say for functions on $\R^2$. However, as it turns out, this problem is significantly more subtle, in particular due to the existence of ``pathological" objects such as Besicovitch sets. Indeed, even the problem of defining an adequate sequence of averages around a point is far from trivial since the geometry of the sets over which we take the averages is critical for the well-posedness of this problem.

In light of these challenging aspects of higher-dimensional differentiation problem, an alternative line of inquiry is offered by studying the problem of differentiation for averages along (variable) one-dimensional sets (curves) in $\R^2$. The most representative example in this context is given by Zygmund's conjecture, which, informally, asks about differentiability of averages along families of lines whose directions are described by a Lipschitz vector field.  The formal statement is given by:

\begin{conj} \textbf{(Zygmund)} Assume one modifies the definition of the variable curve $\Gamma$ as given in \eqref{gam} by introducing an extra $y-$dependence, \textit{i.e.} $\g(x,y,t)=u(x,y)\, t$ where here $u:\,\R^2\,\rightarrow\,\R$ is a Lipschitz vector field. Then, taking $\ep_0$ small enough depending on $\|
u\|_{Lip}$ and defining the maximal operator
\beq\label{maxZ}
M_{\Gamma} f(x,y)=M_{u,\ep_0} f(x,y):=\sup_{0<\ep<\ep_0} \frac{1}{2\ep}\,\int_{-\ep}^{\ep}|f(x-t,y-u(x,y)t)|\,dt\,,
\eeq
we have that $M_{\Gamma}$ is bounded on $L^p(\R^2)$ for any $1<p<\infty$.
\end{conj}

\medskip
One can of course formulate a singular integral analogue of the above:
\medskip

\begin{conj} \textbf{(Stein)} With $\g$, $u$ and $\ep_0$ as above, let us define the Hilbert transform $H_{\G}$ along $\G$ as
\beq\label{HilbZ}
H_{\Gamma} f(x,y)=H_{u,\ep_0} f(x,y):=p.v.\,\int_{-\ep_0}^{\ep_0}f(x-t,y-u(x,y)t)\,\frac{dt}{t}\,.
\eeq
Then we have that $H_{\Gamma}$ is bounded on $L^p(\R^2)$ for any $1<p<\infty$.
\end{conj}

At this point, it is worth saying that the Lipschitz condition imposed on the vector field $u$ is in fact required; indeed, a counterexample based on a construction of Besicovitch-Kakeya sets shows that one cannot expect any $L^p$ bounds if $u$ is only assumed to be H\"older continuous of class $C^{\a}$ with any exponent $\a$ strictly smaller than one.

It is also important to notice the following: the entire discussion in the above statements focuses on the case $\g(x,y,t)=u(x,y)\,t$, that is, on the situation in which $\g$ has \emph{no curvature} in $t$ - hence, the ``\emph{flat case}".

Concerning exclusively the historical evolution of this ``flat case" we have the following:

The first major contribution was made by Bourgain in \cite{Bolip}, where he proved the $L^2$-boundedness of \eqref{maxZ} in the case of analytic vector fields $u$. (The $L^p$ case can also be proven with some standard modifications). Further insight for some particular cases of vector fields $u$ is offered in \cite{CSWWv} and in \cite{CNSW}, with the latter covering many other interesting situations for the non-zero curvature case as well.

The analogue of Bourgain's result for the Hilbert transform $H_{\Gamma}$ was proved by Stein and Street in \cite{SS1}.

A key insight into this subject was brought by Lacey and Li in \cite{Lali}. Indeed, departing from the important observation that the multiplier associated with $H_{\Gamma}$ has a modulation symmetry similar in nature to that of the Carleson operator, the authors introduced time-frequency methods in their study of $H_{\Gamma}$. Using
wave-packet techniques inspired by the influential work of Lacey and Thiele on the Bilinear Hilbert transform, \cite{lt1}, \cite{lt2}, Lacey and Li proved a conditional result: if a suitable Kakeya type maximal operator obeys some $L^2$ bounds, then, assuming $u$ is $C^{1+\ep}$, the corresponding  $H_{\Gamma}$ is $L^2$-bounded. This remains to date the best result in the realm of a genuine two-variable dependent vector field $u$.

In the last decade further model problems, with additional simplifying assumptions on the variable dependence of the underlying vector field have been considered: it is natural to first hope to better understand the situation $u(x,y)=u(x)$, \textit{i.e.} where the vector field $u$ depends only on a single variable. In this context, first Bateman, \cite{Bat}, - in the single annulus case - and then Bateman and Thiele, \cite{BT}, proved that $H_{\Gamma}$ is $L^p$-bounded for $p>\frac{3}{2}$. Similar results were later proved by Guo in \cite{Guoh2} and  \cite{Guohp} for the case in which $u$ is constant along a Lispchitz curve - \textit{i.e.} a Lipschitz perturbation of the situation treated in \cite{Bat} and \cite{BT}.

Notice finally, that all the progress made since Bourgain's result in \cite{Bolip} has exclusively addressed the Hilbert transform case rather than that of the maximal operator. This limitation is due to the fact that the techniques in more recent work to date have crucially exploited the linearity of the operator as well as its good properties under commutation with suitable Littlewood-Paley projections - both of which are absent in the case of the maximal operator.

$\newline$
\noindent\textsf{The non-zero curvature case}
$\newline$

As already described before, this case originates in the study of constant coefficient parabolic equations and is a further natural development of the work on singular integral operators presented in Section \ref{sgo}.

All of the existing results presented below required some sort of non-zero curvature in the parameter $t$ for a generic $\Gamma\equiv\Gamma_{x,y}=(t,\,-\g(x,y,t))$ relative to $H_{\G}$ and $M_{\G}$ described in \eqref{defHVT} - \eqref{defMVT} (with the extra $y$-dependence). As one may expect, the historical evolution of this topic slowly transitions from requiring smoothness in all of the three variables $x,\,y,\,t$, with successive generalizations over the class of curves $\g$ toward very recently imposing smoothness only in $t$ and requesting simply measurable dependence on $x$ and $y$.
\medskip

\noindent A. \textsf{Smooth $x,\,y$-dependence}. This case is treated in a vast array of papers, and thus we will not aim to be exhaustive but limit ourselves here only to those which are the most relevant for our purposes.

Following the earlier work of Stein, Wainger and Nagel described in Section \ref{sgo}, the same authors proved in \cite{NSWhc} the $L^2-$boundedness of $H_{\Gamma}$ and $M_{\Gamma}$ assuming that $\g$ is a globally $C^{\infty}$ curve and has uniform non-vanishing curvature in $t$.\footnote{\textit{I.e.}, the second-order derivative in $t$ is non-zero for every $x,\,y$.}  Next, in a series of papers \cite{NVWW1}, \cite{NVWW2}, \cite{NVWW3}, Nagel, Vance, Wainger and Weinberg, and in two others also joined by Cordoba, \cite{CNVWW1}, \cite{CNVWW2}, provided various necessary and sufficient conditions (in general dimensions) for the $L^p$ boundedness of $H_{\Gamma}$ and $M_{\Gamma}$ for the case of $\g$ smooth, depending only on $t$ and having non-zero curvature. Passing back to the setting of multivariable dependence of $\g$, in \cite{CWWhf}, the authors provided general $L^p$ bounds for both $H_{\Gamma}$ and $M_{\Gamma}$ in the case $\g(x,y,t)=x\,\g(t)$, where here $\g$ obeys some suitable nonvanishing curvature condition that nonetheless allows the case $\g(t)=e^{-\frac{1}{t^2}}$ - that is, vanishing of infinite order at the origin.\footnote{It is worth noticing here that the curve $\g(t)=e^{-\frac{1}{t^2}}$ is non-doubling, and thus it cannot belong to the class $\n\f$ or $\mathbf{NF}$. Accordingly, the result in \cite{CWWhf} is not covered by our Corollary  \ref{Crossprod}. On the other hand, Corollary  \ref{Crossprod} holds for any $\g(x,t)=u(x)\,\g(t)$ with $\g\in \n\f$ and $u$ only measurable. For more on this, please see the Final Remarks section.}

A very extensive and influential study appears in \cite{CNSW}. There, the authors proved $L^p$ bounds for both singular Radon transforms and their maximal analogues with the integration performed over general differential submanifolds of $\R^d$. A helpful insight into this more general theory had previously been offered by Christ in \cite{ch1}, where he investigated the behavior of the Hilbert transform along curves in the nilpotent setting. Further results were obtained as follows: in \cite{SWs}, the authors proved $L^p$ bounds under some suitable $t$-convexity and doubling hypothesis that are uniform in $x,y$; another direction, extending the work in \cite{CWWhf}, considers  the situation $\g(x,y,t)=P(x)\,\g(t)$ with $P$ a polynomial and $\g$ a suitable smooth convex curve with some non-zero curvature condition - for this, see \cite{Benh} and, more recently, \cite{CZhc} and \cite{LY2}.
\medskip

\noindent B. \textsf{Rough $x,\,y$-dependence}: $\g(x,y,t)=u(x,y)\,t^{\a}$ with $\a>0,\,\a\not=1$, and minimal smoothness assumptions on $u$.

A first result in this setting was obtained in \cite{MRi} where, in particular, it was proven that for generic $u:\,\R^2\,\rightarrow\,\R$ measurable one has that $M_{\Gamma}$ is $L^p$-bounded if and only if $p>2$. Further advancements on this topic have only recently been obtained: indeed, the first such result is obtained in \cite{GHLR} where the authors prove the $L^p-$boundedness with $1<p\leq 2$ of $M_{\Gamma}$ if $u$ is Lipschitz (any less smoothness is known to produce a counterexample). In the same paper, a treatment of the singular integral variant is also provided: preserving the Lipschitz assumption on $u$, one obtains the $L^p-$boundedness, $1<p<\infty$, for the single annulus case $H_{\Gamma}\circ P_k^2$ with $P_k^2$ the Littlewood-Paley projection in the second variable. Relying on this and on a further square function estimate, global $L^p$ bounds for $H_{\Gamma}$ are proved in \cite{DGTZ}.
\medskip

\noindent C. \textsf{Rough (minimal smoothness) $x$-dependence}: $\g(x,y,t)=\g(x,t)=u(x)\,\tilde{\g}(t)$ for $u$ measurable and $\tilde{\g}$ smooth and obeying a suitable non-zero curvature.

The first result in this direction is obtained in \cite{GHLR}. There, the authors prove the $L^p$-boundedness of $H_{\Gamma}$ and $M_{\Gamma}$ for $1<p<\infty$  and $\tilde{\g}(t)=t^{\a}$ for $\a>0,\,\a\not=1$. In \cite{LY2}, this result is extended to more general functions $\tilde{\g}$ but only for $H_{\Gamma}$. In an unpublished note from 2016, \cite{lvBLC}, that will soon be made available, we proved a similar result for $H_{\Gamma}$ in the $L^2$ case but for more general functions $\tilde{\g}$. Another, more recent extension, is provided in \cite{GRSY}.

\subsection{Historical background; motivation (II)}

\subsubsection{\underline{Maximal singular oscillatory integral operators}}

In the 1910's, Luzin, \cite{Luz}, formulated a foundational conjecture whose impressive history has deeply influenced the evolution of harmonic analysis in the last century: \emph{the Fourier series of any $f\in L^2(\TT)$ converges to $f$ almost everywhere}. Within a decade of the formulation of this conjecture, Kolmogorov - who was a Ph.D. student of Luzin's -  proved in \cite{Kol1}, \cite{Kol2} that there are functions $f\in L^1(\TT)$ whose Fourier series diverges (almost) everywhere. From this point on, it was widely believed that by possibly modifying Kolomogorov's counterexample one should be able to disprove Luzin's conjecture.  It took another almost fifty years until L. Carleson, \cite{c1}, surprised the math community by providing the positive answer to this conjecture.

Carleson's result turns out to be equivalent to the (weak) $L^2$-boundedness of the following maximal singular integral operator - called Carleson's operator:
\beq\label{carld}
Cf(x):=\sup_{a\in\R}\,\left|\int_{\R}f(x-y)\,\frac{e^{i\,a\,y}}{y}\,dy\right|\:.
\eeq
In 1969,  R. Hunt, \cite{hu}, proved that $C:\:L^p(\R)\,\rightarrow\,L^p(\R)$ for $1<p<\infty$ while
Sj\"olin, \cite{sj2}, extended this result to higher dimensions.

Motivated by 1) the work on the Hilbert transform along curves previously discussed in Section 1.4.2 and 2) the study of singular integrals on the Heisenberg group discussed in Section 1.4.3 (a) regarding estimates for suitable subelliptic partial differential operators, E. Stein proposed the following generalization of Carleson's result:
\medskip

\noindent \textbf{Conjecture (Polynomial Carleson operator, \cite{stk}, \cite{sw})} \textit{Let $\Q_{d,n}$ be the class of all real-coefficient polynomials in $d$ variables with no constant term and of degree less than or equal to $n$, and let $K$ be a suitable Calder\'on-Zygmund kernel on $\R^d$. Then the Polynomial Carleson operator defined as
\beq\label{polcarld}
C_{d,n}f(x):=\sup_{Q\in\Q_{d,n}}\left|
\,\int_{\R^d}e^{i\,Q(y)}\,K(y)\,f(x-y)\,dy\,\right|\:
\eeq
obeys, for any $1<p<\infty$, the bound
\beq\label{polcarldes}
\|C_{d,n} f\|_{L^p(\R^d)} \lesssim \|f\|_{L^p(\R^d)}\:.
\eeq}

In 2001, relying on Van der Corput estimates and $TT^{*}$-method, Stein and Wainger, \cite{sw}, verified the above conjecture in the \emph{non-zero curvature} case - that is, when the supremum in \eqref{polcarld} ranges only over polynomials having no linear term. Notice that this result thus does not extend Carleson's Theorem, which addresses precisely the zero-curvature case (\textit{i.e.} the linear term). In his PhD thesis, \cite{lv1}, developing an approach for treating the interaction of quadratic wave-packets and adapting to it the time-frequency analysis developed by C. Fefferman in \cite{f}, the author settled this conjecture in the affirmative for the case $d=1$, $n=2$, $p=2$. In 2011, \cite{lv3}, we completely solved the one dimensional. The latter is based on a local analysis developed around the concepts of mass and counting function that motivates a new discretization algorithm of the time-frequency plane and which has as a consequence the elimination of the so-called exceptional sets appearing in all of the previous approaches regarding the boundedness of the Carleson operator. This elimination of exceptional sets proved in turn quintessential to providing for the first time direct $L^2$ bounds for the Carleson operator - without using interpolation - thus answering an open question raised by Fefferman in \cite{f}. In Fall 2017, based on the methods developed by the author in \cite{lv1} and \cite{lv3}, Zorin-Kranich, \cite{zk1}, proved the higher dimensional case of this conjecture for $p\geq 2$ and general Calder\'on-Zygmund kernels that are not necessarily translation invariant. Shortly afterwards, in \cite{lv3n}, appealing to quite standard modifications of the one dimensional result, to which we added the Van der Corput estimates proved in \cite{sw} and utilized in \cite{zk1}, we provided the full range of $p$ for general $d$ in Stein's conjecture within the original class of translation invariant Calder\'on-Zygmuns kernels. Finally, in \cite{zk1}, Zorin-Kranich updated his initial argument using a different $L^p$-interpolation approach completing the case $1<p<2$ in the original version of \cite{zk1}.

We end this discussion on the Carleson operator by mentioning one of the most tantalizing open questions in the area of time-frequency analysis, which can be formulated at a heuristic level as follows:

\medskip
\noindent \textbf{Open question} [\textsf{Informal}]. \textit{What is the behavior of the almost everywhere pointwise convergence of Fourier series \emph{between} the two known cases for the Lebesgue-scale spaces $L^p(\TT)$, namely:  divergence for $p=1$ (Kolmogorov) and convergence for $p>1$ (Carleson--Hunt)?}
\medskip

To clarify the connections between this question and our discussion of Carleson operator, one can appeal to Stein's maximal principle, \cite{s1}, which allows one to recast questions of pointwise convergence in terms of weak-type bounds for maximal operators:

\medskip
\noindent \textbf{Open question} [\textsf{Formal}].  (1) \textit{Give a satisfactory description of the Lorentz spaces $Y\subseteq L^1(\TT)$ for which the Carleson operator $C$ obeys the condition
\beq\label{carlopl1}
\textrm{C maps}\:Y\:\textrm{boundedly to}\: L^{1, \infty}(\TT)\:.
\eeq
If such exists, describe the maximal Lorentz $Y$.}  (2) \textit{More generally, let $Y$ be a rearrangement-invariant (quasi-)Banach space. Provide necessary and sufficient conditions on $Y$ under which \eqref{carlopl1} holds.}
\smallskip

This question's relevance extends far beyond its specific formulation: a successful resolution of this question requires the development of significant new methods and ideas that lie at the interface between additive combinatorics and time-frequency analysis. For more on this and on the rich history of the origin and progress on this question, please consult the most recent results in \cite{lv9} and \cite{lvLac} and the bibliography therein.
\medskip

\subsubsection{\underline{Connections between Hilbert transform along curves and the }\\ \underline{Polynomial Carleson operator}}\label{HCa}
\medskip

Departing from the definition of our Hilbert transform along $\G\equiv (t,-\gamma(x,t))$ - see \eqref{defHVT} - we choose now a particular form of $\G$ by taking
\beq\label{gxt}
\g(x,t):=\sum_{j=1}^{n} a_j(x)\,t^j\,,
\eeq
with $\{a_j(\cdot)\}_j$ \emph{any} real measurable functions.

Next, let us notice the following:
\begin{itemize}
\item On the one hand, the $L^2-$boundedness of $H_{\Gamma}$ is equivalent via Parseval to the inequality
\beq\label{l2hgam}
\int_{\R^2}\left|\int_{\R}f(x-t,\eta)\frac{e^{i\,\eta\,\g(x,t)}}{t}\,dt\right|^2\,dx\,d\eta\lesssim \|f(x,\eta)\|_{L^2(\R^2)}^2\:.
\eeq
\item On the other hand, applying Kolmogorov's linearization procedure, one can rewrite \eqref{polcarldes} in the setting $d=1$ and $p=2$ as
\beq\label{l2carl}
\|C_{1,d}f\|_2^2:=\int_{\R}\left|\,\int_{\R}e^{i\,Q(x,t)}\,f(x-t)\,\frac{dt}{t}\,\right|^2\,dx\lesssim \|f\|_2^2\:,
\eeq
where here $Q(x,t):=\sum_{j=1}^{n} \tilde{a}_j(x)\,t^j$ with $\{\tilde{a}_j(\cdot)\}_j$ real measurable functions.
\end{itemize}

Since the bounds in both \eqref{l2hgam} and \eqref{l2carl} are independent of the choice of $\{a_j(\cdot)\}_j$ and $\{\tilde{a}_j(\cdot)\}_j$ respectively, one immediately notices that the $L^2$ boundedness of $H_{\Gamma}$ with $\Gamma$ derived from \eqref{gxt} is in fact \emph{equivalent} to the $L^2$ boundedness of the Polynomial Carleson operator in dimension $d=1$.

Finally, notice that if we modify \eqref{gxt} by taking the linear term to be zero, that is, we take $\g(x,t):=\sum_{j=2}^{n} a_j(x)\,t^j$ with $n\geq 2$ and $\{a_j\}$ as before, then the $L^2$ boundedness of $H_{\Gamma}$ becomes \emph{equivalent} to the one dimensional version of the result proved by Stein and Wainger in \cite{sw}.

As part of the refined analysis in the present paper, we will be able to show suitable $L^2$-decay relative to the size of the multiplier's phase of $H_{\Gamma}$. In view of the aforementioned equivalences, this will suffice in order to provide a new proof for the $L^p$-bounds, $1<p<\infty$, of the Polynomial Carleson operator with no linear term - see Corollary \ref{Polyncaseap} above.

\subsection{Historical background; motivation (III)}

\subsubsection{\underline{Bilinear Hilbert transform and maximal operators along curves}}\label{conBHT}
$\newline$
The original formulation of this third theme, as with those of the previous two, was
cast in terms of a single variable dependence, \emph{i.e.} for curves $\g(x,t)\equiv \g(t)$:
$\newline$

\noindent \underline{\textsf{General Problem (single variable dependence formulation).}} \textit{Let $\Gamma:=(t, -\g(t))$ be a plane curve with $\g$ a suitable (piecewise) smooth real function. \underline{Goal}: Understand the conditions on the curve $\Gamma$ under which one has that
\begin{itemize}
\item the bilinear Hilbert transform along the curve $\Gamma$ denoted by $H_{\Gamma}^{\B}$ and defined as
\beq\label{nhilb}
H_{\Gamma}^{\B}(f,g)(x):= \textrm{p.v.}\int_{\R} f(x-t)\,g(x+\g(t))\frac{dt}{t}\:,
\eeq
\item the (sub)bilinear maximal operator along the curve $\Gamma$ denoted by $M_{\Gamma}^{\B}$ and defined as
\beq\label{nmax}
M_{\Gamma}^{\B}(f,g)(x):= \sup_{\ep>0}\frac{1}{2\ep}\int_{-\ep}^{\ep}|f(x-t)\,g(x+\g(t))|\,dt\:,
\eeq
\end{itemize}
each map $L^{p}(\R)\times L^{q}(\R)\:\rightarrow\: L^{r}(\R)$ boundedly for some $p,\,q,\,r\geq 1$ with $\frac{1}{p}+\frac{1}{q}=\frac{1}{r}$.}
$\newline$

Early analogues of this problem have been studied in ergodic theory, particularly in relation with the fundamental problem of understanding the $L^p$-norm convergence of (non-)conventional bilinear averages. In the harmonic analysis setting we split our discussion of
historical evolution and motivation according to our present focus - the key concept of curvature:

$\newline$
\noindent\textsf{The zero-curvature/flat case}: $\g(t)=a\,t$ with $a\in\R\setminus\{-1,0\}$.
$\newline$

Historically this theme arose in connection with the study of the Cauchy transform along Lipschitz curves, \cite{Cal}, \cite{CMM}. Indeed, this study led Calder\'on to conjecture the $L^p$-boundedness of the \emph{Bilinear Hilbert transform} (BHT)
$H_{\Gamma_a}^{\B}$ with $\g(t)=a\,t$ and $a\in\R\setminus\{-1,0\}$ for H\"older exponents within the ``Banach triangle'' $p,\,q,\,r\geq 1$.

As explained in Section \ref{Dich}, in this situation $H_{\Gamma}$ obeys the modulation symmetry relation given by  \eqref{bihilbcsym}. This invites a method of proof relying on time-frequency analysis; using this key insight, and developing the ideas used by Carleson, \cite{c1}, and Fefferman, \cite{f},  M.\ Lacey and C.\ Thiele provided in \cite{lt1} and \cite{lt2} the affirmative resolution of Calder\'on's conjecture.\footnote{Strictly speaking, the saga of the original Calder\'on's conjecture is completed with the publication of \cite{GL04} in which the authors prove uniform bounds in the parameter $a$ for the Bilinear Hilbert transform $H_{\Gamma_a}^{\B}$. Also, there is still an open problem about the maximal range for $p,\,q,\,r$ that guarantees the boundedness of $H_{\Gamma_a}^{\B}$.} The analogous  result for the maximal operator \eqref{nmax} was proved by  M.\ Lacey in \cite{la3}.

$\newline$
\noindent\textsf{The nonzero-curvature/non-flat case}: $\g(t)=\sum_{j=2}^{n} a_j\,t^j$, with $n>1$.
$\newline$

This direction can be motivated in several ways: i) one route from a pure harmonic analysis perspective is to proceed via analogy with the non-flat case for the linear Hilbert transform; ii) another route, already alluded to at the beginning of this section, arises from ergodic theory and is based on the long-studied problem, \cite{Fu}, \cite{HKr}, of understanding the $L^p$-norm convergence of non-conventional bilinear averages, \textit{e.g.}\ $\frac{1}{N}\sum_{n=1}^N f_1(T^n) f_2(T^{n^2})$ for $T$ an invertible measure-preserving transformation of a finite measure space and $f_j\in L^{p_j}$. This can be interpreted as a \emph{discrete} version of our non-zero curvature direction here for suitable $\g$, although to date there is no satisfactory transference principle between the continuous and discrete cases in this situation; iii)  finally, yet another route, is offered by number theory, in relation with various non-linear extensions of Roth's theorem for sets of positive density - also known in the literature as Ergodic Roth Theorem(s) (see \textit{e.g.} \cite{Bou86}, \cite{Bou88}, \cite{DGR}, \cite{K19} and the bibliography therein).

Moving now to the concrete harmonic analysis setting offered by \eqref{nhilb}, the first such result was studied in \cite{li}, in the special case $\Gamma(t) = (t, t^d)$, $2 \leq d \in \mathbb N$. There, Li proved that $H_{\Gamma}^{\B}:\:L^{2}(\R)\times L^{2}(\R)\:\rightarrow \:L^{1}(\R)$ continuously by relying on the concept of $\sigma$-uniformity introduced in \cite{cltt} and inspired by  Gowers's work in \cite{gowers}.

In \cite{lv4}, \cite{lv10} the author proved boundedness in the maximal possible range of exponents, excluding potential end-points, for $H_{\Gamma}^{\B}$ with $\g$ belonging to a suitable class of curves $\n\f$ that includes in particular any Laurent polynomial with no term of degree $\pm 1$, as well as any finite linear combination of terms $|t|^{\a}\,|\log |t||^{\b}$ with $\a\not\in\{-1,0,1\}$.\footnote{As it turns out, the condition $\a\not=-1$ is not actually required.} Our results improved over \cite{li} both quantitatively and qualitatively, by providing for the first time a scale-type decay and by significantly extending the class of curves. The proof of our results combines elements of time-frequency analysis (Gabor frames) with orthogonality methods and relies on a subtle discretization procedure for the operator that simultaneously: (1) separates the variables on the frequency side and (2) preserves the high oscillation and smoothness in one variable of the multiplier's phase function.

Regarding the maximal operator analogue $M_{\Gamma}^{\B}$, in \cite{LX} the authors prove the expected H\"older range in the case $\g(t)=\sum_{j=2}^{n} a_j\,t^j$ with bounds that are uniform in the polynomial's coefficients. In a joint work with A. Gaitan, \cite{GL}, by extending the ideas in \cite{lv4}, \cite{lv10}, we prove the boundedness within the maximal range of exponents of $M_{\Gamma}^{\B}$ for $\g\in \n\f$. The work in \cite{GL} prepares us for the significantly more complex situation that will be treated in the second part of our study in \cite{lvUA3}. A key feature of both \cite{GL} and \cite{lvUA3} is that they bring to light a \emph{unitary} treatment of both the singular operator $H_{\Gamma}^{\B}$ and the maximal operator $M_{\Gamma}^{\B}$.
\medskip

\subsubsection{\underline{Connections between the bilinear and linear Hilbert transform along}\\ \underline{curves}}
\medskip

We focus our discussion by presenting an antithesis between the following objects:
\begin{itemize}
\item the bilinear Hilbert transform along the curve $\g(t)$:
\beq\label{defBVT1}
H_{\G}^{\B}(f,g)(x):= \textrm{p.v.}\int_{\R}f(x-t)\,g(x+\g(t))\,\frac{dt}{t}\,,
\eeq
which regarded as a multiplier becomes
\beq\label{defBVT2}
H_{\G}^{\B}(f,g)(x):= \textrm{p.v.}\int_{\R^2}\hat{f}(\xi)\,\hat{g}(\eta)\,\left( \int_{\R}e^{-i\,t\,\xi}\,e^{i\,\g(t)\,\eta}\,\frac{dt}{t}\right)\, e^{i x \xi}\,e^{i x \eta}\,d\xi\,d\eta\:.
\eeq
\item the Hilbert transform along the curve $\g(t)$:
\beq\label{defHVT1}
H_{\G}(f)(x,y):= \textrm{p.v.}\int_{\R}f(x-t,\,y+\g(t))\,\frac{dt}{t}\,,
\eeq
which regarded as a multiplier becomes
\beq\label{defHVT2}
H_{\G}(f)(x,y):= \textrm{p.v.}\int_{\R^2}\hat{f}(\xi,\eta)\,\left( \int_{\R}e^{-i\,t\,\xi}\,e^{i\,\g(t)\,\eta}\,\frac{dt}{t}\right)\, e^{i x \xi}\,e^{i y \eta}\,d\xi\,d\eta\:.
\eeq
\end{itemize}

By inspecting \eqref{defBVT2} and \eqref{defHVT2} one realizes that these two expressions have the \emph{same multiplier}.\footnote{The same holds if the curve $\g$ depends on both $x$ and $t$ instead of only on the $t$ variable.} Thus, it comes
as no surprise that, from the perspective offered by the analysis of the multiplier, one sees many similarities between the approaches of \eqref{defBVT1} and \eqref{defHVT1}, respectively. This is indeed the case when comparing our present paper with the corresponding \cite{lv4}, \cite{lv10}, or \cite{lvUA3}. However, while helpful in providing some intuition about the subtleties in the oscillatory behavior of the multiplier, these similarities dilute at the moment in which one's attention shifts from the multiplier towards the input objects(s). On the one hand, on the input side, \eqref{defBVT2} can be regarded as the tensor-product case of \eqref{defHVT2}, hinting that the treatment of the latter should be more difficult. On the other hand, when focusing on the $\xi, \eta$ variables, one notices that \eqref{defBVT2} is more singular than \eqref{defHVT2}, with the former playing the role of a diagonal projection of the latter. As a result of these two competing aspects, the \emph{resemblance} between the corresponding approaches for the two problems transitions from concrete - in the multiplier analysis - to merely philosophical as the proof moves its focus to the input objects.\footnote{A prime example of the latter situation - only philosophical resemblance - would be provided in \cite{lvUA3} when antithetically discussing the $L^p$ approach - see also Section \ref{LPg2} in our present paper.}

\subsection{Structure of the paper}  In this final subsection of the Introduction we present the structure of our paper:

\begin{itemize}

\item In Section 2 we detail the definition of the newly introduced class of curves $\mathbf{M}_x\mathbf{NF}_{t}$.

\item  Several of the key notations used in this paper are introduced in Section 3.

\item  Sections 4 to 8 focus entirely on the proof of our Main Theorem, Part (I). Thus the most extensive part of the present study is dedicated to analysis of the behavior of $H_{\G}$ and $M_{\G}$.

    That being said, in Section 4 we remove the trivial components of the curve $\g$ and formulate our main task - see Theorem \ref{pmtha}. 

\item The analysis of the multiplier is discussed at length in Section 5; this analysis will be performed according to three cases: in Section 5.1 the low frequency case; in Section 5.2 the off-diagonal, non-stationary phase case; and finally in Section 5.3 the diagonal, stationary phase case. Based on this case-discussion we split accordingly each of our operators $\H_{\G}$ and $M_{\G}$ into three components. The first two components corresponding to the first two cases discussed above are solved in Theorems \ref{lowf} and \ref{hofdiag}. Of particular interest here is the unified treatment of the singular and maximal operator via the square function argument provided in Lemma \ref{translk}.

    The control over the main component of our operators resulting from the diagonal case is stated in Theorem \ref{Diagp}. Its proof covers the next two sections.

\item Section 6 treats the $L^2$-bound of the main piece $\L_{\G,j,m}$ appearing in the discretization of the diagonal term and is the central pillar in the construction of our paper. The key estimate is the exponential $m-$decay bound stated in Theorem \ref{l2dec}. The proof of this result is based on the earlier described LGC-methodology and involves an array of techniques that at each stage need to be compatible with one another: from discretization arguments having as an effect the phase linearization of the multiplier passing through Gabor frame decompositions, non-stationary phase, orthogonality and $T T^{*}-$arguments, time-frequency correlation, the implementation of the non-degeneracy condition \eqref{ndeg0} etc.

\item The $L^p$-bound, $1<p<\infty$, providing a unified treatment of both the maximal and the singular operators under discussion is provided in Section 7 and is the content of Theorem \ref{hmdiagp}. The key insight here is provided again by our central Lemma \ref{translk}, proved in Section 5, that applies immediately to the main desired estimate \eqref{apest}.

\item Section 8 provides several alternative approaches for the $L^p$ bounds of our operators, unraveling multiple other techniques, many of which rely on the author's work in \cite{lv10}, \cite{lvBLC}, \cite{lvUA3} and jointly \cite{GHLR}.

\item Section 9 treats Part (II) of our Main Theorem and reveals how various topics — treated until now separately in the harmonic analysis literature - can be brought under the same umbrella. Thus the techniques provided in this paper offer a universal treatment of topics such as linear and bilinear Hilbert and maximal operators along curves - for the latter this will be shown in the follow up study \cite{lvUA3} - as well as Polynomial Carleson type operators; in particular we obtain a new treatment of Stein and Wainger's result on the polynomial Carleson operator with no linear term.

\item In Section \ref{richclass} we present the proof of Theorem \ref{Gencurvpolyn}.  The central result of this section is represented by Proposition \ref{fewhe}. In order to prove the latter, we rely on Lemma \ref{coeffcontr} - a result
    that can be of independent interest.

\item The proofs of the remaining Corollaries \ref{Polyncaseap}, \ref{Crossprod} and \ref{Crossprod1} are provided in Section 11

\item Our paper ends with Section 12, in which various final remarks are presented.
\end{itemize}

$\newline$

\textbf{Acknowledgements:} I would like to thank \'Arp\'ad B\'enyi for providing me with several helpful comments and bibliographical materials. Also I would like to thank my Ph.D. student Alejandra Gaitan, for her patience and care in reading earlier drafts of the manuscript and correcting a number of typos therein. Finally, I'm grateful to my friend Zubin Gautam for elevating and improving the English presentation of this paper.

\section{Introducing the class of general curves $\mathbf{M}_x\mathbf{NF}_{t}$}\label{Curv}

In this section we intend to carefully define the class of curves introduced and utilized in the statements of our main results from Section \ref{Mres}. The reader is advised not to take too seriously the technical nature of the (extensive intended) definition but instead to picture as a point of reference our main model for the class $\mathbf{M}_x\mathbf{NF}_{t}$, that is $\g(x,t)=\sum_{k=2}^d a_k(x)\,t^k$ with $\{a_k(\cdot)\}_{k=2}^d$ measurable, $d\in\N$ with $d\geq 2$.

\begin{d0}\label{defgam} Given a function
\beq\label{mp}
\g:\:\R^2\,\rightarrow\,\R
\eeq
we say that $\g$ defines a \textbf{variable ($x$-measurable) family of twisted non-flat curves} and write
\beq\label{defnot}
\g\in \mathbf{M}_x\mathbf{NF}_{t}\,,
\eeq
if $\g$ is measurable on $\R^2$ such that\footnote{Below, the class $C^{2+}$ simply means the standard $C^{2+\d}$ where here  $\d$ can be any number strictly greater than zero. During our proof, for convenience and clarity, we will assume in fact that $\d=2$; however, this extra-assumption can be easily removed as long as we require $\d>0$. For more on this issue, please see the Final Remarks section.}
\beq\label{mpp}
\eeq
\begin{itemize}
\item $\g_t(\cdot):=\g(\cdot,\,t)$ is \emph{$x-$measurable} for every $t\in\R\setminus\{0\}$;

\item $\g_x(\cdot):=\g(x,\,\cdot)$ is \emph{$t-$smooth} within the class $C^{2+}(\R\setminus\{0\})$ for almost everywhere $x\in\R$\,,
\end{itemize}

and,\footnote{We can make this class translation invariant by adding - for free - $x$-measurable functions, that is, to allow combinations of the form $\g(x,t)+\mu(x)$ with $\mu$ measurable and $\g\in \mathbf{M}_x\mathbf{NF}_{t}$. However, one can reduce this case to the non-translation invariant case by a simple change of variable.}
iff there exist $0<c_0(\g)<c_1(\g)$, $N(\g)\in\N$, $A=A(\g)\in\N$ and $B=B(\g)\in\N$ such that the following conditions hold:

\begin{itemize}
\item \underline{\textit{$x$-fiber decomposition}}

One can partition
\beq\label{u}
\R=\bigcup_{\a=0}^{A}\R_{\a}\,,
\eeq
with each $\R_{\a}$ (possibly empty) being a Lebesgue measurable set and obeying

\begin{enumerate}
\item if $x\in\R_{0}$ then
\beq\label{mp2}
\g(x,t)\equiv 0\:\:\:\:\:\forall\:\:t\in\R\,;
\eeq

\item given $1\leq \a\leq A$, then for any $x\in\R_{\a}$ the map
$$\g_x:\;\R\,\rightarrow\,\R$$
 defined by $\g_x(t)=\g(x,t)$ has the following properties:
\end{enumerate}

\item \underline{\textit{smoothness, pointwise non-zero curvature, variation}}

One can partition\footnote{We allow the possibility for some of the components $\Z_{\b}^{x}$ to be empty.}
\beq\label{Zpart}
\Z=\bigcup_{\b=0}^{B} \Z_{\b}^{x}\,,
\eeq
with the following properties

- the set $\Z_{0}^{x}$ has bounded cardinality, \textit{i.e.}
\beq\label{contr0}
\# \Z_{0}^{x}<N(\g)\:.
\eeq

- each set $\{\Z_{\b}^{x}\}_{\{\b>0\}}$ is a convex set of integers;

- if we let $\Z_{\b}^{x}=:[a_{\b}^x, b_{\b}^x)\cap \Z$ and set $J_{\b}^{x}=(2^{-(b_{\b}^x+10)}, 2^{-(a_{\b}^x-10)})$ then, for any $1\leq \b\leq B$, one has that\footnote{Deduce immediately from \eqref{deriv0} that one can assume wlog that one also has $|\g'_x(\cdot)|>0\:\:\:\textrm{on}\:\:J_{\b}^{x}$.}
\beq\label{deriv0}
\g_x\in C^{2+}(J_{\b}^{x})\,,\:\:\:\textrm{with}\:\:\:|\g''_x(\cdot)|>0\:\:\:\textrm{on}\:\:J_{\b}^{x}\;.
\eeq
Moreover, setting $\Z^{x}=\bigcup_{\b=1}^{B} \Z_{\b}^{x}$, one has
\beq\label{variation0}
\sup_{x\in\R_{\a}}\sup_{c\in\R_{+}}\#\{j\in\Z^{x}\,|\,|2^{-j}\,\g'_x(2^{-j})|\in[c,2c]\}<N(\g)\:.
\eeq

\item \underline{\textit{doubling, uniform non-zero curvature (non-flatness)}}

Fix $1\leq \a\leq A$ and $x\in\R_{\a}$. Take now $1\leq \b\leq B$.

Then, for any $t\in I:=\{s\,|\,\frac{1}{10}\leq|s|\leq 10\}$ and $j\in \Z^{x}_{\b}$, we have
\beq\label{asymptotic0}
Q_{j}(x,t):=\frac{\g_x(2^{-j}\,t)}{2^{-j}\,\g'_x(2^{-j})}\in L^{\infty}_{x}C_t^{2+}(\R_{\a}\times I)\,.
\eeq
Moreover, one has the uniform bounds
\beq\label{fstterma0}
 \sup_{{t\in I}\atop{x\in\R_{\a}}} |Q_{j}(x,t)|<c_1(\g)\:,
\eeq
and\footnote{Throughout the paper, whenever we speak about expressions like $Q^{'}_{j}(x,t)$ we only refer to the $t-$derivative, that is $Q^{'}_{j}(x,t):=\frac{d}{dt}\,Q_{j}(x,t)$. This makes sense since we never assume in this paper any kind of smoothness in the $x-$variable but only measurability/boundedness.}
\beq\label{fstterm0}
 \inf_{{t\in I}\atop{x\in\R_{\a}}} |Q^{''}_{j}(x,t)|>c_0(\g)\:.
\eeq

\item \underline{\textit{non-degeneracy}}\footnote{For more on the significance of this condition one is invited to consult the Final Remarks section.}

Let $\phi$ be a positive Schwartz function supported in $\{\frac14\le |\xi|\le 4\}$ with $\sum_{n\in\Z} \phi(\xi/2^n)=1$ for all $\xi\not=0$ and set $\lfloor x\rfloor:=|x|+1$. Then, there exists $\bar{\ep}>0$ such that for any $m\in\N$ and $1\leq \a\leq A$ one has

\beq\label{ndeg0}
\sup_{j,\,k,\,n\in\Z} \int_{1<|s|<4}\,\sup_{t\in\R}\,\left(\frac{1}{2^{-j}}\int_{(k-\frac{1}{2})\,2^{-j}}^{(k+\frac{1}{2})2^{-j}}
\frac{\chi_{\R_\a}(x)\,\chi_{\Z^{x}}(j)\,\phi(\frac{\g'_x(2^{-j})}{2^n})}{\left\lfloor2^{\frac{m}{2}}
\left(\frac{\g'_x(2^{-j}(s+2^{j}\,x-k))}{2^n}-t\right)\right\rfloor^2}\,dx\right)\,ds\lesssim_{\g} 2^{-2\,\bar{\ep}\,m}\:.
\eeq
\end{itemize}

\end{d0}

\begin{o0}\label{rinv} From \eqref{asymptotic0}-\eqref{fstterm0} one immediately notices that for $x\in\R_{\a}$ and $j\in \Z_{x}^{\b}$, defining
\beq\label{asymptotic00}
q_{j}(x,t):=Q_{j}^{'}(x,t)=\frac{\g'_x(2^{-j}\,t)}{\g'_x(2^{-j})}\,,
\eeq
one has that there exist $0<c_0^{'}(\g)\leq c_1^{'}(\g)$ depending only on $c_0(\g)$ and $c_1(\g)$ such that for any $t\in I_0:=\{s\,|\,\frac{1}{5}\leq|s|\leq 5\}$ the following holds:
\beq\label{fstterm0q}
c_0^{'}(\g)<|q_{j}(x,t)|<c_1^{'}(\g)\:.
\eeq
Deduce thus that on $I_0$, there exists the inverse function of $q_{j}$, denoted with $r_{j}$, such that $r_{j}(x,t)\in L^{\infty}_{x}C_t^{1}$ in the corresponding domain of existence.
\end{o0}

\begin{o0}\label{W} The class $\mathbf{M}_x\mathbf{NF}_{t}$ contains any of the following\footnote{The first two items are straightforward while the next two are consequences of Theorem \ref{Gencurvpolyn}. The last item is very briefly discussed in the Final Remarks section.}:
\begin{itemize}
\item the set of all the real polynomial of degree $\geq 2$ with no constant and no linear term;

\item more generally, any element in $\n\f$ (for its definition, see \cite{lv4});

\item any $P(x,t)=\sum_{k=2}^d a_k(x)\,t^k$ with $\{a_k(\cdot)\}_{k=2}^d$ measurable where here $d\in\N$ with $d\geq 2$;

\item more generally, any $P(x,t)=\sum_{k=1}^d a_k(x)\,t^{\a_k}$ with $\{a_k(\cdot)\}_{k=2}^d$ measurable and $\a_k\in\mathbb{R}\setminus \{0,\,1\}$, where $d\in\N$.

\item even more so, any $P(x,t)=\sum_{k=1}^d a_k(x)\,t^{\a_k}\,\log^{\b_k} |t|$ with $\{a_k(\cdot)\}_{k=2}^d$ measurable, $\a_k\in\mathbb{R}\setminus \{0,\,1\}$ and $\b_k\in\mathbb{R}$, where $d\in\N$.
\end{itemize}
\end{o0}

\section{Notation}\label{Not}

For any smooth real function $\phi$ and $j\in\Z$ we set $\phi_{j}(\xi):=\phi(\frac{\xi}{2^j})$.

We next introduce the following convention: if $f$ is a function of two variables and $\phi$ is a single variable function, we write
\beq\label{con1}
(f{*}_{x}\phi)(x,y):=\int_{\R}f(x-s,y)\,\phi(s)\,ds\,,
\eeq
and similarly
\beq\label{con2}
(f{*}_{y}\phi)(x,y):=\int_{\R}f(x,y-s)\,\phi(s)\,ds\:.
\eeq

Let $M$ be the standard (one-dimensional) Hardy-Littlewood maximal operator defined as
\beq\label{Hmm}
Mf(x):=\sup_{{I\ni x}\atop{I\:\textrm{finite interval}}}\,\frac{1}{|I|}\,\int_{I}|f(s)|\,ds\:,
\eeq
where here $f\in L^1_{loc}(\R)$.

If $a\in \R$ is a given parameter, we set the ($a-$)shifted Hardy-Littlewood maximal operator as
\beq\label{Hmm1}
M^{(a)}f(x):=\sup_{{I\ni x}\atop{I\:\textrm{finite interval}}} \frac{1}{|I|}\,\int_{I+a\,|I|}|f(s)|\,ds\:,
\eeq
where as expected $I+a\,|I|:=\{x+a\,|I|\,|\,x\in I\}$. Also, we let
\beq\label{Hmm10}
M^{[a]}f:=M^{(-a)}f\,+\,M^{(a)}f\:.
\eeq

We set $M_{1}$ be the standard Hardy-Littlewood maximal function applied in the first variable and similarly $M_{2}$ be the standard Hardy-Littlewood maximal function applied in the second variable.

Throughout the paper, unless otherwise specified, the constant $C(\g)>0$ is a constant depending only on the properties of $\g$ that is allowed to change from line to line.

\section{Preparatives for the Main Theorem, Part (I)}\label{Prep}

We start the proof of our Main Theorem by focusing on Part (I) - the treatment of this first theme will cover the most consistent part of our paper.

With these said, in this section we perform few reductions in order to isolate the main term(s) for both  $H_{\Gamma}$ and $M_{\Gamma}$.

We first focus on the Hilbert transform $H_{\G}$. Observe that the Calder\'on-Zygmund kernel $\frac{1}{t}$ on $\R$ has two singularities: at zero and at infinity.  Using the dilation invariance of our kernel, we apply a Whitney type decomposition relative to our singularities
\beq\label{wit}
\frac1t = \sum_{j\in\Z} 2^{j}\rho(2^{j} t)\:,\:\:\:t\in\R\setminus\{0\}\,,
\eeq
where here $\rho$ is a smooth, odd, compactly supported function with its support in the set $\{t\in\R\,:\,\frac14<|t|<1\}$.

Set now
\beq\label{bdHVTpmm}
\rho_j(t):= 2^{j}\rho(2^{j} t)\,.
\eeq

With this we have
\beq\label{bdHVTpm}
H_{\G}(f)(x,y):= \sum_{j\in\Z}\int_{\R}f(x-t,\,y+\g(x,t))\,\rho_j(t)\,dt\:.
\eeq
In a similar fashion, letting
\beq\label{bdHVTpmm11}
\underline{\r}_j(t):=2^{j}\,\underline{\r}(2^{j} t):=2^j\,|\rho(2^{j} t)|\,,
\eeq
we deduce that
\beq\label{bdMVTpm}
M_{\G}(f)(x,y):= \sup_{j\in\Z}\int_{\R}|f(x-t,\,y+\g(x,t))|\,\underline{\r}_j(t)\,dt\:.
\eeq

Now our Main Theorem, part (I), follows from the following

\begin{t1}\label{pmth} With the above notations, for any $1<p<\infty$, one has
\beq\label{hm1}
\left\|H_{\G}(f)\right\|_{L^p(\R^2)} \lesssim_{\g,p} \left\|f\right\|_{L^p(\R^2)}\:,
 \eeq
 and
\beq\label{mhm1}
\left\|M_{\G}(f)\right\|_{L^p(\R^2)} \lesssim_{\g,p} \left\|f\right\|_{L^p(\R^2)}\:.
\eeq
\end{t1}

With this we appeal to Definition \ref{defgam}: from \eqref{u}, we know that
$$1=\sum_{\a=0}^{A}\chi_{\R_{\a}}(x)\,,$$
and thus
\beq\label{udec}
\left\|H_{\G}(f)\right\|_{L^p(\R^2)}^p=\sum_{\a=0}^{A}\,\left\|H_{\G}(f)\right\|_{L^p_x L^p_y(\R_{\a}\times\R)}^p\,,
\eeq
and
\beq\label{mudec}
\left\|M_{\G}(f)\right\|_{L^p(\R^2)}^p\leq\sum_{\a=0}^{A}\,\left\|M_{\G}(f)\right\|_{L^p_x L^p_y(\R_{\a}\times\R)}^p\,.
\eeq

From the definition of $\R_{0}$, it is trivial to notice that
$$\left\|H_{\G}(f)\right\|_{L^p_x L^p_y(\R_{0}\times\R)}^p=\int_{\R_{0,+}}\int_{\R}
\left| \textrm{p.v.}\int_{\R}f(x-t,\,y)\,
\frac{1}{t}\,dt\right|^p\,dy\,dx$$
$$\leq \|H_1f\|_p^p\lesssim_{p}\|f\|_p^p\,,$$
as a consequence of the $L^p-$boundedness of the (one dimensional) Hilbert transform. Similarly,
$$\left\|M_{\G}(f)\right\|_{L^p_x L^p_y(\R_{0}\times\R)}^p=\int_{\R_{0,+}}\int_{\R}
\left( \sup_{j\in\Z}\int_{\R}|f(x-t,\,y)|\,2^{j}\,\underline{\r}(2^{j} t)\,dt\right)^p\,dy\,dx$$
$$\leq\|M_1f\|_p^p\lesssim_{p}\|f\|_p^p\,,$$
as a consequence of the $L^p-$boundedness of the (one dimensional) Hardy-Littlewood maximal function.

Since $A$ is a fixed natural number depending only on $\g$ it is enough to focus on proving $L^p$-bounds for
$\chi_{\R_{\a}}(x)\,H_{\G}(f)(x,y)$ and $\chi_{\R_{\a}}(x)\,M_{\G}(f)(x,y)$ for a fixed $1\leq \a\leq A$. Thus from now on we can assume wlog that
$\a=1$ and that $\R_{\a}=\R$.

Define the following multipliers:
\beq\label{tfirstl}
\underline{m}_j(x,\xi,\eta):= \int_\R e^{-i\xi 2^{-j} t+i\eta \gamma_x(2^{-j}t)} \underline{\r}(t) dt\,,
\eeq
\beq\label{firstl}
m_j(x,\xi,\eta):= \int_\R e^{-i\xi 2^{-j} t+i\eta \gamma_x(2^{-j}t)} \rho(t) dt\,,
\eeq
and
\beq\label{firstlt}
m := \sum_{j\in\Z} m_j,
\eeq
or equivalently
\[\label{eqn:mult} m(x,\xi,\eta) := p.v. \int_\R e^{-i\xi t + i\eta \gamma_x(t)}\, \frac{1}{t}\,dt. \]

Deduce thus, that if regarded from the Fourier side, our operators are given by
$$H_{\G}=\int_{\R^2} e^{i\xi x+i\eta y} \widehat{f}(\xi,\eta) m(x,\xi,\eta) d(\xi,\eta)=\sum_{j\in\Z}H_{\G,j}\,,$$
with
\beq\label{hgj}
 H_{\G,j} f(x,y) = \int_{\R^2} e^{i\xi x+i\eta y} \widehat{f}(\xi,\eta) m_{j}(x,\xi,\eta) d(\xi,\eta)\,,
\eeq
and, assuming from now on wlog that $f\geq 0$, that
\beq\label{mgj}
M_{\G} f(x,y)=\sup_{j\in\Z}\L_{\G,j} f(x,y):= \sup_{j\in\Z}\left(\int_{\R^2} e^{i\xi x+i\eta y} \widehat{f}(\xi,\eta) \underline{m}_{j}(x,\xi,\eta) d(\xi,\eta)\right)\,.
\eeq

Now, recalling \eqref{Zpart}, we have that
$$H_{\G} f(x,y)=\sum_{\b=0}^{B} H_{\G}^{\b} f(x,y)\,,$$
$$M_{\G} f(x,y)\leq\sum_{\b=0}^{B} M_{\G}^{\b} f(x,y)\,,$$
where
$$H_{\G}^{\b}f(x,y):=\sum_{j\in\Z}\chi_{\Z_{\b}^{x}}(j)\, H_{\G,j}f(x,y)\,,$$
and
$$M_{\G}^{\b}f(x,y):=\sup_{j\in\Z}\, \L_{\G,j}(\chi_{\Z_{\b}^{x}}(j)f)(x,y)\,.$$

Let now $j(\cdot)$ be a real measurable function and set
\beq\label{meas}
\L_{\G,j(x)}f(x,y):=  \int_{\R^2} e^{i\xi x+i\eta y} \widehat{f}(\xi,\eta) \underline{m}_{j(x)}(x,\xi,\eta) d(\xi,\eta)\,.
\eeq

\begin{p1}\label{odem} With the above notations, for any $1<p<\infty$ and any real measurable function $j(\cdot)$, one has
\beq\label{h0}
\left\| \L_{\G,j(x)} f(x,y)\right\|_{L^p(\R^2)} \lesssim_{p} \left\|f\right\|_{L^p(\R^2)}\:.
 \eeq
\end{p1}
\begin{proof}

Fix $x\in\R$ and $1<p<\infty$. The key relation to prove is the following inequality:
\beq\label{keyh0}
\left\| \L_{\G,j(x)} f(x,\cdot)\right\|_{L_y^p(\R)} \leq M_1\left(\|f(\cdot,y)\|_{L_y^p(\R)}\right)(x)\:.
 \eeq
Indeed, if we assume for the moment this, then raising \eqref{keyh0} to the power $p$, integrating the result in the $x$ variable and using the standard Hardy-Littlewood maximal theorem we conclude the veracity of \eqref{h0}.

Returning now to  \eqref{keyh0}, by a simply application of Minkowski inequality we get

$$\left\| \L_{\G,j(x)} f(x,\cdot)\right\|_{L_y^p(\R)}
\leq \int_{\R}\left(\int_{\R}\left| f(x-t,\,y-\g(x,t))\right|^p\,dy\right)^{\frac{1}{p}}\,\,2^{j(x)}\,\underline{\r}(2^{j(x)}\,t)\,dt$$
$$\leq \sup_{j\in\Z}\int_{\R}\left(\int_{\R}\left| f(x-t,\,y)\right|^p\,dy\right)^{\frac{1}{p}}
\,2^{j}\,|\r(2^{j}\,t)|\,dt = M_1\left(\|f(\cdot,y)\|_{L_y^p(\R)}\right)(x)\:.$$
\end{proof}
Applying now Proposition \ref{odem} and appealing to \eqref{contr0}, we immediately deduce that
\begin{t1}\label{otrm}
With the previous notations, we have that
\beq\label{0trm}
\|H_{\G}^{0}f\|_{p},\,\|M_{\G}^{0}f\|_p\lesssim_{p} N(\g)\,\|f\|_p\:.
 \eeq
\end{t1}

We are thus left with controlling the $L^p$ bounds of $H_{\G}^{\b}$ and $M_{\G}^{\b}$  for $1\leq \b\leq B$. Since there are bounded many terms, we can assume wlog that $\b=1$.

For notational simplicity we set
$$\H_{\G}:=H_{\G}^{\b}\,,\:\:\M_{\G}:=M_{\G}^{\b}\:\:\:\textrm{and}\:\:\:\Z^{x}:=\Z_{\b}^{x}\,.$$

Thus, our Main Theorem follows now from

\begin{t1}\label{pmtha} For any $1<p<\infty$, one has
\beq\label{hm11}
\left\|\H_{\G}(f)\right\|_{L^p(\R^2)},\,\left\|\M_{\G}(f)\right\|_{L^p(\R^2)} \lesssim_{\g,p} \left\|f\right\|_{L^p(\R^2)}\:.
 \eeq
\end{t1}

Our entire work within the next four sections will focus on proving Theorem \ref{pmtha}.

\section{Analysis of the multiplier}\label{anmult}

In this section our goal is to isolate the main component(s) of our multiplier(s) that will have as a consequence the reduction of our operators $\H_{\G}$ and $\M_{\G}$ to their corresponding main terms.

Recalling now \eqref{firstl}, we set
\beq\label{mdef}
\mathfrak{m}_{j}(x,\xi,\eta):=\chi_{\Z^{x}}(j)\,m_{j}(x,\xi,\eta)=\chi_{\Z^{x}}(j)\,\int_\R e^{-i\xi 2^{-j} t+i\eta \gamma_x(2^{-j}t)} \rho(t) dt\,.
\eeq

We start our journey with performing a detailed analysis of the multiplier.

We first notice that we are dealing with a highly oscillatory integrand, and thus it is natural to expect an analysis of the phase appearing in \eqref{mdef} according to the principle of stationary phase.

Denoting the phase function by
\beq\label{firstl1}
\varphi_{x,\xi,\eta,j}(t):=-\frac{\xi}{2^j}\,t+\eta\,\g_x(\frac{t}{2^j})\:,
\eeq
and isolating its derivative
\beq\label{firstl2}
\frac{d}{dt}\,\varphi_{x,\xi,\eta,j}(t):=-\frac{\xi}{2^j}\,+\,\eta\,2^{-j}\,\g'_x(\frac{t}{2^j})\:,
\eeq
one observes, based on the properties obeyed by $\g$ - in particular that of $\g$ being doubling, that \eqref{firstl2} can be regarded at the heuristic level as
\beq\label{firstl3}
\frac{d}{dt}\,\varphi_{x,\xi,\eta,j}(t)\approx-\frac{\xi}{2^j}\,+\,\eta\,2^{-j}\,\g'_x(2^{-j})\:.
\eeq

Thus, in the view of \eqref{firstl3}, it becomes natural to apply a further decomposition relative to the
size of the terms involved in the phase derivative. Concretely, considering
\begin{eqnarray}\label{defph}
\:\:\:\:\:\:\:\:\:\:\phi\:\textrm{a positive} & \textrm{Schwartz function supported in}\:\{\frac14\le |\xi|\le 4\}\\
\nonumber &\textrm{with}\:\sum_{n\in\Z} \phi(\xi/2^n)=1\: \textrm{for all}\:\xi\not=0\,.
\end{eqnarray}
we use this partition of unity to write
\beq\label{firstl4}
1=\sum_{m,n,k\in\Z}\phi(\frac{\xi}{2^{m+j}})\,\phi(\frac{\g'_x(2^{-j})}{2^{n+j-k}})\,\phi(\frac{\eta}{2^k})\:,
\eeq

With these done, we notice that writing
\beq\label{firstl5}
 \mathfrak{m}_{j,m,n,k}(x,\xi,\eta):= \mathfrak{m}_{j}(x,\xi,\eta)\,\phi(\frac{\xi}{2^{m+j}})\,\phi(\frac{\g'_x(2^{-j})}{2^{n+j-k}})\,\phi(\frac{\eta}{2^k})\:,
\eeq
we have the following equality (in the distributional sense)
\beq\label{firstl6}
  \mathfrak{m}_{j}(x,\xi,\eta):=\sum_{m,n,k\in\Z}  \mathfrak{m}_{j,m,n,k}(x,\xi,\eta)\:.
\eeq
In some moments of our analysis it will be convenient to group the terms involved in the summation over $k$; for such situations it is advantageous to define
\beq\label{firstl7}
 \mathfrak{m}_{j,m,n}(x,\xi,\eta):=\sum_{k\in\Z}  \mathfrak{m}_{j,m,n,k}(x,\xi,\eta)\:.
\eeq
We notice that  $ \mathfrak{m}_{j,m,n}(x,\xi,\eta)$ can be written as
\beq\label{firstl9}
  \mathfrak{m}_{j,m,n}(x,\xi,\eta)=  \mathfrak{m}_j(x,\xi,\eta)\, \phi\left(\frac{\xi}{2^{m+j}}\right)\,\varrho\left(\eta,\,\frac{\g'_x(2^{-j})}{2^{n+j}}\right)\:,
\eeq
where above we set
\beq\label{defvr}
\varrho\left(\eta,\,\frac{\g'_x(2^{-j})}{2^{n+j}}\right):=
\sum_{k\in\Z}\,\phi(\frac{\g'_x(2^{-j})}{2^{n+j-k}})\,\phi(\frac{\eta}{2^k})\,,
\eeq
and further notice that $\varrho\in \C^{\infty}$ with $\textrm{supp}\,\varrho\subset \{(t,s)\,|\,\frac{1}{100}<|t\cdot s|<100\}$.

Following the ideas in \cite{lv4} and guided by the representation \eqref{firstl9}, we split our multiplier's analysis in three regions corresponding to the following situations\footnote{Throughout the paper $\Z_{-}:=\Z\setminus\N$ with $\N:=\{0,\,1,\,2\ldots\}$.}:
\begin{itemize}
\item  \textbf{(I)} \textit{the \underline{low frequency} case} - no oscillation present:
\beq\label{LF}
  \mathfrak{m}^{L}_j = \sum_{(m,n)\in (\Z_-)^2}  \mathfrak{m}_{j,m,n}\,;
\eeq

\item \textbf{(II)} \textit{the \underline{high frequency far from diagonal} case} - no stationary points present:

\beq\label{HFND}
  \mathfrak{m}^{H\not\Delta}_j = \sum_{(m,n)\in \Z^2\setminus ((\Z_-)^2 \cup \Delta)}  \mathfrak{m}_{j,m,n}\,,
\eeq

\item \textbf{(III)} \textit{the \underline{high frequency diagonal} case} - stationary points present:

\beq\label{HFND}
 \mathfrak{m}^{H\Delta}_j = \sum_{(m,n)\in\Delta}  \mathfrak{m}_{j,m,n}\,,
\eeq
\end{itemize}
where here $\Delta=\{(n,m)\in\Z^2\,:\,n,m\ge 0,\,|n-m|\le C(\gamma)\}$ with $C(\gamma)\ge 1$ a large constant depending only on $\gamma$.

With these, from \eqref{firstl6} - \eqref{HFND}, we have that
\beq\label{firstl8s}
  \mathfrak{m}_{j}=  \mathfrak{m}^{L}_j+ \mathfrak{m}^{H\not\Delta}_j+ \mathfrak{m}^{H\Delta}_j\:.
\eeq

Setting $\mathfrak{m}=\sum_{j}\mathfrak{m}_{j}$ and appealing to the obvious correspondences, we have
\beq\label{firstl8ss}
  \mathfrak{m}=  \mathfrak{m}^{L}+ \mathfrak{m}^{H\not\Delta}+ \mathfrak{m}^{H\Delta}\:,
\eeq
and
\beq\label{firstl8sss}
  \mathfrak{\underline{m}}_{j}=  \mathfrak{\underline{m}}^{L}_j+ \mathfrak{\underline{m}}^{H\not\Delta}_j+ \mathfrak{\underline{m}}^{H\Delta}_j\:,
\eeq
which ends our preliminary decomposition of the multiplier.

\begin{o0}\label{meanzer1}
The only difference between \eqref{firstl8s} and \eqref{firstl8sss} relies on the mean zero condition, that is
 $\mathfrak{m}_{j}(x,0,0)=0$ derived from $\int\rho=0$ as opposed to $\mathfrak{\underline{m}}_{j}(x,0,0)>0$. As we will soon see, \textit{this difference will be important \underline{only} at the level of the low-frequency term}.
\end{o0}

\subsection{(I) The low-frequency case.} In this section we will prove the following

\begin{t1}\label{lowf} Set
\beq\label{ml}
\M_{\G}^{L}(f)(x,y):=\sup_{j\in \Z} \left|\int_{\R^2} e^{i\xi x+i\eta y} \widehat{f}(\xi,\eta)\, \underline{m}_{j}^{L}(x,\xi,\eta) d(\xi,\eta)\right|
\eeq
and
\beq\label{hl}
\H_{\G}^{L}(f)(x,y):=\int_{\R^2} e^{i\xi x+i\eta y} \,\widehat{f}(\xi,\eta)\, m^{L}(x,\xi,\eta)\, d(\xi,\eta)\,.
\eeq
Then, the following holds:
\beq\label{ml1}
|\M_{\G}^{L}(f)(x,y)|\lesssim_{\g} M_1 M_2 f(x,y)\:,
\eeq
and
\beq\label{hl1}
|\H_{\G}^{L}(f)(x,y)|\lesssim_{\g} \left(\sum_{{j\in\Z}\atop{k\in\Z}}|\phi(\frac{\g'_x(2^{-j})}{2^{j-k-1}})\,M_1(f*_{y}\check{\phi}_{k})(x,y)|^2\right)^{\frac{1}{2}}\:.
\eeq
This further implies that, for any $1<p<\infty$, one has
\beq\label{mhl1p}
\|\M_{\G}^{L}(f)\|_{p},\,\|\H_{\G}^{L}(f)\|_{p}\lesssim_{\g,p} \|f\|_p\:.
\eeq
\end{t1}

The proof of the theorem above will be given in several steps below.

\subsubsection{Decomposing the low-frequency multiplier(s) into elementary building blocks}

As mentioned in the multiplier itemization above, see \eqref{LF}, in the low-frequency situation the phase has essentially no oscillation. Consequently, the main role will be played by the properties of the $t-$integrant in the \textit{absence} of the complex exponential (phase) - this last step will be rigorously justified via a Taylor series argument.
$\newline$

\noindent\textbf{1. The multiplier $\mathfrak{\underline{m}}^{L}_j$.}
$\newline$

In this setting, recalling \eqref{tfirstl}, \eqref{LF} and \eqref{firstl8sss}, we notice that
\beq\label{sm1}
\left|\frac{\xi}{2^j}\right|,\,\left|\frac{\g'_x(2^{-j})\,\eta}{2^{j}}\right|\lesssim 1\:.
\eeq
Deduce that in this regime the equality below is well defined
\beq\label{mhl}
\eeq
$$\int_\R  e^{-i\xi 2^{-j}t + i\eta \gamma_x(2^{-j}t)} \underline{\rho}(t) dt$$
$$= \sum_{p,\,\ell\in\N} \frac{i^{p+\ell}(-1)^\ell}{p!\ell !} \left( \frac{\eta \gamma_x'(2^{-j})}{2^j}\right)^p \left(\frac{\xi}{2^j}\right)^\ell \int_\R \left( \frac{\gamma_x(2^{-j}t)}{2^{-j}\gamma'_x(2^{-j})}\right)^p t^\ell \underline{\rho}(t) dt\:.$$

Denote with
\beq\label{constm}
\underline{C}_{p,\ell,j,x}:=\int_\R \left( \frac{\gamma_x(2^{-j}t)}{2^{-j}\gamma'_x(2^{-j})}\right)^p t^\ell \underline{\rho}(t) dt\:,
\eeq
and notice that based on the hypothesis imposed on $\g$, we have that
\beq\label{curvegM}
\|\underline{C}_{p,\ell,j,x}\|_{L_x^{\infty}(\R)}\lesssim_{p,\ell,\g} 1\:\:\:\:\:\forall\:p,\,\ell\in\N,\,j\in\Z\:.
\eeq

Setting $\tilde{\phi}_p(\xi):=\xi^p\,\phi(\xi)$ and the standard $\tilde{\phi}_{p,j}(\xi):=\phi_p(\frac{\xi}{2^j})$, one concludes that
\beq\label{mhl}
\mathfrak{\underline{m}}^{L}_j(x,\xi,\eta)
=  \chi_{\Z^{x}}(j)\,\sum_{{m,n\in\Z_{-}}\atop{k\in\Z}}\sum_{p,\,\ell\in\N} 2^{m l}\,2^{n p}\,\frac{i^{p+\ell}(-1)^\ell}{p!\ell !} \underline{C}_{p,\ell,j,x}\,\tilde{\phi}_{p}(\frac{\g'_x(2^{-j})}{2^{n+j-k}})\,\tilde{\phi}_{l,m+j}(\xi)\,\tilde{\phi}_{p,k}(\eta)\:.
\eeq

\noindent\textbf{2. The multiplier $\mathfrak{m}^{L}_j$.}
$\newline$

In this setting, recalling Observation \ref{meanzer1}, we appeal to \eqref{firstl} and make essential use - the only moment in the present paper - of the mean zero condition  $\hat{\rho}(0)=0$; consequently, in the same regime instituted by \eqref{sm1}, we have that as opposed to \eqref{mhl}, the $p=l=0$ term below is trivial:

\beq\label{inde}
\eeq
$$\int_\R  e^{-i\xi 2^{-j}t + i\eta \gamma_x(2^{-j}t)} \rho(t) dt$$
$$= \sum_{{p+\ell>0}\atop{p,\,\ell\in\N}} \frac{i^{p+\ell}(-1)^\ell}{p!\ell !} \left( \frac{\eta \gamma_x'(2^{-j})}{2^j}\right)^p \left(\frac{\xi}{2^j}\right)^\ell \int_\R \left( \frac{\gamma_x(2^{-j}t)}{2^{-j}\gamma'_x(2^{-j})}\right)^p t^\ell \rho(t) dt\:.$$
Using now the same notations as above with the obvious adaptation
\beq\label{consth}
C_{p,\ell,j,x}:=\int_\R \left( \frac{\gamma_x(2^{-j}t)}{2^{-j}\gamma'_x(2^{-j})}\right)^p t^\ell \rho(t) dt\:,
\eeq
one has
\beq\label{hhl}
\mathfrak{m}^{L}_j(x,\xi,\eta)
=  \chi_{\Z^{x}}(j)\,\sum_{{m,n\in\Z_{-}}\atop{k\in\Z}}\sum_{{p+\ell>0}\atop{p,\,\ell\in\N}} 2^{m l}\,2^{n p}\,\frac{i^{p+\ell}(-1)^\ell}{p!\ell !} C_{p,\ell,j,x}\,\tilde{\phi}_{p}(\frac{\g'_x(2^{-j})}{2^{n+j-k}})\,\tilde{\phi}_{l,m+j}(\xi)\,\tilde{\phi}_{p,k}(\eta)\:.
\eeq

\subsubsection{The maximal operator $\M_{\G}^{L}$ case}

Based on \eqref{mhl} and making use of \eqref{curvegM}, we deduce that
$$\M_{\G}^{L}(f)(x,y)\lesssim_{\g}$$
$$\sup_{j\in\Z^{x}}\left(\sum_{{m,n\in\Z_{-}}\atop{k\in\Z}}\sum_{p,\,\ell\in\N} \left|\frac{2^{m l}\,2^{n p}}{p!\ell !}\,
\underline{C}_{p,\ell,j,x}\, \,\tilde{\phi}_p(\frac{\g'_x(2^{-j})}{2^{n+j-k}})\,
(|f|*_{x}\check{\tilde{\phi}}_{l,m+j}*_{y}\check{\tilde{\phi}}_{p,k})(x,y)\right|\right)$$
$$\lesssim\sum_{m,n\in\Z_{-}}\sum_{p,\,\ell\in\N} \frac{2^{m l}\,2^{n p}}{p!\ell !}\, C^l\,C^p\,
 \sup_{j\in\Z}\left(\sum_{k\in\Z}\tilde{\phi}_p(\frac{\g'_x(2^{-j})}{2^{n+j-k}})\right)\,
M_1 M_2 f(x,y)$$
$$\lesssim_{\g} M_1 M_2 f(x,y)\,,$$
where above $C=C(\g)>0$ is some fixed positive constant depending only on $\g$ and for the last inequality we used property \eqref{variation0} of $\g$ to deduce that
\beq\label{jcgam}
\sup_{j\in\Z}\left(\sum_{k\in\Z}\tilde{\phi}_p(\frac{\g'_x(2^{-j})}{2^{n+j-k}})\right)\lesssim_{\g} C^p\,.
\eeq
This proves \eqref{ml1}.

\subsubsection{The Hilbert transform $\H_{\G}^{L}$ case}

From \eqref{hhl}, we deduce that
\beq\label{mlre}
\eeq
$$ \mathfrak{m}^{L}(x,\xi,\eta)=\sum_{j\in\Z}  \mathfrak{m}^{L}_j(x,\xi,\eta)=$$
$$\sum_{{p+\ell>0}\atop{p,\,\ell\in\N}} \sum_{j\in\Z} \chi_{\Z^{x}}(j)\sum_{(m,n)\in (\Z_-)^2} \sum_{k\in\Z} \frac{i^{p+\ell}(-1)^\ell}{p!\ell !}\,2^{m\,l}\,2^{n\,p}\,C_{p,\ell,j,x} \,\tilde{\phi}_l(\frac{\xi}{2^{m+j}})\,\tilde{\phi}_{p}(\frac{\g'_x(2^{-j})}{2^{n+j-k}})\,\tilde{\phi}_p(\frac{\eta}{2^k})\:.$$

Now due to the very fast decay of the coefficients in the above expansion, it becomes transparent that the main two terms to treat are those corresponding to the cases:
$\newline$

\noindent \textbf{Case 1} $l=0$, $p=1$,  $m\in \Z_{-}$ and $n=-1$

\noindent and

\noindent \textbf{Case 2} $l=1$, $p=0$,  $m=-1$ and $n\in \Z_{-}$

$\newline$
\noindent\textit{Treatment of Case 1}
$\newline$

\begin{p1}\label{ca1} Set $\psi(\xi):=\sum_{m\in \Z_-} \phi(\frac{\xi}{2^{m}})$ and with the above notations, define the multiplier
\beq\label{mul1}
 \mathfrak{m}^{L,1}(x,\xi,\eta):= \sum_{j\in\Z} \sum_{k\in\Z}\chi_{\Z^{x}}(j)\, C_{1,0,j,x} \,\psi(\frac{\xi}{2^{j}})\,\tilde{\phi}_{1}(\frac{\g'_x(2^{-j})}{2^{j-k-1}})\,\tilde{\phi}_1(\frac{\eta}{2^k})\,.
\eeq

Then, for any $1<p<\infty$, the operator
\beq\label{TL1}
 T_{ \mathfrak{m}^{L,1}} f(x,y) := \int_{\R^2} e^{i\xi x+i\eta y} \widehat{f}(\xi,\eta)  \mathfrak{m}^{L,1}(x,\xi,\eta) d(\xi,\eta),
\eeq
obeys the bound
\beq\label{TL1p}
 \|T_{ \mathfrak{m}^{L,1}} f\|_{p}\lesssim_{\g,p} \|f\|_p\,.
\eeq
\end{p1}

\begin{proof}
Let
\beq\label{f1}
\Lambda_{ \mathfrak{m}^{L,1}}(f,g):=<T_{ \mathfrak{m}^{L,1}} f,\,g>\,,
\eeq
where here $f\in L^p$ and $g\in L^{p'}$ with $\frac{1}{p}+\frac{1}{p'}=1$.

In light of the above, recalling our notation $\psi_j(\xi):=\psi(\frac{\xi}{2^j})$ and $\tilde{\phi}_{1,k}(\eta):=\tilde{\phi}_1(\frac{\eta}{2^k})$, we rewrite \eqref{f1} (ignoring conjugation) as\footnote{For notational simplicity, throughout the paper, when convenient we will think at, say, $\phi_k$ as equivalent with $\phi_k^2$ so that we can distribute the frequency location to both input functions via Parseval. Alternatively, one can appeal to the following standard modification:  letting $\vartheta\in C_0^{\infty}(\R)$ be one on the support of $\phi$ and zero outside the region $\frac{1}{20}<|\eta|<20$ one has the exact identity $\phi_k=\phi_k\,\vartheta_k$.}
\beq\label{ff1}
\eeq
$$\Lambda_{ \mathfrak{m}^{L,1}}(f,g)= $$
$$\sum_{{j\in\Z}\atop{k\in\Z}} \int_{\R^2}
\chi_{\Z^{x}}(j)\,C_{1,0,j,x}\,\tilde{\phi}_{1}(\frac{\g'_x(2^{-j})}{2^{j-k-1}})\,(f*_{x} \check{\psi}_j*_{y}\check{\phi}_{k})(x,y)\,(g*_{y}\check{\tilde{\phi}}_{1,k})(x,y)\,dx\,dy\:.$$

Thus, applying \eqref{curvegM}, Cauchy-Schwarz and H\"older in \eqref{ff1}, we deduce that
\beq\label{ff2}
\eeq
$$|\Lambda_{ \mathfrak{m}^{L,1}}(f,g)|\lesssim_{\g}$$
$$
\left\|\left(\sum_{{j\in\Z}\atop{k\in\Z}}|\tilde{\phi}_{1}(\frac{\g'_x(2^{-j})}{2^{j-k-1}})|\,|(f*_{x} \check{\psi}_j*_{y}\check{\phi}_{k})(x,y)|^2\right)^{\frac{1}{2}}\right\|_{p}\,\times\,$$
$$\left\|\left(\sum_{{j\in\Z}\atop{k\in\Z}}|\tilde{\phi}_{1}(\frac{\g'_x(2^{-j})}{2^{j-k-1}})|\,|(g*_{y}
\check{\tilde{\phi}}_{1,k})(x,y)|^2\right)^{\frac{1}{2}}\right\|_{p'}\:.$$
Now Proposition \ref{ca1} follows immediately from the two lemmas below.
\end{proof}

\begin{l1}\label{trm1lp} With the previous notations, for any $1<p<\infty$, one has
\beq\label{t1cont}
\left\|\left(\sum_{{j\in\Z}\atop{k\in\Z}}|\tilde{\phi}_{1}(\frac{\g'_x(2^{-j})}{2^{j-k-1}})|\,|(g*_{y}
\check{\tilde{\phi}}_{1,k})(x,y)|^2\right)^{\frac{1}{2}}\right\|_{p}\lesssim_{p,\g} \|g\|_{p}\:.
\eeq
\end{l1}
\begin{proof}
Recalling the definition/properties of $\g$ and that $\tilde{\phi}_{1}$ is compactly supported with $\tilde{\phi}_{1}(0)=0$ we immediately deduce that
\beq\label{gs}
\sum_{j\in\Z}|\tilde{\phi}_{1}(\frac{\g'_x(2^{-j})}{2^{j-k-1}})|\lesssim_{\g} 1\:.
\eeq
Thus, the LHS of \eqref{t1cont} is bounded from above by
$$\left\|\left(\sum_{k\in\Z}|(g*_{y}
\check{\tilde{\phi}}_{1,k})(x,y)|^2\right)^{\frac{1}{2}}\right\|_{p}\,,$$
which in turn, by standard Littlewood-Paley theory is bounded by $\|g\|_{p}$ proving our lemma.
\end{proof}

\begin{l1}\label{trm2lp} With the previous notations, for any $1<p<\infty$, one has
\beq\label{t2cont}
\left\|\left(\sum_{{j\in\Z}\atop{k\in\Z}}|\tilde{\phi}_{1}(\frac{\g'_x(2^{-j})}{2^{j-k-1}})|\,|(f*_{x} \check{\psi}_j*_{y}\check{\phi}_{k})(x,y)|^2\right)^{\frac{1}{2}}\right\|_{p}\lesssim_{p,\g} \|f\|_{p}\:.
\eeq
\end{l1}
\begin{proof}
The proof of this result is in the same spirit with the lemma above, with one modification. Indeed, recalling the notations in Section \ref{Not}, one has
 $$\left\|\left(\sum_{{j\in\Z}\atop{k\in\Z}}|\tilde{\phi}_{1}(\frac{\g'_x(2^{-j})}{2^{j-k-1}})|\,|(f*_{x} \check{\psi}_j*_{y}\check{\phi}_{k})(x,y)|^2\right)^{\frac{1}{2}}\right\|_p$$
 $$\lesssim\left\|\left(\sum_{{j\in\Z}\atop{k\in\Z}}|\tilde{\phi}_{1}(\frac{\g'_x(2^{-j})}{2^{j-k-1}})|\,|M_{1}(f*_{y}\check{\phi}_{k})(x,y)|^2\right)^{\frac{1}{2}}\right\|_{p}$$
which, after applying \eqref{gs}, is further dominated from above by
 $$\lesssim_{p,\g}
 \left\|\left(\sum_{k\in\Z}|M_1(f*_{y}\check{\phi}_{k})(x,y)|^2\right)^{\frac{1}{2}}\right\|_{p}$$
$$\lesssim_{p}
 \left\|\left(\sum_{k\in\Z}|(f*_{y}\check{\phi}_{k})(x,y)|^2\right)^{\frac{1}{2}}\right\|_{p}\lesssim\|f\|_{p}\,,$$
where for the second to the last relation we used Fefferman-Stein's inequality (\cite{FeSt71}).
\end{proof}

$\newline$
\noindent\textit{Treatment of Case 2}
$\newline$

\begin{p1}\label{ca2} Define the multiplier
\beq\label{mul2}
 \mathfrak{m}^{L,2}(x,\xi,\eta):= \sum_{j\in\Z} \sum_{k\in\Z}\chi_{\Z^{x}}(j)\,C_{0,1,j,x} \,\tilde{\phi}_1(\frac{\xi}{2^{j-1}})\,\psi(\frac{\g'_x(2^{-j})}{2^{j-k}})\,\phi(\frac{\eta}{2^k})\,.
\eeq

Then, for any $1<p<\infty$, the operator
\beq\label{TL2}
 T_{ \mathfrak{m}^{L,2}} f(x,y) := \int_{\R^2} e^{i\xi x+i\eta y} \widehat{f}(\xi,\eta) \mathfrak{m}^{L,2}(x,\xi,\eta) d(\xi,\eta),
\eeq
obeys the bound
\beq\label{TL2p}
 \|T_{ \mathfrak{m}^{L,2}} f\|_{p}\lesssim_{\g,p} \|f\|_p\,.
\eeq
\end{p1}

\begin{proof}
As in the proof of Proposition \ref{ca1}, we start by dualizing the problem, and write
\beq\label{f11f}
\Lambda_{ \mathfrak{m}^{L,2}}(f,g):=<T_{ \mathfrak{m}^{L,2}} f,\,g>\,.
\eeq
Further, we let $\psi_{\g,j,k}(x):=  C_{0,1,j,x}\,\psi(\frac{\g'_x(2^{-j})}{2^{j-k}})= C\,\psi(\frac{\g'_x(2^{-j})}{2^{j-k}})$.

Now, recalling Section \ref{Not}, we express \eqref{f11f} as
\beq\label{ff20}
\eeq
$$\Lambda_{ \mathfrak{m}^{L,2}}(f,g)= $$
$$\sum_{k\in\Z} \int_{\R^2}
\left(\sum_{j\in\Z}\chi_{\Z^{x}}(j)\,\psi_{\g,j,k}(x)\,(f*_{x} \check{\tilde{\phi}}_{1,j-1}*_{y}\check{\phi}_{k})(x,y)\right)\,(g*_{y}\check{\phi}_{k})(x,y)\,\,dx\,dy\:.$$

Proceeding as in the proof of Proposition \ref{ca1}, we apply Cauchy-Schwarz and H\"older in \eqref{ff20}, to deduce that
\beq\label{ff21}
\eeq
$$|\Lambda_{ \mathfrak{m}^{L,2}}(f,g)|\lesssim_{\g}\left\|\left(\sum_{k\in\Z}|(g*_{y}\check{\phi}_{k})(x,y)|^2\right)^{\frac{1}{2}}\right\|_{p'}\,\times\,$$
$$\left\|\left(\sum_{k\in\Z}\left|\sum_{j\in\Z}\chi_{\Z^{x}}(j)\,\psi_{\g,j,k}(x)\,(f*_{x} \check{\tilde{\phi}}_{1,j-1}*_{y}\check{\phi}_{k})(x,y)\right|^2\right)^{\frac{1}{2}}\right\|_{p}\:.$$
Since we know that $\phi(0)=0$, from standard Littlewood-Paley theory we deduce that
$$\left\|\left(\sum_{k\in\Z}|(g*_{y}\check{\phi}_{k})(x,y)|^2\right)^{\frac{1}{2}}\right\|_{p'}\lesssim_{p}\|g\|_{p'}\:.$$
Our proposition follows now from Lemma \ref{trm22p} below.
\end{proof}

\begin{l1}\label{trm22p} With the previous notations, for any $1<p<\infty$, one has\footnote{Recall that from Definition \ref{defgam}, one may assume wlog that $\Z^{x}=[a^x, b^x)$ with $a^x, b^x\in\Z$ measurable functions in the $x$ parameter.}
\beq\label{t3cont}
\left\|\left(\sum_{k\in\Z}\left|\sum_{j=a^x}^{b^x}\,\psi_{\g,j,k}(x)\,(f*_{x} \check{\tilde{\phi}}_{1,j-1}*_{y}\check{\phi}_{k})(x,y)\right|^2\right)^{\frac{1}{2}}\right\|_{p}\lesssim_{p,\g} \|f\|_{p}\:.
\eeq
\end{l1}
\begin{proof}
First we notice that by a standard density argument we may assume wlog that both sums in $j$ and $k$ are finite. Thus, we will assume that $k$ stays within the set $\{-N,\ldots,N\}$ for some $N\in\N$ and that $a^x,\, b^x$ are bounded measurable functions and prove that our estimates are independent of $N$, $a^x$ and $b^x$.

Define now $S(\tilde{\phi})_{j}:=\sum_{l=-\infty}^{j} \check{\tilde{\phi}}_{1,l-1}$. With this, applying Abel summation for the inner sum in \eqref{t3cont}, one has
$$\sum_{j=a^x}^{b^x}\psi_{\g,j,k}\, (f*_{x}\check{\tilde{\phi}}_{1,j-1}*_{y}\check{\phi}_{k})
=\sum_{j=a^x}^{b^x}\psi_{\g,j,k}\, (f*_{x}(S(\tilde{\phi})_{j}-S(\tilde{\phi})_{j-1})*_{y}\check{\phi}_{k})$$
$$= \psi_{\g,b^x,k}\, (f*_{x}S(\tilde{\phi})_{b^x}*_{y}\check{\phi}_{k})- \psi_{\g,a^x,k}\, (f*_{x}S(\tilde{\phi})_{a^x-1}*_{y}\check{\phi}_{k})$$
$$+\sum_{j=a^x}^{b^x-1}(\psi_{\g,j,k}-\psi_{\g,j+1,k})\, (f*_{x}S(\tilde{\phi})_{j}*_{y}\check{\phi}_{k})$$
For $N(x),\,M(x)$ positive integer measurable functions, define now
$$I(f):=\int_{\R^2}\left(\sum_{k=-N}^{N}\left|\psi_{\g,N(x),k}\, (f*_{x}S(\tilde{\phi})_{M(x)}*_{y}\check{\phi}_{k})\right|^2\right)^{\frac{p}{2}}\,dx\,dy\,,$$
and
$$II(f):=\int_{\R^2}\left(\sum_{k=-N}^{N}\left|\sum_{j=a^x}^{b^x-1}(\psi_{\g,j,k}-\psi_{\g,j+1,k})\, (f*_{x}S(\tilde{\phi})_{j}*_{y}\check{\phi}_{k})\right|^2\right)^{\frac{p}{2}}\,dx\,dy\:.$$
Notice now that \eqref{t3cont} follows from
\beq\label{p4cs}
I(f)+II(f)\lesssim_{\g,p} \|f\|_{p}^p\:.
\eeq
The first term is easy to treat. Indeed, we first notice that
\beq\label{sp1}
\|\psi_{\g,N(x),k}\|_{L^{\infty}_x}\lesssim_{\g} 1\:,
\eeq
and record the key condition $\phi(0)=0$. Thus, applying standard Littlewood-Paley and Calder\'on-Zygmund theory (including a variant of Cotlar's lemma for CZ operators), we have
$$I(f)\lesssim_{\g}\int_{\R^2}\left(\sum_{k=-N}^{N}\left| (f*_{x}S(\tilde{\phi})_{M(x)}*_{y}\check{\phi}_{k})\right|^2\right)^{\frac{p}{2}}\,dx\,dy$$
$$\lesssim\int_{\R^2}\left| f*_{x}S(\tilde{\phi})_{M(x)}\right|^p\,dx\,dy\lesssim \int_{\R^2}\left(\left| M_1(f*_{x}S(\tilde{\phi})_{\infty})\right|^p+ |M_1 f|^p\right)\,dx\,dy$$
$$\lesssim_{p}\|f\|_p^{p}\:.$$

Passing now to the second term, we have
$$II(f)\lesssim\int_{\R^2}\left(\sum_{k=-N}^{N}\left|\left(\sum_{j=a^x}^{b^x-1}|\psi_{\g,j,k}-\psi_{\g,j+1,k}|\right)\, \sup_j|f*_{x}S(\tilde{\phi})_{j}*_{y}\check{\phi}_{k}|\right|^2\right)^{\frac{p}{2}}\,dx\,dy$$
Next, we claim that
\beq\label{gamprop}
\left\|\sup_{k}\sum_{j=a^x}^{b^x-1}\left|\psi_{\g,j,k}-\psi_{\g,j+1,k}\right|\right\|_{L^{\infty}_x}\lesssim_{\g} 1\:.
\eeq
If we believe this for the moment, applying a similar argument with the one for $I(f)$ above at which we add Fefferman-Stein and Littlewood-Paley, we conclude
$$II(f)\lesssim_{\g}$$
$$\int_{\R^2}\left(\sum_{k=-N}^{N}\left| M_1(f*_{y}\check{\phi}_{k})\right|^2\right)^{\frac{p}{2}}\,dx\,dy\,+\,
\int_{\R^2}\left(\sum_{k=-N}^{N}\left| M_1(f*_{x}S(\tilde{\phi})_{\infty}*_{y}\check{\phi}_{k})\right|^2\right)^{\frac{p}{2}}\,dx\,dy$$
$$\lesssim_{p} \|f\|_p^p\:.$$

We are left now with proving \eqref{gamprop}. As expected, this will be based on the properties of $\g$, specifically
\eqref{deriv0} - \eqref{asymptotic0}. To see this, we use the fundamental theorem of calculus, in order to notice that
for $j\in \Z^{x}$ one has
\beq\label{ftc}
\eeq
$$\psi(\frac{\g'_x(2^{-j})}{2^{j-k-1}})-\psi(\frac{\g'_x(2^{-j-1})}{2^{j-k}})$$
$$=(\frac{\g'_x(2^{-j})}{2^{j-k-1}}-\frac{\g'_x(2^{-j-1})}{2^{j-k}})\,\int_{0}^{1} \psi'\left((1-t)\frac{\g'_x(2^{-j-1})}{2^{j-k}}+t\frac{\g'_x(2^{-j})}{2^{j-k-1}}\right)\,dt$$
$$=\frac{\g'_x(2^{-j-1})}{2^{j-k}}\,\left(\frac{2\,\g'_x(2^{-j})}{\g'_x(2^{-j-1})}-1\right)
\int_0^1\psi'\left(\frac{\g'_x(2^{-j-1})}{2^{j-k}}\,\left(1+t(\frac{2\,\g'_x(2^{-j})}{\g'_x(2^{-j-1})}-1)\right)\right)\,dt$$
$$=\frac{\g'_x(2^{-j})}{2^{j-k-1}}\,\left(1-\frac{\g'_x(2^{-j-1})}{2\,\g'_x(2^{-j})}\right)
\int_0^1\psi'\left(\frac{\g'_x(2^{-j})}{2^{j-k-1}}\,\left(1+(1-t)(\frac{\g'_x(2^{-j-1})}{2\,\g'_x(2^{-j})}-1)\right)\right)\,dt\:.$$
Adding now the fact that $|\psi'(t)|\lesssim \frac{1}{1+t^2}$, one concludes that
$$\sum_{j=a^x}^{b^x-1}|\psi_{\g,j,k}-\psi_{\g,j+1,k}|\lesssim_{\g} 1\,+\,
\sum_{j\in\Z}\frac{|\frac{\g'_x(2^{-j-1})}{2^{j-k}}|}{1+(\frac{\g'_x(2^{-j-1})}{2^{j-k}})^2} + \sum_{j\in\Z}\frac{|\frac{\g'_x(2^{-j})}{2^{j-k-1}}|}{1+(\frac{\g'_x(2^{-j})}{2^{j-k-1}})^2}\lesssim_{\g} 1$$
where in the last line we used that
$$\sup_{l,k\in\Z}\#\{j\in\Z\,|\, |\frac{\g'_x(2^{-j-1})}{2^{j-k}}|\in [2^{-l-1},\,2^{-l}]\}\lesssim_{\g} 1\,.$$
\end{proof}

We end this section with a reference to Observation \ref{meanzer1} that takes the following form:

\begin{o0}\label{meanzer} As mentioned, the mean zero condition of the function $\rho$ appearing in the definition of the multiplier $\mathfrak{m}_{j}$ is only used when dealing with the low frequency component $\mathfrak{m}_{j}^{L}$. Once this term is treated, one can completely dismiss this property. As a consequence, from now on, throughout the remaining part of the paper, we will identify the components $\mathfrak{m}^{H\not\Delta}_j\equiv \mathfrak{\underline{m}}^{H\not\Delta}_j$ and $\mathfrak{m}^{H\Delta}_j\equiv\mathfrak{\underline{m}}^{H\Delta}_j$ respectively. Alternatively, by decomposing $\r$ into the positive and negative part, one can consider from now on that $\r\geq 0$.
\end{o0}

\subsection{The off-diagonal, non-stationary phase case: $ \mathfrak{m}^{H\not\Delta}_{j}$}

As in the treatment of the Bilinear Hilbert transform along ``non-flat" curves, \cite{lv4}, the multiplier $\mathfrak{m}^{H\not\Delta}_j$ corresponding to the off-diagonal term deals with the situation when the phase of the $t-$integrant has no stationary points. As a consequence, we expect to see decay in the $m,\,n$ parameters.

Indeed, in our regime $(m,n)\in \Z^2\setminus ((\Z_-)^2 \cup \Delta)$, due to the lack of stationary points, we will be able to use a careful integration by parts which will have as a result the following evocative relation\footnote{This will be made precise during the proof of Theorem \ref{hofdiag} below.}
\beq\label{muldec}
 \mathfrak{m}_{j,m,n}=\frac{1}{2^{\max\{m,n\}}}\,\tilde{ \mathfrak{m}}_{j,m,n}\,,
\eeq
where here $\tilde{ \mathfrak{m}}_{j,m,n}$ is a multiplier having the same nature as $ \mathfrak{m}_{j,m,n}$.

Based on this observation, we will prove the following

\begin{t1}\label{hofdiag} Set
\beq\label{mhnd}
\M_{\G}^{H\not\Delta}(f)(x,y):=\sup_{j\in \Z} \left|\int_{\R^2} e^{i\xi x+i\eta y} \widehat{f}(\xi,\eta)\, \mathfrak{m}_{j}^{H\not\Delta}(x,\xi,\eta) d(\xi,\eta)\right|
\eeq
and
\beq\label{hnnd}
\H_{\G}^{H\not\Delta}(f)(x,y):=\int_{\R^2} e^{i\xi x+i\eta y} \,\widehat{f}(\xi,\eta)\, \mathfrak{m}^{H\not\Delta}(x,\xi,\eta)\, d(\xi,\eta)\,.
\eeq
Also, let
\beq\label{phig}
\phi_{\g,j,k,n}(x):=  \phi(\frac{\g'_x(2^{-j})}{2^{n+j-k}})\, \chi_{\Z^{x}(j)}\,.
\eeq
We then have
\beq\label{hmhndp}
\eeq
$$\|\M_{\G}^{H\not\Delta}(f)\|_p,\,\|\H_{\G}^{H\not\Delta}(f)\|_p\lesssim_{\g,p}\sum_{(m,n)\in \Z^2\setminus ((\Z_-)^2 \cup \Delta)} \frac{1}{2^{\max\{m,n\}}}\,\times$$
$$\left\|\left(\sum_{{j\in\Z}\atop{k\in\Z}} \phi_{\g,j,k,n}(x)\,\left|\left(f*_{x} \check{\phi}_{j+m}*_{y}\check{\phi}_{k}\right)\left(x-\frac{t}{2^j},y+\g_x(\frac{t}{2^j})\right)\right|^2 \right)^{\frac{1}{2}}\right\|_{L^{\infty}_t L^p_{x,y}([\frac{1}{10},10]\times\R^2)}$$
$$\lesssim_{\g,p} \|f\|_p\:.$$
\end{t1}

We start by making precise relation \eqref{muldec}.

As in \cite{lv4}, we define the differential operator
\beq\label{difop}
L:=\frac{-i}{\varphi'_{x,\xi,\eta,j}(t)}\,\partial_{t}=
\frac{-i}{-\frac{\xi}{2^j}+\frac{\eta}{2^j}\,\g'_x(\frac{t}{2^j})}\,\partial_{t}\:,
\eeq
in order to exploit its key feature
\beq\label{fixpoint}
L(e^{i\,\varphi_{x,\xi,\eta,j}(t)})=e^{i\, \varphi_{x,\xi,\eta,j}(t)}\:.
\eeq
Defining the adjoint of $L$ as
$$L^{\tau}(\r(t)):= \partial_{t}\left( \frac{i}{\varphi'_{x,\xi,\eta,j}(t)}\,\r(t)\right)\:,$$
we use integration by parts to deduce
\beq\label{partintegr}
\int_{\R} L(e^{-i\, \frac{\xi}{2^j}\, t}\: e^{i \,\eta \g_x(\frac{t}{2^j})})\,\r(t)\,dt=
\int_{\R} e^{-i\, \frac{\xi}{2^j}\, t}\: e^{i \,\eta \g_x(\frac{t}{2^j})}\,L^{\tau}(\r(t))\,dt\:.
\eeq

Define
\beq\label{Ajmn}
\eeq
$$A_{j,m,n}(x,\xi,\eta):=$$
$$\left(\int_{\R} e^{-i\, \frac{\xi}{2^j}\, t}\: e^{i \,\eta \g_x(\frac{t}{2^j})}\,\frac{i\,\r'(t)}{-\frac{\xi}{2^j}+\frac{\eta}{2^j}\,\g'_x(\frac{t}{2^j})}\,dt\right)\,
\phi(\frac{\xi}{2^{m+j}})\,\sum_{k\in\Z}\,\phi(\frac{\g'_x(2^{-j})}{2^{n+j-k}})\,\phi(\frac{\eta}{2^k})\:, $$
and respectively
\beq\label{Bjmn}
\eeq
$$B_{j,m,n}(x,\xi,\eta):=$$
$$\left(\int_{\R} e^{-i\, \frac{\xi}{2^j}\, t}\: e^{i \,\eta \g_x(\frac{t}{2^j})}\,\frac{-i\,\frac{\eta}{2^{2j}}\,\g''_x(\frac{t}{2^{j}})}{(-\frac{\xi}{2^j}+\frac{\eta}{2^j}\,\g'_x(\frac{t}{2^j}))^2}\,\r(t)\,dt\right)\,
\phi(\frac{\xi}{2^{m+j}})\,\sum_{k\in\Z}\,\phi(\frac{\g'_x(2^{-j})}{2^{n+j-k}})\,\phi(\frac{\eta}{2^k})\:.$$
With this, we have
$$ \mathfrak{m}_{j,m,n}(x,\xi,\eta)= \chi_{\Z^{x}}(j)\,A_{j,m,n}(x,\xi,\eta)\,+\,\chi_{\Z^{x}}(j)\,B_{j,m,n}(x,\xi,\eta)\;.$$

\noindent\textbf{Case 1}. $m>|n|+ C(\g)$.
$\newline$

Applying Taylor series, we have
\beq\label{T11}
\eeq
$$\frac{1}{-\frac{\xi}{2^j}+\frac{\eta}{2^j}\,\g'_x(\frac{t}{2^j})}=-\frac{1}{2^m}\,\frac{1}{\frac{\xi}{2^{j+m}}}\,\sum_{l=0}^{\infty} \frac{1}{2^{l(m-n)}}\,\left(\frac{\frac{\eta}{2^{n+j}}\,\g'_x(\frac{1}{2^j})}{\frac{\xi}{2^{j+m}}}\right)^l\,\left(\frac{\g'_x(\frac{t}{2^j})}{\g'_x(\frac{1}{2^j})}\right)^l\;,$$
and
\beq\label{T12}
\eeq
$$\frac{\frac{\eta}{2^{2j}}\,\g''_x(\frac{t}{2^{j}})}{(-\frac{\xi}{2^j}+\frac{\eta}{2^j}\,\g'_x(\frac{t}{2^j}))^2}
=\frac{1}{2^m}\,\frac{1}{\frac{\xi}{2^{j+m}}}\,\frac{\frac{1}{2^{2j}}\g''_x(\frac{t}{2^j})}{\frac{1}{2^j}\g'_x(\frac{t}{2^j})}\,\sum_{l=1}^{\infty} \frac{l}{2^{l(m-n)}}\,\left(\frac{\frac{\eta}{2^{n+j}}\,\g'_x(\frac{1}{2^j})}
{\frac{\xi}{2^{j+m}}}\right)^l\,\left(\frac{\g'_x(\frac{t}{2^j})}{\g'_x(\frac{1}{2^j})}\right)^l\;.$$
Deduce that
\beq\label{dec}
A_{j,m,n}=\sum_{l\in\N} A_{j,m,n,l}\:\:\:\textrm{and}\:\:\:B_{j,m,n}=\sum_{l\in\N^{*}} B_{j,m,n,l}\,,
\eeq
with
\beq\label{decm}
\eeq
$$A_{j,m,n,l}(x,\xi,\eta),\:B_{j,m,n,l}(x,\xi,\eta)\approx\frac{1}{2^m}\,\frac{1}{2^{l(m-n)}}\,\left(\int_{\R} e^{-i\, \frac{\xi}{2^j}\, t}\: e^{i \,\eta \g_x(\frac{t}{2^j})}\,\r_{j,l}(x,t)\,dt\right)\times$$
$$\tilde{\phi}_{-(l+1)}\left(\frac{\xi}{2^{m+j}}\right)\,
\sum_{k\in\Z}\,\phi_l(\frac{\g'_x(2^{-j})}{2^{n+j-k}})\,\phi_l(\frac{\eta}{2^k})\,$$
where, recalling \eqref{asymptotic0} and assuming $\chi_{\Z^{x}}(j)\not=0$, we have
\begin{itemize}
\item  $\r_{j,l}(x,t):=
\begin{cases}
-i\,\r'(t)\,(q_j(x,t))^l&\quad \textrm{for}\:\: A_{j,m,n,l}\,, \\
-i\,\r(t)\,q'_j(x,t)\,(q_j(x,t))^{l}&\quad\textrm{for}\:\:B_{j,m,n,l}\,.\\
\end{cases}$
\item from the first item and the properties of $\g$, we deduce that $\r_{j,l}(x,\cdot)$ is smooth and compactly supported in $t$ in the same region as the original function $\r(\cdot)$, and is an $L^{\infty}$-function in $x$ with $$\|\r_{j,l}(x,t)\|_{L_x^{\infty}C^{1}_t}\lesssim (C_{\g})^l\,,$$
  for some suitable $C_{\g}>0$ depending only on $\g$.
\item $\phi_l,\,\tilde{\phi}_l$ smooth, compactly supported away from the origin and $\|\phi_l\|_{C^{\a}},\,$
$\|\tilde{\phi}_l\|_{C^{\a}}\lesssim \a!\,|l|^{\a}\,|C|^l$ for some absolute constant $C\in\R$.
\end{itemize}
Inspecting \eqref{dec} and \eqref{decm} it is clear that in order to prove \eqref{hmhndp} for $\H_{\G}^{H\not\Delta}$
it is enough to show the $L^p$-boundedness of the operator whose multiplier is of the form $\sum_{j\in\Z} A_{j,m,n,0}$ for fixed $m,n$ with $m>|n|+ C(\g)$ with the formula for $A_{j,m,n,0}$ as described by \eqref{decm}.

In light of the above discussion, fixing $l\in\N$ and dropping for notational simplicity the sub-index $l$ (since we anyhow get a fast decay in the parameter $l$), we can assume wlog that the multiplier
$$ \mathfrak{m}_{m,n}(x,\xi,\eta)=\sum_{j\in\Z^{x}}\mathfrak{m}_{j,m,n}(x,\xi,\eta)\,,$$ is given by the expression
\beq\label{mfdmex}
\eeq
 $$\mathfrak{m}_{m,n}(x,\xi,\eta)= $$
$$\frac{1}{2^m}\,\sum_{{j\in\Z}\atop{k\in\Z}} \chi_{\Z^{x}}(j)\,\left(\int_{\R} e^{-i\, \frac{\xi}{2^j}\, t}\: e^{i \,\eta \g_x(\frac{t}{2^j})}\,\r_j(x,t)\,dt\right)\,
\phi\left(\frac{\xi}{2^{m+j}}\right)\,\phi(\frac{\g'_x(2^{-j})}{2^{n+j-k}})\,\phi(\frac{\eta}{2^k})\,.$$

Notice that $\mathfrak{m}^{H\not\Delta}(x,\xi,\eta)=\sum_{(m,n)\in \Z^2\setminus ((\Z_-)^2 \cup \Delta)} \mathfrak{m}_{m,n}(x,\xi,\eta)$.

With these, we have the following

\begin{p1}\label{lpcont} With the previous notations and conventions, let
\beq\label{tmn}
T_{ \mathfrak{m}_{m,n}} f(x,y)= T_{m,n}f(x,y):= \int_{\R^2} e^{i\xi x+i\eta y} \widehat{f}(\xi,\eta) \mathfrak{m}_{m,n}(x,\xi,\eta) d(\xi,\eta),
\eeq
where here $\mathfrak{m}_{m,n}$ is given by \eqref{mfdmex}.

Recalling \eqref{phig}, we have that for any $1<p<\infty$ the following holds:
\beq\label{Tmnpp}
\eeq
$$ \|T_{m,n} f\|_{p}\lesssim_{\g,p} $$
$$\frac{1}{2^{\max\{m,n\}}}\,
\left\|\left(\sum_{{j\in\Z}\atop{k\in\Z}} \phi_{\g,j,k,n}(x)\,\left|\left(f*_{x} \check{\phi}_{j+m}*_{y}\check{\phi}_{k}\right)\left(x-\frac{t}{2^j},y+\g_x(\frac{t}{2^j})\right)\right|^2 \right)^{\frac{1}{2}}\right\|_{L^{\infty}_t L^p_{x,y}([\frac{1}{10},10]\times\R^2)}$$
$$\lesssim_{\g,p}  m\,n_{+}\,2^{-\max \{m,n\}}\,\|f\|_p\,,$$
where here we set $n_{+}=\max\{n,\,1\}$.
\end{p1}

\begin{o0}\label{mhfdiag} Notice that relation \eqref{Tmnpp} immediately implies our Theorem \ref{hofdiag}. The proof of \eqref{Tmnpp} will be performed in two stages: the first inequality is detailed immediately below while the second one becomes the content of Lemma \ref{translk} afterwards.
\end{o0}

\begin{proof}

As before, we prefer to dualize our expression and thus define
$$\Lambda_{m,n}(f,g):= \int_{\R^2}T_{m,n} f(x,y)\,g(x,y)\,dx\,dy $$
$$=\frac{1}{2^m}\,\sum_{{j\in\Z}\atop{k\in\Z}}\,\int_{\R}\int_{\R^2}
\phi_{\g,j,k,n}(x)\,\left(f*_{x} \check{\phi}_{j+m}*_{y}\check{\phi}_{k}\right)(x-\frac{t}{2^j},y+\g_x(\frac{t}{2^j}))$$\,
$$\times(g*_{y}\check{\phi}_{k})(x,y)\,\r_j(x,t)\,dx\,dy\,dt\:,$$
where above we used the previous notations and conventions (see Section \ref{Not}).

Using again Cauchy-Schwarz and H\"older we have
\beq\label{cshp}
\eeq
$$|\Lambda_{m,n}(f,g)|\lesssim \frac{1}{2^m}\,\left(\int_{\R^2} \left(\sum_{{j\in\Z}\atop{k\in\Z}} \phi_{\g,j,k,n}(x)\,|(g*_{y}\check{\phi}_{k})(x,y)|^2\right)^{\frac{p'}{2}}\,dx\,dy\right)^{\frac{1}{p'}}\times $$
$$\int_{\R}\left(\int_{\R^2}\left(\sum_{{j\in\Z}\atop{k\in\Z}} \phi_{\g,j,k,n}(x)\,|\left(f*_{x} \check{\phi}_{j+m}*_{y}\check{\phi}_{k}\right)\left(x-\frac{t}{2^j},y+\g_x(\frac{t}{2^j})\right)|^2 |\r_j(x,t)|^2\right)^{\frac{p}{2}}\,dx\,dy\right)^{\frac{1}{p}} \,dt\:.$$
Now the first expression is easy to treat. Indeed, we only need to notice that $\|\sup_{k\in \Z}\sum_{j\in\Z} \phi_{\g,j,k,n}(x)\|_{L^{\infty}_x(\R)}\lesssim_{\g} 1$ and hence, using standard Littlewood-Paley, to conclude
$$\left(\int_{\R^2} \left(\sum_{{j\in\Z}\atop{k\in\Z}} \phi_{\g,j,k,n}(x)\,|(g*_{y}\check{\phi}_{k})(x,y)|^2\right)^{\frac{p'}{2}}\,dx\,dy\right)^{\frac{1}{p'}}$$
$$\lesssim_{\g}\left(\int_{\R^2} \left(\sum_{k\in\Z} |(g*_{y}\check{\phi}_{k})(x,y)|^2\right)^{\frac{p'}{2}}\,dx\,dy\right)^{\frac{1}{p'}}\lesssim\|g\|_{p'}^{p'}\:.$$
For the second term we first notice that $\|\sup_{j} \chi_{\Z^{x}}(j)\,|\r_j(x,t)|\|_{L^{\infty}_xL^1_{t}}\lesssim_{\g} 1$.
Once at this point, our proof follows by appealing to the lemma below.
\end{proof}

\begin{l1}\label{translk} Let $m,\,n\in\N$ and $1<p<\infty$. Then, for any $\frac{1}{10}\leq |t|\leq 10$, one has that
the following inequality holds uniformly in $t$:
\beq\label{trankey}
\left\|\left(\sum_{{j\in\Z}\atop{k\in\Z}} \phi_{\g,j,k,n}(x)\,\left|\left(f*_{x} \check{\phi}_{j+m}*_{y}\check{\phi}_{k}\right)(x-\frac{t}{2^j},y+\g_x(\frac{t}{2^j}))\right|^2 \right)^{\frac{1}{2}}\right\|_{p}\lesssim_{p,\g} m\,n_{+}\,\|f\|_{p}\:.
\eeq
\end{l1}
\begin{proof}
In what follows we will make use in an essential way of the result below appearing previously in various, slightly different forms in the math literature. The proof of the precise form displayed below, can be found in \cite{lv10} - see Lemma 3 therein:

\begin{l1}\label{Shiftedsq} (\cite{lv10}) Let $l\in \Z$. We define the $l-$shifted square function by
\beq\label{Sql}
S_l h(x)=\left(\sum_{j\in\Z} \left|(h*\check{\phi}_{j})(x-\frac{l}{2^j})\right|^2\right)^{\frac{1}{2}}\:,
\eeq
where here $h$ is a one-variable function and the meaning of $\phi_{j}$ is the same as before.
$\newline$
\indent Then one has
\beq\label{Sqp}
\|S_l h\|_p\lesssim_p \log (|l|+1)\,\|h\|_p\:.
\eeq
\end{l1}

In light of the above quoted result, we first provide a \textit{heuristic} for our claim \eqref{trankey}. Indeed one should think that

$$\left\|\left(\sum_{{j\in\Z}\atop{k\in\Z}} \phi_{\g,j,k,n}(x)\,\left|\left(f*_{x} \check{\phi}_{j+m}*_{y}\check{\phi}_{k}\right)\left(x-\frac{t}{2^j},y+\g_x(\frac{t}{2^j})\right)\right|^2 \right)^{\frac{1}{2}}\right\|_{p}$$
$$\approx \left\|\left(\sum_{{j\in\Z}\atop{k\in\Z}} |\left(f*_{x} \check{\phi}_{j+m}*_{y}\check{\phi}_{k}\right)(x-\frac{2^m}{2^{j+m}},y+\frac{2^n}{2^k})|^2 \right)^{\frac{1}{2}}\right\|_{p}$$
which, if we believe a vector valued inequality variant of \eqref{Sqp}, is further bounded from above by
$$n_{+}\,\left\|\left(\sum_{j\in\Z} |\left(f*_{x} \check{\phi}_{j+m}\right)(x-\frac{2^m}{2^{j+m}},y)|^2 \right)^{\frac{1}{2}}\right\|_{p}\lesssim_{p}n_{+}\,m\,\|f\|_{p}\:.$$

Now, in order to rigorously complete our proof, we must show that the following holds:
$\newline$

\noindent\textbf{Claim.} Let $1<p<\infty$ and $\{h_j\}_{j\in\Z}\subset L^{p}(\R)$. Then, with the above notations, one has uniformly for $\frac{1}{10}\leq |t|\leq 10$ and $x\in\R$ that
\beq\label{cla}
\eeq
$$\left\|\left(\sum_{{j\in\Z}\atop{k\in\Z}} \phi_{\g,j,k,n}(x)\,\left|\left(h_j*_{y}\check{\phi}_{k}\right)\left(y+\g_x(\frac{t}{2^j})\right)\right|^2 \right)^{\frac{1}{2}}\right\|_{L^p_y}\lesssim_{p,\g} n_{+}\,\left\|\left(\sum_{j\in\Z} |h_j(y)|^2 \right)^{\frac{1}{2}}\right\|_{L^p_y}\:. $$

In order to prove our claim we appeal to a standard vector valued result that appears in several forms in the math literature - and whose most convenient form for us is given by Theorem 1.1. in \cite{GLY}. Embracing for simplicity the notations from \cite{GLY}, we will now verify the hypothesis in the corresponding theorem.

Firstly, in our setting, the space of homogenous type $(\mathcal{H}, d, \mu)$ stands for $\R$ with the standard induced metric and Lebesgue measure. The Banach spaces $\B_1$, $\B_2$ correspond in our situation to the (isomorphic) Hilbert spaces
$l^2(\Z)$ and $l^2(\Z^2)$ respectively. Next, we consider a vector valued kernel $\vec{K}$ defined on $\R\times \R$ and such that $\vec{K}(y,s)$ is an element of $\mathcal{L}(l^2(\Z), l^2(\Z^2))$ - the space of all bounded linear operators from $l^2(\Z)$ to  $l^2(\Z^2)$ - and given in our context by
\beq\label{K}
\vec{K}(y,s)(\{a_j\}_{j\in\Z}):=\left\{a_j\,\phi_{\g,j,k,n}(x)\,\check{\phi}_{k}(y-s+\g_x(\frac{t}{2^j}))\right\}_{{j\in\Z}\atop{k\in\Z}}\,,
\eeq
where here $x,\,t$ are fixed real parameters with $\frac{1}{10}\leq |t|\leq 10$.

Given $y\in\R$, we set now
\beq\label{T}
\vec{T}(F)(y):=\int_{\R}\vec{K}(y,s)(F(s))\,d s
\eeq
and notice that this is a well defined element of $l^2(\Z^2)$ for any $F\in L^{\infty}(\R, l^2(\Z))$.

We state now the desired vector-valued result that we want to appeal to:
$\newline$

\noindent \textbf{Theorem (\cite{GLY}).} \textit{Assume that the operator $\vec{T}$ defined in \eqref{T} via \eqref{K} is a bounded linear operator from $L^r(\R,l^2(\Z))$ to $L^r(\R,l^2(\Z^2))$ for some $r\in (1,\infty)$ with norm $A_r>0$. Assume that $\vec{K}$ satisfies H\"ormander's condition
\beq\label{KHor}
\int_{|y-s|>2|z-s|}\|\vec{K}(y,s)-\vec{K}(y,z)\|_{l^2(\Z)\rightarrow l^2(\Z^2)}\,d y \leq C_{\vec{K}}\,,
\eeq
for some constant $C_{\vec{K}}>0$.
$\newline$
\noindent Then, for any $1<p<\infty$, $\vec{T}$ can be extended to an $L^p$-bounded operator, that is
\beq\label{Tlp}
\|\vec{T}(F)\|_{L^p(\R, l^2(\Z^2))}\leq C_p (C_{\vec{K}}+A_r)\,\|F\|_{L^p(\R, l^2(\Z))}\,,
\eeq
where here $C_p$ is a positive constant depending only on $p$.}

Deduce now that \eqref{cla} follows from the above theorem once we prove uniformly in $x, t$ that
\begin{itemize}
\item $\vec{T}$ is an $L^2$-bounded operator (thus $r=2$) with $A_r\approx 1$, that is
\beq\label{Tl2}
\|\vec{T}(F)\|_{L^2(\R, l^2(\Z^2))}\lesssim \|F\|_{L^2(\R, l^2(\Z))}\,.
\eeq
\item in \eqref{KHor} one has
\beq\label{ckb}
C_{\vec{K}}\lesssim_{\g} n_{+}\,.
\eeq
\end{itemize}
Now the first item is straightforward from Parseval since, for $F:=\{h_j\}_j\in L^2(\R, l^2(\Z))$, one has
$$\|\vec{T}(F)\|_{L^2(\R, l^2(\Z^2))}^2=\int_{\R}\sum_{{j\in\Z}\atop{k\in\Z}} \phi_{\g,j,k,n}(x)\,\left|\left(h_j*_{y}\check{\phi}_{k}\right)\left(y+\g_x(\frac{t}{2^j})\right)\right|^2\,dy$$
$$\lesssim\sum_{j\in\Z}\int_{\R}\sum_{k\in\Z} |\left(h_j*_{y}\check{\phi}_{k}\right)(y)|^2\,dy\lesssim \sum_{j\in\Z}\int_{\R}|h_j(y)|^2\,dy= \|F\|_{L^2(\R, l^2(\Z))}^2\:.$$
We pass now to the proof of \eqref{ckb}. In what follows we only discuss the case $n\in\N$ since if $n\in\Z_{-}$ the reasonings below become much easier.

Based on the translation invariant property of our kernel $\vec{K}$, it is enough to estimate the LHS of \eqref{KHor} when $s=0$:
\beq\label{KHor1}
\eeq
$$\int_{|y|>2|z|}\|\vec{K}(y,0)-\vec{K}(y,z)\|_{l^2(\Z)\rightarrow l^2(\Z^2)}\,d y $$
$$=\int_{|y|>2|z|}\left( \sup_{\|\{a_j\}\|_{l^2(\Z)}\leq 1}\sum_{{j\in\Z}\atop{k\in\Z}}|a_j|^2\,|\phi_{\g,j,k,n}(x)|^2\,
|\check{\phi}_{k}(y+\g_x(\frac{t}{2^j}))-\check{\phi}_{k}(y-z+\g_x(\frac{t}{2^j}))|^2\right)^{\frac{1}{2}}\,dy$$
$$\leq \int_{|y|>2|z|}\left( \sup_{j\in\Z} \sum_{k\in\Z}|\phi_{\g,j,k,n}(x)|\,
|\check{\phi}_{k}(y+\g_x(\frac{t}{2^j}))-\check{\phi}_{k}(y-z+\g_x(\frac{t}{2^j}))|\right)\,dy$$
$$=\int_{|y|>2|z|}\left( \sup_{j\in\Z} \sum_{k\in\Z}|\phi(\frac{2^k\,\g'_x(2^{-j})}{2^{n+j}})|\,\chi_{\Z^{x}(j)}\times\right.$$
$$\left.\left(2^{k}\,|\check{\phi}(2^k\,y+\frac{2^k\,\g'_x(2^{-j})}{2^{j}}\,Q_j(x,t))-\check{\phi}(2^k\,(y-z)+ \frac{2^k\,\g'_x(2^{-j})}{2^{j}}\,Q_j(x,t))|\right)\right)\,dy\,,$$

Next, recalling the $\g$-properties \eqref{asymptotic0} and \eqref{fstterma0}, we record the latter below
\beq\label{k0k1b}
|Q_j(x,t)|\leq c_1(\g)\,,
\eeq
where here $x,\,t$ are fixed real parameters with $|t|\approx 1$.

The last integral term can be further decomposed in
$$=\int_{|y|>2|z|}\left(\sup_{j\in\Z} \sum_{k>k_0}(\ldots)\right)\,+\,\int_{|y|>2|z|}\left(\sup_{j\in\Z} \sum_{k_1\leq k\leq k_0}(\ldots)\right)$$
$$+\,\int_{|y|>2|z|}\left(\sup_{j\in\Z} \sum_{k<k_1}(\ldots)\right)=:S_0+S_1+S_2\,,$$
where $k_1\leq k_0\in\Z$ depend solely on $z$, $n$ and $\g$, and are given by
\beq\label{k0k1}
2^{n+10}\,c_1(\g)<2^{k_0}\,|z|\leq  2^{n+11}\,c_1(\g)\:\:\:\textrm{and}\:\:\:k_1=k_0-n\:.
\eeq
For the first term we use the decay of $\check{\phi}$ to deduce
$$S_0\lesssim \int_{|y|>2|z|}\left(\sum_{k>k_0} \frac{2^k}{(2^k\,(|y|-|z|)-c_1(\g)\,2^{n+5})^4}\right)\,dy$$
$$\lesssim \sum_{k>k_0}\int_{|y|>2|z|} \frac{2^k dy}{|2^k y|^4}
\lesssim \sum_{k>k_0} \frac{1}{(2^k |z|)^3}\lesssim \frac{1}{(2^{k_0} |z|)^3}\lesssim_{\g}\frac{1}{2^{3n}}\:.$$

For the last term we apply the mean value theorem
$$S_2\leq\int_{|y|>2|z|}\left( \sup_{j\in\Z} \sum_{k\leq k_1}|\phi(\frac{2^k\,\g'_x(2^{-j})}{2^{n+j}})|\,\chi_{\Z^{x}(j)}\,\times\right.$$
$$\left.\left(2^{k}\,2^{k}\,|z|\,\int_{0}^1|\check{\phi}'(2^k\,y -s 2^k z +\frac{2^k\,\g'_x(2^{-j})}{2^{j}}\,Q_j(x,t))
|\,ds\right)\right)\,dy$$
$$\leq \sum_{k\leq k_1} \sum_{j\in\Z}\left( |\phi(\frac{2^k\,\g'_x(2^{-j})}{2^{n+j}})|\,\chi_{\Z^{x}(j)}\,\times\right.$$
$$\left.\left(2^{k}\,2^{k}\,|z|\,\int_{0}^1\int_{\R}|\check{\phi}'(2^k\,y -s 2^k z +\frac{2^k\,\g'_x(2^{-j})}{2^{j}}\,Q_j(x,t))
|\,dy\,ds\right)\right)$$
$$\lesssim_{\g}  \sum_{k\leq k_1} 2^{k}\,|z|\leq 2^{k_1}\,|z|\lesssim_{\eqref{k0k1}} 1\:.$$

Finally, for the middle term, we have
$$S_1\lesssim \int_{\R}\left( \sup_{j\in\Z} \sum_{k=k_1}^{k_0}|\phi(\frac{2^k\,\g'_x(2^{-j})}{2^{n+j}})|\,\chi_{\Z^{x}(j)}\,\times\right.$$
$$\left.\left(2^{k}\,|\check{\phi}(2^k\,y+\frac{2^k\,\g'_x(2^{-j})}{2^{j}}\,Q_j(x,t))|+|\check{\phi}(2^k\,(y-z)+ \frac{2^k\,\g'_x(2^{-j})}{2^{j}}\,Q_j(x,t))|\right)\right)\,dy$$
$$\lesssim \sum_{k=k_1}^{k_0}\int_{\R}\left( \sum_{j\in\Z} |\phi(\frac{2^k\,\g'_x(2^{-j})}{2^{n+j}})|\,\chi_{\Z^{x}(j)}\,
\left(2^{k}\,|\check{\phi}(2^k\,y+\frac{2^k\,\g'_x(2^{-j})}{2^{j}}\,Q_j(x,t))|\right)\right)\,dy$$
$$\lesssim_{\g} n\:.$$
With these our claim \eqref{cla} is verified.  This ends the proof of Lemma \ref{translk}.
\end{proof}

$\newline$
\noindent\textbf{Case 2}. $n>|m|+C(\g)$
$\newline$
In this situation our estimates for the terms defined in \eqref{Ajmn} and \eqref{Bjmn} change in the obvous fashion. Indeed,
applying Taylor series (again only for large values of $|j|$), one has
\beq\label{T21}
\eeq
$$\frac{1}{-\frac{\xi}{2^j}+\frac{\eta}{2^j}\,\g'_x(\frac{t}{2^j})}=
\frac{1}{2^n}\,\frac{1}{\frac{\eta}{2^{n+j}}\,\g'_x(\frac{1}{2^j})}\,\sum_{l=0}^{\infty} \frac{1}{2^{l(n-m)}}\,\left(\frac{\frac{\xi}{2^{j+m}}}{\frac{\eta}{2^{n+j}}\,\g'_x(\frac{1}{2^j})}\right)^l\,
\left(\frac{\g'_x(\frac{1}{2^j})}{\g'_x(\frac{t}{2^j})}\right)^{l+1}\;,$$
and
\beq\label{T22}
\eeq
$$\frac{\frac{\eta}{2^{2j}}\,\g''_x(\frac{t}{2^{j}})}{(-\frac{\xi}{2^j}+\frac{\eta}{2^j}\,\g'_x(\frac{t}{2^j}))^2}
=\frac{1}{2^n}\,\frac{1}{\frac{\xi}{2^{j+m}}}\frac{\frac{1}{2^{2j}}\g''_x(\frac{t}{2^j})}{\frac{1}{2^j}\g'_x(\frac{1}{2^j})}\,\sum_{l=1}^{\infty} \frac{l}{2^{(l-1)(n-m)}}\,\left(\frac{\frac{\xi}{2^{j+m}}}{\frac{\eta}{2^{n+j}}\,\g'_x(\frac{1}{2^j})}
\right)^l\,\left(\frac{\g'_x(\frac{1}{2^j})}{\g'_x(\frac{t}{2^j})}\right)^l\;.$$

From here on, one can apply similar methods with the one described at Case 1. We leave further details to the interested reader.

$\newline$

\subsection{The diagonal, stationary phase case: $ \mathfrak{m}^{H\Delta}_{j}$ - main term.}\label{statph}

This is the central case of our Main Theorem, part (I), due to the presence of the stationary points within the phase appearing in the definition of our multiplier $\mathfrak{m}^{H\Delta}_j$.

We will split our section in two subsections:
\begin{itemize}
\item the first one addresses the statements of the theorems treating our main terms.
Proving these theorems will be our focus for the next four sections.

\item the second one describes the analysis/properties of our main multiplier $\mathfrak{m}^{H\Delta}_j$.
\end{itemize}

\subsubsection{Stating the key results for the main terms}

Based on Theorems \ref{otrm}, \ref{lowf} and \ref{hofdiag}, we see that our Main Theorem, part (I), is a direct consequence of the following

\begin{t1}\label{Diagp} Let
\beq\label{mhnd1}
\M_{\G}^{H\Delta}(f)(x,y):=\sup_{j\in \Z} \left|\int_{\R^2} e^{i\xi x+i\eta y} \widehat{f}(\xi,\eta)\, \mathfrak{m}_{j}^{H\Delta}(x,\xi,\eta) d(\xi,\eta)\right|
\eeq
and
\beq\label{hnnd}
\H_{\G}^{H\Delta}(f)(x,y):=\int_{\R^2} e^{i\xi x+i\eta y} \,\widehat{f}(\xi,\eta)\, \mathfrak{m}^{H\Delta}(x,\xi,\eta)\, d(\xi,\eta)\,.
\eeq
Then, for any $1<p<\infty$, we have that
\beq\label{MHP}
\|\M_{\G}^{H\Delta}(f)\|_p,\,\|\H_{\G}^{H\Delta}(f)\|_p\lesssim_{\g,p} \|f\|_p\:.
\eeq
\end{t1}

Further on, our Theorem \ref{Diagp} will be a direct consequence of Theorems \ref{hmdiag2} and \ref{hmdiagp} below. Now, in order to be able to state these last results, we will need to introduce several notations.

We start by recalling the definition of $\mathfrak{m}^{H\Delta}_j$ in \eqref{HFND}. Since $\Delta=\{(n,m)\in\Z^2\,:\,n,m\ge 0,\,|n-m|\le C(\gamma)\}$, for notational simplicity, we will assume from now on wlog that
\beq\label{HFNDs}
 \mathfrak{m}^{H\Delta}_j = \sum_{m\in\N}  \mathfrak{m}_{j,m}\,,
\eeq
where here, recalling \eqref{firstl7}, we set for notational simplicity
\beq\label{mnm}
\mathfrak{m}_{j,m}:= \mathfrak{m}_{j,m,m}\,.
\eeq

Let now
\beq\label{lgjmk}
\L_{\G,j,k,m}f(x,y):= \int_{\R^2} e^{i\xi x+i\eta y} \widehat{f}(\xi,\eta) \mathfrak{m}_{j,m,m,k}(x,\xi,\eta) d(\xi,\eta)\,,
\eeq
and
\beq\label{lgjmn}
\L_{\G,j,m} f(x,y):= \sum_{k\in\Z} \L_{\G,j,k,m} f(x,y)\:.
\eeq

\begin{t1}\label{hmdiag2} [\textsf{The $L^2$-case}] Set
\beq\label{mhnm}
\M_{\G,m}(f)(x,y):=\sup_{j\in \Z} \left|\int_{\R^2} e^{i\xi x+i\eta y} \widehat{f}(\xi,\eta)\, \mathfrak{m}_{j,m}(x,\xi,\eta) d(\xi,\eta)\right|
\eeq
and
\beq\label{hnnd0}
\H_{\G,m}(f)(x,y):=\int_{\R^2} e^{i\xi x+i\eta y} \,\widehat{f}(\xi,\eta)\, \left(\sum_{j\in\Z} \mathfrak{m}_{j,m}(x,\xi,\eta)\right)\, d(\xi,\eta)\,.
\eeq

Recalling our hypothesis \eqref{ndeg0}, we then have that there exists $\ep=\ep(\bar{\ep})>0$ such that
\beq\label{HMm2}
\|\M_{\G,m}(f)\|_2,\,\|\H_{\G,m}(f)\|_2\lesssim_{\g}
\left\|\left(\sum_{{j\in\Z}\atop{k\in\Z}} \left|\L_{\G,j,k,m} f\right|^2 \right)^{\frac{1}{2}}\right\|_{2}\lesssim_{\g} 2^{-\ep\,m}\,\|f\|_2\:.
\eeq
\end{t1}

\begin{t1}\label{hmdiagp} [\textsf{The $L^p$-case}] With the previous notations, we have
\beq\label{HMmp}
\|\M_{\G,m}(f)\|_p,\,\|\H_{\G,m}(f)\|_p\lesssim_{\g,p}
\left\|\left(\sum_{{j\in\Z}\atop{k\in\Z}} \left|\L_{\G,j,k,m} f\right|^2 \right)^{\frac{1}{2}}\right\|_{p}\lesssim_{\g,p} m^2\,\|f\|_p\:.
\eeq
\end{t1}

\subsubsection{Analysis of $ \mathfrak{m}_{j,m}$; key properties}

Our approach will be based on a delicate analysis of our multiplier. The first part of our analysis goes in parallel with the corresponding analysis made in \cite{lv4}. For this reason we will only outline our approach and invite the interested reader to consult for more details the corresponding steps from \cite{lv4}.

Assume throughout the section that
\beq\label{mnz}
\mathfrak{m}_{j,m}(x,\xi,\eta)\not=0\:.
\eeq
After a careful analysis of the multiplier's phase
\beq\label{firstl1d}
\varphi_{x,\xi,\eta,j}(t):=-\frac{\xi}{2^j}\,t+\eta\,\g_x(\frac{t}{2^j})\:,
\eeq
based on the properties of $\g$ and following similar reasonings with the ones in \cite{lv4}, we claim that the following hold:
\begin{itemize}
\item  for $x,\,\xi,\,\eta,\,j\in\Z^{x}\,$ fixed, there exists exactly one critical point
\beq\label{tdefcritpoint}
t_{c,x}=t_c(x,\,\xi\,,\eta,\,j)\in J:=[2^{-k(\g)}, 2^{k(\g)}]
\eeq
with $k(\g)\in\N$ depending only on $\g$ such that
\beq\label{critpoint}
\varphi'_{x,\xi,\eta,j}(t_{c,x})=-\frac{\xi}{2^j}+\frac{\eta}{2^j}\,\g'_x(\frac{t_{c,x}}{2^j})=0\;.
\eeq

Also,  \beq\label{tc}
 t_{c,x}=2^{j}\,(\g'_x)^{-1}\left(\frac{\xi}{\eta}\right)\:\:\textrm{with}\:\: |\g'_x(2^{-j})|\approx_{\g}\left|\frac{\xi}{\eta}\right|\:;
 \eeq
Moreover, one can uniquely extend the functions $q_j(x,\cdot)$ and $r_j(x,\cdot)$ on the interval $[2^{-k(\g)}, 2^{k(\g)}]$. Then, one can rewrite \eqref{tc} as
\beq\label{tcinv}
q_j(x,\,t_{c,x})=\frac{\xi}{\g'_x(2^{-j})\,\eta}\:\:\:\Leftrightarrow\:\:\:r_j(x,\,\frac{\xi}{\g'_x(2^{-j})\,\eta})=t_{c,x}\:.  \eeq

\item  based on \eqref{asymptotic0}, we have that
 $$\varphi''_{x,\xi,\eta,j}(t)=\frac{\eta\,\g'_x(2^{-j})}{2^{j}}\times q'_j(x,t)\,,$$
and hence, using \eqref{fstterm0}, we further deduce that
\beq\label{hess}
|\varphi''_{x,\xi,\eta,j}(t)|\gtrsim_{\g} 2^m \:\:\:\:\:\:\;\;\:\:\:\:\:\:\forall\:(x,\,\xi,\,\eta)\in\textrm{supp}\,\mathfrak{m}_{j,m}\;.
\eeq
\item consider wlog $t>0$; also we let $\va,\,\tilde{\va}\in C_{0}^{\infty}(\R)$ such that $\textrm{supp}\,\va\subseteq[-10,10]$,  $\textrm{supp}\,\tilde{\va}\subseteq \{t\,|\,\frac{1}{100}<|t|<100\}$ and
$$1=\va(t)\,+\,\sum_{\k\in\N}\tilde{\va}(2^{-\k}\,t)\:\:\:\:\:\:\:\:\:\:\:\:\forall\:t\in\R\:.$$

Recalling \eqref{defvr}, we write now
\beq\label{mdecomp}
 \mathfrak{m}_{j,m}(x,\xi,\eta)= \chi_{\Z^{x}}(j)\,\A_{j,m}(x,\xi,\eta)\,+\,\sum_{\k=0}^{\infty}\chi_{\Z^{x}}(j)\,\B_{j,m}^\k(x,\xi,\eta)\:,
\eeq
where
\beq\label{Aterm}
\eeq
$$\A_{j,m}(x,\xi,\eta):= \,\left(\int_{\R} \va\left(2^{\frac{m}{2}}(s-\frac{\xi}{\g'_x(2^{-j})\,\eta})\right)\,e^{i\,\varphi_{x,\xi,\eta,j}( r_j(x,\,s))}\:\r(r_j(x,\,s))\,r'_j(x,\,s)\,ds\right)$$
$$\times\phi\left(\frac{\xi}{2^{m+j}}\right)\,
\varrho\left(\eta,\,\frac{\g'_x(2^{-j})}{2^{m+j}}\right)\:,$$
and
\beq\label{Bkterm}
\eeq
$$\B_{j,m}^\k(x,\xi,\eta):= \,\left(\int_{\R} \tilde{\va}\left(2^{\frac{m}{2}-\k}(s-\frac{\xi}{\g'_x(2^{-j})\,\eta})\right)\,e^{i\,\varphi_{x,\xi,\eta,j}( r_j(x,\,s))}\:\r(r_j(x,\,s))\,r'_j(x,\,s)\,ds\right)$$
$$\times\phi\left(\frac{\xi}{2^{m+j}}\right)\,
\varrho\left(\eta,\,\frac{\g'_x(2^{-j})}{2^{m+j}}\right)\:.$$

\item there exist the functions
$\{\zeta_{\k}\}_{\k\geq 0}$
with the properties
\beq\label{zeta}
\zeta_{\k}:\,[\frac{1}{10},10]\times [\frac{1}{10},10]\,\rightarrow\,\R\:\:\:\textrm{with}\:\:\:\;\:\:\|\zeta_{\k}(\cdot,\cdot)\|_{C^{1+}}\lesssim_{\g} 2^{-\k}
\eeq
such that
\beq\label{C}
\eeq
$$ \mathfrak{m}_{j,m,m}(x,\xi,\,\eta)= \chi_{\Z^{x}}(j)\,\times$$
$$\sum_{\k\geq 0} 2^{-\frac{m}{2}}\,e^{i\,\varphi_{x,\xi,\eta,j}(t_{c,x})}\,
\zeta_{\k}\left(\frac{\xi}{2^{m+j}},\,\frac{\eta\,\g'_x(2^{-j})}{2^{m+j}}\right)\,
\phi\left(\frac{\xi}{2^{m+j}}\right)\,\varrho\left(\eta,\,\frac{\g'_x(2^{-j})}{2^{m+j}}\right)\:,$$
\end{itemize}

\begin{o0}\label{redmult} 1) Strictly speaking, each of the above functions $\zeta_{\k}$ is in fact depending on $j,m$. However we have chosen not to write this explicit dependence, since the norm $\|\zeta_{\k}\|_{C^{1+}}$ is independent of the parameters $j$ and $m$.

2) Based on \eqref{mnm}, \eqref{zeta} and \eqref{C}, we can assume wlog that
\beq\label{vjm}
\eeq
$$\mathfrak{m}_{j,m}(x,\xi,\eta)= \chi_{\Z^{x}}(j)\,\times$$
$$2^{-\frac{m}{2}}\,e^{i\,\varphi_{x,\xi,\eta,j}(t_{c,x})}\,
\zeta\left(\frac{\xi}{2^{m+j}},\,\frac{\eta\,\g'_x(2^{-j})}{2^{m+j}}\right)\,
\phi\left(\frac{\xi}{2^{m+j}}\right)\,
\varrho\left(\eta,\,\frac{\g'_x(2^{-j})}{2^{m+j}}\right)\,,$$
where here we set $\zeta:=\zeta_0$.
\end{o0}

\section{Main Term for (I) - The $L^2$-decay bound}\label{ldecsec}

The core of this section will be to show that that there exists $\ep>0$ depending on $\bar{\ep}$ as displayed in \eqref{ndeg0}, such that
\beq\label{jmngen}
\|\L_{\G,j,m}(f)\|_{L^2(\R^2)} \lesssim_{\g} \,2^{-\ep m}\,\left\|f\right\|_{L^2(\R^2)}\:.
 \eeq

\begin{o0}\label{l2}
 Recall \eqref{asymptotic0}, \eqref{tdefcritpoint} and Observation \ref{rinv}; setting
\beq\label{R}
R_j(x,s):=s\,r_j(x,s)-Q_j(x,r_j(x,s))\:\:\:\:\:\:\forall\:s\in J\:,
\eeq
then for any $\xi,\:\eta\in I$ with $\frac{\xi}{\g'_x(2^{-j})\,\eta}\in J$ we have
\beq\label{Psi}
\varphi_{x,\xi,\eta,j}(t_{c,x})= -\,\frac{\g'_x(2^{-j})}{2^{j}}\,\eta\,R_j(x,\frac{\xi}{\g'_x(2^{-j})\,\eta})\:.
\eeq
\end{o0}

Now, from \eqref{vjm} and the above observation, we need to estimate the $L^2$ bound of
\beq\label{Bjmintdec}
\eeq
$$\L_{\G,j,m}(f)(x,y)=2^{-\frac{m}{2}}\,\chi_{\Z^{x}}(j)\,\int_{\R^2}\,\hat{f}(\xi,\eta)\,
e^{-i\,\frac{\g'_x(2^{-j})}{2^{j}}\,\eta\,R_j(x,\frac{\xi}{\g'_x(2^{-j})\,\eta})}$$
$$\,e^{i\,x\,\xi}\,e^{i\,y\,\eta}\,\zeta\left(\frac{\xi}{2^{m+j}},\,\frac{\eta\,\g'_x(2^{-j})}{2^{m+j}}\right)\,\phi\left(\frac{\xi}{2^{m+j}}\right)\,
\varrho\left(\eta,\,\frac{\g'_x(2^{-j})}{2^{m+j}}\right)\,d\xi\,d\eta$$

Recall now \eqref{defvr} and the properties of the function $\phi$, in particular that $\textrm{supp}\,\phi\subset \{\frac{1}{4}<|x|<4\}$; we decompose\footnote{For notational simplicity we will use the same function $\phi$ on both sides of equality \eqref{xloc}.}
\beq\label{xloc}
\phi(\frac{\g'_x(2^{-j})}{2^{m+j-k}})=\sum_{p_1=2^m}^{2^{m+1}} \phi(\frac{\g'_x(2^{-j})}{2^{j-k}}-p_1)\:.
\eeq
Write now
$$\frac{\g'_x(2^{-j})}{2^{j}}\,\eta\,R_j(x,\frac{\xi}{\g'_x(2^{-j})\,\eta})=
\frac{\g'_x(2^{-j})}{2^{j-k}}\,\frac{\eta}{2^k}\,
R_j\left(x,\frac{\frac{\xi}{2^{m+j}}}{\frac{\g'_x(2^{-j})}{2^{m+j-k}}\,\frac{\eta}{2^k}}\right)\,,$$
and define
$$\Psi_{x,\frac{\xi}{2^{m+j}},\frac{\eta}{2^k}}(s):= - s\,\frac{\eta}{2^k}\,
R_j\left(x,\frac{\frac{\xi}{2^{m+j}}}{\frac{\eta}{2^k}}\cdot\frac{2^{m}}{s}\,\right)\:.$$
Assuming that $\frac{\g'_x(2^{-j})}{2^{j-k}}-p_1,\,\frac{\xi}{2^{m+j}},\,\frac{\eta}{2^k}\in \textrm{supp}\,\phi$
and noticing that
$$\Psi_{x,\frac{\xi}{2^{m+j}},\frac{\eta}{2^k}}(\frac{\g'_x(2^{-j})}{2^{j-k}})-\Psi_{x,\frac{\xi}{2^{m+j}},\frac{\eta}{2^k}}(p_1)=$$
$$\int_0^1 (\frac{\g'_x(2^{-j})}{2^{j-k}}-p_1)\,
\Psi'_{x,\frac{\xi}{2^{m+j}},\frac{\eta}{2^k}}(\frac{\g'_x(2^{-j})}{2^{j-k}}\,t\,+(1-t)\,p_1)\,dt\,,$$
with $$\left|\Psi'_{x,\frac{\xi}{2^{m+j}},\frac{\eta}{2^k}}(s)\right|=\left|\frac{d}{ds}\Psi_{x,\frac{\xi}{2^{m+j}},\frac{\eta}{2^k}}(s)\right|\lesssim_{\g}1\,,$$
we deduce that one can decompose our exponential phase into an absolutely convergent (Taylor) series
\beq\label{xlocc}
e^{i\,\left(\Psi_{x,\frac{\xi}{2^{m+j}},\frac{\eta}{2^k}}(\frac{\g'_x(2^{-j})}{2^{j-k}})-
\Psi_{x,\frac{\xi}{2^{m+j}},\frac{\eta}{2^k}}(p_1)\right)}=
\sum_{l=0}^{\infty}\,\frac{i^l}{l!}\,\Psi_l\left(x,\frac{\g'_x(2^{-j})}{2^{j-k}}-p_1,\frac{\xi}{2^{m+j}}, \frac{\eta}{2^k}\right),
\eeq
with each $\Psi_l(x,\cdot,\cdot,\cdot)$ being a smooth function (relative to the last three variables) with $\|\Psi_l(x,\frac{\g'_x(2^{-j})}{2^{j-k}}-p_1,\cdot,\cdot)\|_{C^{1+}((\frac{1}{10},10)^{2})}\lesssim_{\g} l^{1+}$ uniformly in the $x$-parameter.

As a consequence, based on \eqref{defvr} and \eqref{xlocc}, we deduce that
\beq\label{xlocfor}
\eeq
$$\L_{\G,j,m}(f)(x,y)=2^{-\frac{m}{2}}\,\chi_{\Z^{x}}(j)\,\sum_{l\geq 0}\frac{i^l}{l!}\,\sum_{k\in\Z}\,\sum_{p_1=2^m}^{2^{m+1}}
\phi(\frac{\g'_x(2^{-j})}{2^{j-k}}-p_1)\,$$
$$\times\int_{\R^2}\,\hat{f}(\xi,\eta)\,
e^{-i\,p_1\,\frac{\eta}{2^k}\,
R_j\left(x,\frac{\frac{\xi}{2^{m+j}}}{\frac{p_1}{2^{m}}\,\frac{\eta}{2^k}}\right)}\,e^{i\,x\,\xi}\,e^{i\,y\,\eta}\,
\phi\left(\frac{\xi}{2^{m+j}}\right)\,\phi\left(\frac{\eta}{2^{k}}\right)$$
$$\times\zeta\left(\frac{\xi}{2^{m+j}},\,\frac{\eta\,\g'_x(2^{-j})}{2^{m+j}}\right)\,
\Psi_l\left(x,\frac{\g'_x(2^{-j})}{2^{j-k}}-p_1,\,\frac{\xi}{2^{m+j}},\,\frac{\eta}{2^k}\right)\,d\xi\,d\eta$$
Now due to absolute summability in the $l-$parameter, it is enough to treat only one term given say by the value $l=1$ (or $l=0$) and hence, by abusing the notation (for simplicity) we will refer to  $\L_{\G,j,m}(f)(x,y)$ as given by the expression
$$\L_{\G,j,m}(f)(x,y)=\sum_{k\in\Z} \L_{\G,j,k,m}(f)(x,y)\,,$$
where
\beq\label{hjkmnd}
\eeq
$$\L_{\G,j,k,m}(f)(x,y):=2^{-\frac{m}{2}}\,\chi_{\Z^{x}}(j)\,\sum_{p_1=2^m}^{2^{m+1}}\,\int_{\R^2}\,\hat{f}(\xi,\eta)$$
$$\times\,
e^{-i\,p_1\,\frac{\eta}{2^k}\,
R_j\left(x,\frac{\frac{\xi}{2^{m+j}}}{\frac{p_1}{2^{m}}\,\frac{\eta}{2^k}}\right)}\,e^{i\,x\,\xi}\,e^{i\,y\,\eta}\,
\Phi\left(x,\frac{\g'_x(2^{-j})}{2^{j-k}}-p_1,\frac{\xi}{2^{m+j}},\,\frac{\eta}{2^k}\right)\,d\xi\,d\eta\,,$$
with $\Phi(x,\cdot,\cdot,\cdot)$ being a smooth compactly supported function, with
$\textrm{supp}\,\Phi(x,\cdot,\cdot,\cdot)\subset \{\frac{1}{10}<|x|<10\}^3$ and $$\|\Phi(x,\frac{\g'_x(2^{-j})}{2^{j-k}}-p_1,\cdot,\cdot)\|_{C^2((\frac{1}{10},10)^{2})}\lesssim_{\g} 1\,,$$
all uniformly in the $x$-parameter.\footnote{It is precisely this point that footnote 26 is addressing it by making the choice of $\d=2$; this is a harmless assumption as condition \eqref{mpp} would transfer in our case into $\|\Phi(x,\frac{\g'_x(2^{-j})}{2^{j-k}}-p_1,\cdot,\cdot)\|_{C^{1+}((\frac{1}{10},10)^{2})}\lesssim_{\g} 1\,,$ which is in fact all that we need for the proper summability of the bounds in \eqref{Keyestim2}.}

Applying the Fourier transform in the $y-$variable we deduce:
\beq\label{hjmmain1}
\eeq
$$\hat{\L}_{\G,j,k,m}(f)(x,\eta):=2^{-\frac{m}{2}}\,\chi_{\Z^{x}}(j)\,\times\,$$
$$\sum_{p_1=2^m}^{2^{m+1}}\int_{\R}\,\hat{f}(\xi,\eta)\,
e^{-i\,p_1\,\frac{\eta}{2^k}\,
R_j\left(x,\frac{\frac{\xi}{2^{m+j}}}{\frac{p_1}{2^{m}}\,\frac{\eta}{2^k}}\right)}\,e^{i\,x\,\xi}\,
\Phi\left(x,\frac{\g'_x(2^{-j})}{2^{j-k}}-p_1,\frac{\xi}{2^{m+j}},\,\frac{\eta}{2^k}\right)\,d\xi\,.$$

It becomes now transparent from Parseval that providing $L^2\,\rightarrow\,L^2$ bounds for the initial operator $\L_{\G,j,m}f(x,y)$ is in fact equivalent with providing $L^2\,\rightarrow\,L^2$ bounds for the operator $\hat{\L}_{\G,j,k,m}f(x,\eta)$.

With all these being said, we will prove the following

\begin{t1}\label{l2dec} With the previous notations, there exists $\ep=\ep(\bar{\ep})>0$ such that the following holds:
\beq\label{dualityop}
\|\hat{\L}_{\G,j,k,m}\|_{L^2(\R^2)}\lesssim_{\g}\,2^{-\ep m}\,\|f\|_{L^2(\R^2)}\:.
\eeq
\end{t1}

\begin{proof}

Our proof is based on three key steps:
\begin{itemize}
\item firstly, an adequate discretization in the $\xi$ and $\eta$ variables whose main purpose is to localize the multiplier's phase oscillation and allow a further decomposition into \emph{linearized} wave-packets;

\item secondly, a Gabor frame decomposition in both $\xi$ and $\eta$ of the function $\hat{f}$ adapted to the previous discretization;

\item thirdly, an extraction of the cancellation encapsulated within the phase multiplier by 1) appealing first to a $T T^{*}$ argument followed by 2) a subtle iteration of the non-stationary phase principle applied to the integral kernel of the resulting operator (as a consequence of 1)) and who reveals a correlation among the time and frequency parameters involved in the discretization(s) performed at the previous two items.
\end{itemize}

With these, we now initiate the algorithm described above:
$\newline$

\noindent\textbf{Step 1} \textsf{The $\xi$ and $\eta $ discretization}.
$\newline$

Recalling the properties of the function $\phi$, we decompose\footnote{For simplicity we maintain the same notation/function $\phi$ on both sides of the equalities below.}
\beq\label{xiloc}
\phi(\frac{\xi}{2^{m+j}})=\sum_{w=2^{\frac{m}{2}}}^{2^{\frac{m}{2}+1}} \phi(\frac{\xi}{2^{j+\frac{m}{2}}}-w)\:,
\eeq
and
\beq\label{etaloc}
\phi(\frac{\eta}{2^{k}})=\sum_{v=2^{\frac{m}{2}}}^{2^{\frac{m}{2}+1}}  \phi(\frac{\eta}{2^{k-\frac{m}{2}}}-v)\:.
\eeq

Based on \eqref{hjmmain1}, \eqref{xiloc} and \eqref{etaloc}, we have that

\beq\label{redT}
\eeq
$$\hat{\L}_{\G,j,k,m}(f)(x,\eta)=2^{-\frac{m}{2}}\,\chi_{\Z^{x}}(j)\,\sum_{p_1=2^m}^{2^{m+1}}
\sum_{v,w=2^{\frac{m}{2}}}^{2^{\frac{m}{2}+1}}\,
\phi(\frac{\eta}{2^{k-\frac{m}{2}}}-v)$$
$$\times \int_{\R}\,\hat{f}(\xi,\eta)\,
e^{-i\,p_1\,\frac{\eta}{2^k}\,
R_j\left(x,\frac{\frac{\xi}{2^{m+j}}}{\frac{p_1}{2^{m}}\,\frac{\eta}{2^k}}\right)}\,e^{i\,x\,\xi}\,
\Phi\left(x,\frac{\g'_x(2^{-j})}{2^{j-k}}-p_1,\frac{\xi}{2^{m+j}},\,\frac{\eta}{2^k}\right)\,
\phi(\frac{\xi}{2^{j+\frac{m}{2}}}-w)\,d\xi\,.$$

$\newline$
\noindent\textbf{Step 2} \textsf{The adapted Gabor frame decomposition.}
$\newline$

We first introduce the Gabor frame given by
\beq\label{Gs}
\eeq
$$\{\phi_{l,s}^{w,v}(\xi,\,\eta)\}_{{l,s\in\Z}\atop{w,v\in\{2^{\frac{m}{2}}\ldots2^{\frac{m}{2}+1}\}}}:=$$
$$\left\{\frac{
1}{2^{\frac{j}{2}+\frac{m}{4}}}\,\phi\left(\frac{\xi}{2^{j+\frac{m}{2}}}-w\right)\,e^{i\,l\,\frac{\xi}{2^{j+\frac{m}{2}}}}\,
\times\frac{1}{2^{\frac{k}{2}-\frac{m}{4}}}\,\phi\left(\frac{\eta}{2^{k-\frac{m}{2}}}-v\right)\,
e^{i\,s\,\frac{\eta}{2^{k-\frac{m}{2}}}}\right\}_{{l,s\in\Z}\atop{w,v\in\{2^{\frac{m}{2}}\ldots2^{\frac{m}{2}+1}\}}}\:.$$
We now decompose the function $\hat{f}(\xi,\eta)$ relative to the above Gabor system, \textit{i.e.}
\beq\label{gs1}
\hat{f}(\xi,\eta)\sim \sum_{{{l,s\in\Z}\atop{w,v\in\{2^{\frac{m}{2}}\ldots2^{\frac{m}{2}+1}\}}}} <\hat{f},\,\phi_{l,s}^{w,v}>\,\phi_{l,s}^{w,v}(\xi,\eta)\:,
\eeq
and set
\beq\label{gs11}
\mu^{w,v}(x,p_1,\xi,\eta):=\Phi\left(x,\frac{\g'_x(2^{-j})}{2^{j-k}}-p_1,\frac{\xi}{2^{\frac{m}{2}}},\,
\frac{\eta}{2^{\frac{m}{2}}}\right)\,
\phi\left(\xi-w\right)\,\phi\left(\eta-v\right)\:,
\eeq
\beq\label{gs12}
\tilde{\mu}_{l,s}^{w,v}(x,p_1,\xi,\eta):=
\frac{1}{2^{\frac{j+k}{2}}}\,e^{i\,l\,\frac{\xi}{2^{j+\frac{m}{2}}}}
\,e^{i\,s\,\frac{\eta}{2^{k-\frac{m}{2}}}}\,\mu^{w,v}(x,p_1,\frac{\xi}{2^{j+\frac{m}{2}}},\frac{\eta}{2^{k-\frac{m}{2}}})\,.
\eeq

With these we rewrite \eqref{redT} as
\beq\label{redTT1}
\eeq
$$\hat{\L}_{\G,j,k,m}(f)(x,\eta)=2^{-\frac{m}{2}}\,\chi_{\Z^{x}}(j)\,
\sum_{p_1=2^m}^{2^{m+1}}\sum_{{{l,s\in\Z}\atop{w,v\in\{2^{\frac{m}{2}}\ldots2^{\frac{m}{2}+1}\}}}} <\hat{f},\,\phi_{l,s}^{w,v}>$$
$$\times\left(\int_{\R}\,
 e^{-i\,p_1\,\frac{\eta}{2^k}\,R_j(x,\frac{\frac{\xi}{2^j}}{p_1\,\frac{\eta}{2^k}})}
\,e^{i\,x\,\xi}\,\tilde{\mu}_{l,s}^{w,v}(x,p_1,\xi,\eta)\,d\xi\right)\:.$$

$\newline$
\noindent\textbf{Step 3} \textsf{Cancelation via the $T T^{*}$ method and time-frequency correlation.}
$\newline$
In the first part of our approach we will appeal to the $T T^{*}$ method in order to prepare the ground for exploiting the cancelation offered by the non-zero curvature hypothesis that will be encoded in the time-frequency correlation analysis discussed in the next subsection.

$\newline$
\noindent\textbf{Step 3.1} \textsf{The $T T^{*}$ argument.}
$\newline$

Taking the $L^2$ norm of the expression defined by \eqref{redTT1} we have

\beq\label{l2estim}
\eeq
$$\|\hat{\L}_{\G,j,k,m}(f)\|_{L^2(\R^2)}^2\approx 2^{-m}\,\sum_{p_1=2^m}^{2^{m+1}} \sum_{{{l,s\in\Z}\atop{w,v\in\{2^{\frac{m}{2}}\ldots2^{\frac{m}{2}+1}\}}}}
\sum_{{{l_1,s_1\in\Z}\atop{w_1\in\{2^{\frac{m}{2}}\ldots2^{\frac{m}{2}+1}\}}}}$$
$$\int_{\R}<\hat{f},\,\phi_{l,s}^{w,v}>\,\overline{<\hat{f},\,\phi_{l_1,s_1}^{w_1,v}>}\:\:
\K^{l_1,s_1,w_1}_{l,s,w}(j,p_1,v,x)\,dx\,,$$
where in the last expression
\beq\label{Keer}
\eeq
$$\K^{l_1,s_1,w_1}_{l,s,w}(j,p_1,v,x):=\chi_{\Z^{x}}(j)\int_{\R}\left(\int_{\R}\,\tilde{\mu}_{l,s}^{w,v}(x,p_1,\xi,\eta)\,
 e^{-i\,p_1\,\frac{\eta}{2^k}\,R_j(x,\frac{\frac{\xi}{2^j}}{p_1\,\frac{\eta}{2^k}})}
\,e^{i\,x\,\xi}\,d\xi\right)$$
$$\times\overline{\left(\int_{\R}\,\tilde{\mu}_{l_1,s_1}^{w_1,v}(x,p_1,\xi_1,\eta)\,
 e^{-i\,p_1\,\frac{\eta}{2^k}\,R_j(x,\frac{\frac{\xi_1}{2^j}}{p_1\,\frac{\eta}{2^k}})}
\,e^{i\,x\,\xi_1}\,d\xi_1\right)}\,d\,\eta\:.$$
We rewrite \eqref{Keer} in the explicit form
\beq\label{Keer1}
\eeq
$$\K^{l_1,s_1,w_1}_{l,s,w}(j,p_1,v,x)=\chi_{\Z^{x}}(j)\times$$
$$\frac{1}{2^{j+k}}\,\int_{\R^3}
 e^{-i\,p_1\,\frac{\eta}{2^k}\,R_j(x,\frac{\frac{\xi}{2^j}}{p_1\,\frac{\eta}{2^k}})}
\,e^{i\,p_1\,\frac{\eta}{2^k}\,R_j(x,\frac{\frac{\xi_1}{2^j}}{p_1\,\frac{\eta}{2^k}})}\,\,e^{i\,x\,\xi}\,
\,e^{-i\,x\,\xi_1}\,e^{i\,l\,\frac{\xi}{2^{j+\frac{m}{2}}}}\,e^{-i\,l_1\,\frac{\xi_1}{2^{j+\frac{m}{2}}}}\times$$
$$e^{i\,s\,\frac{\eta}{2^{k-\frac{m}{2}}}}\,e^{-i\,s_1\,\frac{\eta}{2^{k-\frac{m}{2}}}}
\,\mu^{w,v}(x,p_1,\frac{\xi}{2^{j+\frac{m}{2}}},\frac{\eta}{2^{k-\frac{m}{2}}})\,
\overline{\mu^{w_1,v}(x,p_1,\frac{\xi_1}{2^{j+\frac{m}{2}}},\frac{\eta}{2^{k-\frac{m}{2}}})}\,d\xi_1\,d\,\xi\,d\eta\:.$$
At this point, applying the change of variable $\xi\,\rightarrow\,2^{j+\frac{m}{2}}\,\xi$, $\xi_1\,\rightarrow\,2^{j+\frac{m}{2}}\,\xi_1$ and $\eta\,\rightarrow\,2^{k-\frac{m}{2}}\,\eta$, we notice that
\beq\label{Keer2}
\eeq
$$\K^{l_1,s_1,w_1}_{l,s,w}(j,p_1,v,x)=$$
$$\chi_{\Z^{x}}(j)\times2^{j+\frac{m}{2}}\,\int_{\R^3}
 e^{-i\,\frac{p_1}{2^{\frac{m}{2}}}\,\eta\,R_j(x,\frac{\xi}{\frac{p_1}{2^m}\,\eta})}
\,e^{i\,\frac{p_1}{2^{\frac{m}{2}}}\,\eta\,R_j(x,\frac{\xi_1}{\frac{p_1}{2^m}\,\eta})}\,e^{i\,2^{j+\frac{m}{2}}\,x\,\xi}\,
\times$$
$$e^{-i\,2^{j+\frac{m}{2}}\,x\,\xi_1}\,e^{i\,l\,\xi}\,e^{-i\,l_1\,\xi_1}\,e^{i\,s\,\eta}\,e^{-i\,s_1\,\eta}
\,\mu^{w,v}(x,p_1,\xi,\eta)\,\overline{\mu^{w_1,v}(x,p_1,\xi_1,\eta)}\,d\xi_1\,d\,\xi\,d\eta\:.$$
Our goal will be to prove the following

\begin{l1}\label{kerest}[\textsf{Time-frequency correlation}] $\newline$With the previous notations, we have
\beq\label{Keyestim2}
\eeq
$$|\K^{l_1,s_1,w_1}_{l,s,w}(j,p_1,v,x)|\lesssim_{\g}2^{j+\frac{m}{2}}\,
\chi_{\Z^{x}}(j)\,\phi(\frac{\g'_x(2^{-j})}{2^{j-k}}-p_1)\times$$
$$\frac{1}{(l+2^{j+\frac{m}{2}}\,x\,-2^{\frac{m}{2}}\,r_j(x,\frac{w}{\frac{p_1}{2^m}\,v}))^2+1}\,
\frac{1}{(l_1+2^{j+\frac{m}{2}}\,x\,-2^{\frac{m}{2}}\,r_j(x,\frac{w_1}{\frac{p_1}{2^m}\,v}))^2+1}\,$$
$$\times\frac{1}{(s-s_1+2^{\frac{m}{2}}\,Q_j(x,r_j(x,\frac{w}{\frac{p_1}{2^m}\,v}))-2^{\frac{m}{2}}
\,Q_j(x,r_j(x,\frac{w_1}{\frac{p_1}{2^m}\,v})))^2+1}\:.$$
\end{l1}

Assuming for the moment that Lemma \ref{kerest} holds, we proceed with proving Theorem \ref{l2dec}.

Based on \eqref{asymptotic0} - \eqref{asymptotic00} and the Mean Value Theorem, we deduce that

\beq\label{Keyestim21}
\eeq
$$\sum_{p_1=2^m}^{2^{m+1}}|\K^{l_1,s_1,w_1}_{l,s,w}(j,p_1,v,x)|\lesssim_{\g} 2^{j+\frac{m}{2}}\,
\chi_{\Z^{x}}(j)\,\phi(\frac{\g'_x(2^{-j})}{2^{j-k+m}})\times$$
$$\frac{1}{(l+2^{j+\frac{m}{2}}\,x\,-2^{\frac{m}{2}}\,r_j(x,\frac{w}{\frac{\g'_x(2^{-j})}{2^{j-k+m}}\,v}))^2+1}\,
\frac{1}{(l_1+2^{j+\frac{m}{2}}\,x\,-2^{\frac{m}{2}}\,r_j(x,\frac{w_1}{\frac{\g'_x(2^{-j})}{2^{j-k+m}}\,v}))^2+1}\,$$
$$\times\frac{1}{(s-s_1+2^{\frac{m}{2}}\,Q_j(x,r_j(x,\frac{w}
{\frac{\g'_x(2^{-j})}{2^{j-k+m}}\,v}))-2^{\frac{m}{2}}\,Q_j(x,r_j(x,\frac{w_1}{\frac{\g'_x(2^{-j})}{2^{j-k+m}}\,v})))^2+1}\:.$$

With $v$ fixed, for notational convenience we set:

- $a_{v,w}(x):=\frac{1}{\frac{\g'_x(2^{-j})}{2^{j-k+m}}}\,\frac{w}{v}$;

- $Q_{j,m}(x,w,w_1):=2^{\frac{m}{2}}\,Q_j(x,r_j(x,\frac{w}
{\frac{\g'_x(2^{-j})}{2^{j-k+m}}\,v}))-2^{\frac{m}{2}}\,Q_j(x,r_j(x,\frac{w_1}{\frac{\g'_x(2^{-j})}{2^{j-k+m}}\,v}))$.

With these notations, \eqref{Keyestim21} becomes

\beq\label{Keyestim212}
\eeq
$$\sum_{p_1=2^m}^{2^{m+1}}|\K^{l_1,s_1,w_1}_{l,s,w}(j,p_1,v,x)|\lesssim_{\g} 2^{j+\frac{m}{2}}\,
\chi_{\Z^{x}}(j)\,\phi(\frac{\g'_x(2^{-j})}{2^{j-k+m}})\times$$
$$\frac{1}{\left( \lfloor l+2^{j+\frac{m}{2}}\,x\,-2^{\frac{m}{2}}\,r_j(x,\,a_{v,w}(x))\rfloor\right)^2}\,
\frac{1}{\left( \lfloor l_1+2^{j+\frac{m}{2}}\,x\,-2^{\frac{m}{2}}\,r_j(x,\,a_{v,w_1}(x))\rfloor\right)^2}$$
$$\times\frac{1}{\left( \lfloor s-s_1+Q_{j,m}(x,w,w_1)\rfloor\right)^2}\:.$$

Putting now together \eqref{l2estim} and \eqref{Keyestim212} we deduce that
\beq\label{l2estimVv}
\eeq
$$\|\hat{\L}_{\G,j,k,m}(f)\|_{L^2(\R^2)}^2\lesssim_{\g} 2^{-m}\,\sum_{v\in\{2^{\frac{m}{2}}\ldots2^{\frac{m}{2}+1}\}}
\sum_{{l,\,l_1\in\Z}\atop{w,w_1\approx 2^{\frac{m}{2}}}}$$
$$\left(\int_{\R} \frac{ 2^{j+\frac{m}{2}}\,
\chi_{\Z^{x}}(j)\,\phi(\frac{\g'_x(2^{-j})}{2^{j-k+m}})}{\left( \lfloor l+2^{j+\frac{m}{2}}\,x\,-2^{\frac{m}{2}}\,r_j(x,\,a_{v,w}(x))\rfloor\right)^2}\,
\frac{1}{\left( \lfloor l_1+2^{j+\frac{m}{2}}\,x\,-2^{\frac{m}{2}}\,r_j(x,\,a_{v,w_1}(x))\rfloor\right)^2}\right.$$
$$\left.\times\left(\sum_{s,s_1\in\Z}\frac{|<\hat{f},\,\phi_{l,s}^{w,v}>|\,|<\hat{f},\,\phi_{l_1,s_1}^{w_1,v}>|}{\left( \lfloor s-s_1+Q_{j,m,v}(x,w,w_1)\rfloor\right)^2}\right)\,dx\right)$$
Using Cauchy-Schwarz, we further have
\beq\label{l2estimVV}
\eeq
$$\sum_{s,s_1\in\Z}\frac{|<\hat{f},\,\phi_{l,s}^{w,v}>|\,|<\hat{f},\,\phi_{l_1,s_1}^{w_1,v}>|}{\left( \lfloor s-s_1+Q_{j,m,v}(x,w,w_1)\rfloor\right)^2}$$
$$\lesssim_{\g}
(\sum_{s\in\Z}|<\hat{f},\,\phi_{l,s}^{w,v}>|^2)^{\frac{1}{2}}\,
\left(\sum_{s,s_1\in\Z}\,\frac{|<\hat{f},\,\phi_{l_1,s_1}^{w_1,v}>|^2}{\left( \lfloor s-s_1+Q_{j,m,v}(x,w,w_1)\rfloor\right)^2}\right)^{\frac{1}{2}}$$
$$\lesssim_{\g} c_{l}^{w,v}\,c_{l_1}^{w_1,v}\,,$$
where here we set
\beq\label{coeff}
c_{l}^{w,v}=c_{l}^{w,v}(f):=(\sum_{s\in\Z}|<\hat{f},\,\phi_{l,s}^{w,v}>|^2)^{\frac{1}{2}}\:.
\eeq
and notice that for $\textbf{c}:=\{c_{l}^{w,v}\}$ one has
\beq\label{coeff1}
\|\textbf{c}\|_2^2=\sum_{l\in\Z}\sum_{v,w\approx 2^{\frac{m}{2}}}
|c_{l}^{w,v}|^2\lesssim \|f\|_2^2\:.
\eeq

With these, defining
\beq\label{ndeg01}
\eeq
$$\n_{\g,\a,\b}^{v}(a(\cdot),\textbf{c},j,m):=$$
$$\left(\frac{1}{2^{m}}\,\int_{\R_{\a}}\,2^{j+\frac{m}{2}}\,\chi_{\Z^{x}_{\b}}(j)\,\phi(\frac{\g'_x(2^{-j})}{2^{j-k+m}})\,
\left(\sum_{{l\in\Z}\atop{w\approx 2^{\frac{m}{2}}}}\,
\frac{|c_{l}^{w,v}|}{\left( \lfloor l+2^{j+\frac{m}{2}}\,x\,-2^{\frac{m}{2}}\,r_j(x,\,a_{v,w}(x))\rfloor\right)^2}\right)^2\,dx\right)^{\frac{1}{2}}\:.$$
we notice from \eqref{l2estimVv} - \eqref{ndeg01} that we just proved the following\footnote{Recall that throughout this section we considered $\a=1$ and $\b=1$ and that in this setting we often drop the $\a,\,\b-$dependence.}
\beq\label{keybd}
\|\hat{\L}_{\G,j,k,m}(f)\|_{L^2(\R^2)}^2\leq \sum_{v\approx2^{\frac{m}{2}}}\n_{\g,1,1}^{v}(a(\cdot),\textbf{c},j,m)^2\:.
\eeq
Deduce now that in order to prove \eqref{dualityop},  based on \eqref{coeff1} and \eqref{keybd}, it remains to show the following

\begin{l1}\label{nondegcont} There exists $\ep=\ep(\bar{\ep})>0$ such that for any $m\in\N$, $j\in\Z$, $1\leq \a\leq A$ and $1\leq \b\leq B$ one has that
\beq\label{ndegeq0}
\n^v_{\g,\a,\b}(a(\cdot),\textbf{c},j,m)\lesssim_{\g} 2^{-\ep\,m}\,\left(\sum_{l\in\Z}\sum_{w\approx 2^{\frac{m}{2}}}
|c_{l}^{w,v}|^2\right)^{\frac{1}{2}}\:.
\eeq
\end{l1}

\begin{proof}

In what follows, without loss of generality we fix as before $\a=1$, and $\b=1$.

Since for $x\in \R_{1}$ the image of $|r_{j}(x,\cdot)|$ is always within $[C_{1,\g},\,C_{2,\g}]$ with  $C_{2,\g}\geq C_{1,\g}>0$ two absolute constants depending only on $\g$ it is enough to prove our estimate \eqref{ndegeq0} for the integral expression in \eqref{ndeg01} restricted on compact intervals, \textit{i.e.} $x\in [(k-\frac{1}{2})\, 2^{-j},\,(k+\frac{1}{2})\,2^{-j}]$ with $k\in\Z$, since one would get almost orthogonality among blocks of coefficients of the form $\{c_{l}^{w,v}\}_{\{l\in \{(-k-\frac{1}{2}+ C_{1,\g})\,2^{\frac{m}{2}},\,(-k+ \frac{1}{2}+C_{2,\g})\,2^{\frac{m}{2}}\}\}}$.

These being said, choosing $k=0$, letting $\Z^{x}_{1}=\Z^{x}$ and assuming wlog that $\R_1\equiv\R$,  we remain to study the term
\beq\label{refkee}
\eeq
$$\frac{1}{2^{m}}\,\int_{-2^{-j-1}}^{2^{-j-1}}\,2^{j+\frac{m}{2}}\,\chi_{\Z^{x}}(j)\,\phi(\frac{\g'_x(2^{-j})}{2^{j-k+m}})\,
\left(\sum_{l,w\approx 2^{\frac{m}{2}}}\,
\frac{|c_{l}^{w,v}|}{\left( \lfloor 2^{\frac{m}{2}}( l\,2^{-\frac{m}{2}}+\,2^{j}\,x\,-\,r_{j}(x,\,a_{v,w}(x))\rfloor\right)^2}\right)^2\,dx\:.$$
After some elementary reasonings, taking in account \eqref{refkee}, and setting $n=j-k+m$, one reduces \eqref{ndegeq0} to proving that one has uniformly in $n$ and $j$
\beq\label{refkeer}
\eeq
$$\frac{1}{2^{m}}\,\int_{-2^{-j-1}}^{2^{-j-1}}\,2^{j+\frac{m}{2}}\,\chi_{\Z^{x}}(j)\,\phi(\frac{\g'_x(2^{-j})}{2^{n}})\,
\left(\sum_{l,w\approx 2^{\frac{m}{2}}}\,
\frac{|c_{l}^{w,v}|}{\left( \lfloor 2^{\frac{m}{2}}( \frac{v}{2^{\frac{m}{2}}}\,\frac{\g'_x(2^{-j})}{2^{n}}\,q_j(x,\,l\,2^{-\frac{m}{2}}+\,2^{j}\,x)\,-\,\frac{w}{2^{\frac{m}{2}}})\rfloor\right)^2}\right)^2\,dx$$
$$\lesssim 2^{-2\ep m}\,\sum_{l,w\approx 2^{\frac{m}{2}}}
|c_{l}^{w,v}|^2\:.$$

Making now the change of variable $ 2^j\,x\,\rightarrow\,x$  and letting $\mathcal{R}_m$ be a $2^{-\frac{m}{2}}-$refinement of the interval $[\frac{1}{2},\,2]$ and $\textbf{c}=\{c_{l}^{w}\}_{l,w}$ we remark that \eqref{refkeer} follows from the uniform estimate

\beq\label{refkeer111}
\Lambda(m,\tilde{q},\textbf{c})\lesssim 2^{-2\ep m}\,\sum_{l,w\in\mathcal{R}_m}\,|c_{l}^{w}|^2\:,
\eeq
where here we set
\beq\label{refkeer11121}
\Lambda(m,\tilde{q},\textbf{c}):=2^{-\frac{m}{2}}\,\int_{-\frac{1}{2}}^{\frac{1}{2}}\,
\left(\sum_{l,w\in\mathcal{R}_m}\,\frac{|c_{l}^{w}|}{\left( \lfloor 2^{\frac{m}{2}}( \tilde{q}(x,\,l+x)\,-
\,w)\rfloor\right)^2}\right)^2\,dx\:,
\eeq
and
$$\tilde{q}(x,t):=\frac{v}{2^{\frac{m}{2}}}\,\chi_{\Z^{2^{-j}\,x}}(j)\,\tilde{\phi}(\frac{\g'_{2^{-j}\,x}(2^{-j})}{2^{n}})\,\frac{\g'_{2^{-j}\,x}(2^{-j})}{2^{n}}\,q_j(2^{-j}\,x,\,t)\,,$$
with $\tilde{\phi}$ a smooth compactly supported function that is identically $1$ on the support of  $\phi$ and $0$ on the complement of the set $\{\frac{1}{20}\leq |x|\leq 20\}$.

Recalling now the convexity conditions \eqref{fstterm0} and \eqref{fstterm0q}, we notice that
\begin{itemize}
\item for each fixed $x\in [-\frac{1}{2},\,\frac{1}{2}]$ and $l\in \mathcal{R}_m$ on has that
\beq\label{ctrw}
\sum_{w\in\mathcal{R}_m}\,
\frac{1}
{\left( \lfloor 2^{\frac{m}{2}}\,\left(\tilde{q}(x,\,l\,+\,x)-w\right)\rfloor\right)^2}\lesssim 1\:;
\eeq
\item for each fixed $x\in [-\frac{1}{2},\,\frac{1}{2}]$ and $w\in \mathcal{R}_m$ one has that
\beq\label{ctrl}
\sum_{l\in\mathcal{R}_m}\,
\frac{1}
{\left( \lfloor 2^{\frac{m}{2}}\,\left(\tilde{q}(x,\,l\,+\,x)-w\right)\rfloor\right)^2}\lesssim 1\:;
\eeq
\item deduce from the above that
\beq\label{ctrw1}
\sum_{(l,w)\in\mathcal{R}_m^2}\,
\int_{-\frac{1}{2}}^{\frac{1}{2}}\frac{2^{-\frac{m}{2}}}
{\left( \lfloor 2^{\frac{m}{2}}\,\left(\tilde{q}(x,\,l\,+\,x)-w\right)\rfloor\right)^2}\,dx\lesssim 1\:.
\eeq
\end{itemize}

Take now $\ep=\frac{\bar{\ep}}{10}$ a small parameter. We introduce the following:
\begin{itemize}
\item for $(l,w)\in \mathcal{R}_m^2$ we define
\beq\label{A}
A_{\tilde{q},\ep}(l,w):=\{x\in [-\frac{1}{2},\,\frac{1}{2}]\,|\,|\tilde{q}(x,\,l\,+\,x)-w|\leq 2^{-(\frac{1}{2}-2\ep)\,m}\}\:;
\eeq
\item we further define the set of \textit{light} pairs as
\beq\label{L}
\L_{\tilde{q},\ep}:=\{(l,w)\in \mathcal{R}_m^2\,|\,|A_{\tilde{q},\ep}(l,w)|\leq 2^{-2\ep m}\}\:;
\eeq
\item now, the set of \textit{heavy} pairs are simply given by
\beq\label{H}
\H_{\tilde{q},\ep}:=\mathcal{R}_m^2\setminus \L_{\tilde{q},\ep}\:.
\eeq
\end{itemize}

The key fundamental fact here is that the non-degeneracy condition \eqref{ndeg0} is essentially equivalent with the requirement that we only have \emph{few} heavy pairs. Indeed, in order to make this claim precise we first notice that via standard Riemann summation/integral considerations relation  \eqref{ndeg0} is equivalent with

\beq\label{ndeg00}
2^{-\frac{m}{2}}\,\sum_{l\in\mathcal{R}_m} \sup_{w\in\mathcal{R}_m}\,\left(\int_{-\frac{1}{2}}^{\frac{1}{2}}\frac{1}
{\left( \lfloor 2^{\frac{m}{2}}\,\left(\tilde{q}(x,\,l\,+\,x)-w\right)\rfloor\right)^2}\,dx\right)\lesssim_{\g} 2^{-2\,\bar{\ep}\,m}\:.
\eeq
Set now
\beq\label{Hhh1}
\H^2_{\tilde{q}}(l):=\{w\,|\,(l,w)\in\H_{\tilde{q},\ep}\}\:,
\eeq
and notice that, due to \eqref{ctrw}, for every $l\in \mathcal{R}_m$ one has
\beq\label{Hh1}
\#\H^2_{\tilde{q}}(l)\lesssim 2^{4\ep\,m}\:.
\eeq
As a consequence, from \eqref{ndeg00}, we have
\beq\label{ndeg001}
2^{-\frac{m}{2}}\,\sum_{(l,w)\in\H_{\tilde{q},\ep}} \left(\int_{-\frac{1}{2}}^{\frac{1}{2}}\frac{1}
{\left( \lfloor 2^{\frac{m}{2}}\,\left(\tilde{q}(x,\,l\,+\,x)-w\right)\rfloor\right)^2}\,dx\right)\lesssim_{\g} 2^{-2\,\bar{\ep}\,m}\,2^{4\,\ep\,m}\:,
\eeq
from which we deduce
\beq\label{fhe0}
\# \H_{\tilde{q},\ep}\lesssim 2^{(\frac{1}{2}-2\,\bar{\ep}+8\ep)\,m}\leq 2^{(1-\bar{\ep})\,\frac{m}{2}}\:.
\eeq
We are now ready to prove the veridicity of \eqref{refkeer111}. We start by decomposing
\beq\label{d1}
\Lambda(m,\tilde{q},\textbf{c})\lesssim \Lambda_{L}(m,\tilde{q},\textbf{c})\,+\,\Lambda_{H}(m,\tilde{q},\textbf{c})\,,
\eeq
where
\begin{itemize}
\item the \emph{light} component
\beq\label{d1L}
\Lambda_{L}(m,\tilde{q},\textbf{c})=:2^{-\frac{m}{2}}\,\int_{-\frac{1}{2}}^{\frac{1}{2}}\,
\sum_{(l,w)\in\L_{\tilde{q},\ep}}\,\sum_{(l_1,w_1)\in \L_{\tilde{q},\ep}}\,(\ldots)\,dx\:;
\eeq

\item the \emph{heavy} component

\beq\label{d1H}
\Lambda_{H}(m,\tilde{q},\textbf{c}):=2^{-\frac{m}{2}}\,\int_{-\frac{1}{2}}^{\frac{1}{2}}\,
\sum_{(l,w)\in\H_{\tilde{q},\ep}}\,\sum_{(l_1,w_1)\in \mathcal{R}_m^2}\,(\ldots)\,dx
\eeq
\end{itemize}

$\newline$
\noindent\textbf{3.1.1} \textsf{Treating the light component.}
\medskip

By a first application of Cauchy-Schwarz, we have
$$\Lambda_{L}(m,\tilde{q},\textbf{c})\lesssim $$
$$\sum_{{(l,w)\in \L_{\tilde{q},\ep}}\atop{(l_1,w_1)\in \L_{\tilde{q},\ep}}}\,\frac{|c_{l,w}|^2}{2^{m}} \,\int_{-\frac{1}{2}}^{\frac{1}{2}}
\frac{2^{\frac{m}{2}}\,dx}
{\left(\lfloor 2^{\frac{m}{2}}\,\left(\tilde{q}(x,\,l\,+\,x)-w\right) \rfloor\right)^2\,\left(\lfloor 2^{\frac{m}{2}}\,\left(\tilde{q}(x,\,l_1\,+\,x)-w_1\right) \rfloor\right)^2}$$
$$:=\sum_{{(l,w)\in \L_{\tilde{q},\ep}}\atop{(l_1,w_1)\in \L_{\tilde{q},\ep}}}\,\frac{|c_{l,w}|^2}{2^{m}} \,\int_{A_{\tilde{q},\ep}(l,w)}
(\ldots)\,dx \,
+\, \sum_{{(l,w)\in \L_{\tilde{q},\ep}}\atop{(l_1,w_1)\in \L_{\tilde{q},\ep}}}\,\frac{|c_{l,w}|^2}{2^{m}} \,\int_{[-\frac{1}{2},\frac{1}{2}]\setminus A_{\tilde{q},\ep}(l,w)}
(\ldots)\,dx $$
$$=:\Lambda_{L}^{I}(m,\tilde{q},\textbf{c})\,+\,\Lambda_{L}^{II}(m,\tilde{q},\textbf{c})\:.$$
Applying now \eqref{ctrw1} and \eqref{L}, we have
$$\Lambda_{L}^{I}(m,\tilde{q},\textbf{c})\leq
\sum_{(l,w)\in \L_{\tilde{q},\ep}}\,\frac{|c_{l,w}|^2}{2^{m}} \,\int_{A_{\tilde{q},\ep}(l,w)}
\sum_{(l_1,w_1)\in \L_{\tilde{q},\ep}}\frac{2^{\frac{m}{2}}\,dx}
{\left(\lfloor 2^{\frac{m}{2}}\,\left(\tilde{q}(x,\,l_1\,+\,x)-w_1\right) \rfloor\right)^2}$$
$$\lesssim_{d}\sum_{(l,w)\in \L_{\tilde{q},\ep}} |c_{l,w}|^2\,|A_{\tilde{q},\ep}(l,w)|\leq 2^{-2\ep m}\,\|\textbf{c}\|_{l^2(\mathcal{R}_m^2)}^2\:.$$

For the second term we apply \eqref{ctrw1} to conclude
$$\Lambda_{L}^{II}(m,\tilde{q},\textbf{c})\lesssim$$
$$\sum_{{(l,w)\in \L_{\tilde{q},\ep}}\atop{(l_1,w_1)\in \mathcal{R}_m^2}}\,\frac{|c_{l,w}|^2}{2^{m}} \,\int_{-[\frac{1}{2},\frac{1}{2}]\setminus A_{\tilde{q},\ep}(l,w)}
\frac{2^{\frac{m}{2}}}
{2^{4 \ep m}\,\left(\lfloor 2^{\frac{m}{2}}\,\left(\tilde{q}(x,\,l_1\,+\,x)-w_1\right) \rfloor\right)^2}\,dx$$
$$\lesssim 2^{-4\ep m}\,\|\textbf{c}\|_{l^2(\mathcal{R}_m^2)}^2\:.$$

$\newline$
\noindent\textbf{3.1.2} \textsf{Treating the heavy component.}
\medskip

For this term we apply \eqref{ctrw1} and \eqref{fhe0} to deduce

$$\Lambda_{H}(m,\tilde{q},\textbf{c})\lesssim (\sum_{(l,w)\in \H_{\tilde{q},\ep}}\,|c_{l,w}|)\,(\sum_{(l_1,w_1)\in \mathcal{R}_m^2}\,|c_{l_1,w_1}|^2 )^{\frac{1}{2}}$$
$$\frac{1}{2^{m}} \,\int_{-\frac{1}{2}}^{\frac{1}{2}}
2^{\frac{m}{2}}\,
\left(\sum_{(l_1,w_1)\in \mathcal{R}_m^2}\frac{1}{\left(\lfloor 2^{\frac{m}{2}}\,\left(\tilde{q}(x,\,l_1\,+\,x)-w_1\right) \rfloor\right)^2}\right)^{\frac{1}{2}}\,dx$$
$$\lesssim (\#\H_{\tilde{q},\ep})^{\frac{1}{2}}\,2^{-\frac{m}{4}}\,\|\textbf{c}\|_{l^2(\mathcal{R}_m^2)}^2\leq 2^{-(\bar{\ep}-4\ep)\,m}\,\|\textbf{c}\|_{l^2(\mathcal{R}_m^2)}^2\:.$$
Recalling now that $\ep=\frac{\bar{\ep}}{10}$ we conclude that \eqref{refkeer111} holds.
\end{proof}

$\newline$
\noindent\textbf{Step 3.2} \textsf{Time-frequency correlation.}
$\newline$

This is the component of the proof where, based on the phase linearization procedure accomplished earlier, we will be able to 
extract the cancelation via a carefully performed integration by parts process. Of key importance is that the curvature of the phase manifests now at the linear level and that the integration by parts is performed in such a way that 1) each of the variables participates in this process and 2) the transition from one variable integration to the next variable integration is performed with the preservation of the oscillation of the phase.

With this, our main focus will be to prove Lemma \ref{kerest} which is the only thing left in order to conclude
the proof of our Theorem \ref{l2dec}.

We recall below \eqref{Keer2}:
\beq\label{keystTtts1}
\eeq
$$\K^{l_1,s_1,w_1}_{l,s,w}(j,p_1,v,x)=$$
$$\chi_{\Z^{x}}(j)\times2^{j+\frac{m}{2}}\,\int_{\R^3}
 e^{-i\,\frac{p_1}{2^{\frac{m}{2}}}\,\eta\,R_j(x,\frac{\xi}{\frac{p_1}{2^m}\,\eta})}
\,e^{i\,\frac{p_1}{2^{\frac{m}{2}}}\,\eta\,R_j(x,\frac{\xi_1}{\frac{p_1}{2^m}\,\eta})}\,e^{i\,2^{j+\frac{m}{2}}\,x\,\xi}\,
\times$$
$$e^{-i\,2^{j+\frac{m}{2}}\,x\,\xi_1}\,e^{i\,l\,\xi}\,e^{-i\,l_1\,\xi_1}\,e^{i\,s\,\eta}\,e^{-i\,s_1\,\eta}
\,\mu^{w,v}(x,p_1,\xi,\eta)\,\overline{\mu^{w_1,v}(x,p_1,\xi_1,\eta)}\,d\xi_1\,d\,\xi\,d\eta\:.$$

We first isolate the integrand in the $\xi$ variable (of course, due to the symmetry, all the reasonings below will apply unchanged to the integrand in the $\xi_1$ variable). Thus, we set

\beq\label{ita}
I_{\eta}:=\left(\int_{\R}  e^{-i\,\frac{p_1}{2^{\frac{m}{2}}}\,\eta\,R_j(x,\frac{\xi}{\frac{p_1}{2^m}\,\eta})}
\,e^{i\,2^{j+\frac{m}{2}}\,x\,\xi}\,e^{i\,l\,\xi}\,\mu^{w,v}(x,p_1,\xi,\eta)\,d\,\xi\right)\:.
\eeq
Let the phase in \eqref{ita} be defined as
\beq\label{itaphase}
\o_{\eta}(\xi):=-\frac{p_1}{2^{\frac{m}{2}}}\,\eta\,R_j(x,\frac{\xi}{\frac{p_1}{2^m}\,\eta})+
2^{j+\frac{m}{2}}\,x\,\xi\,+\,l\,\xi\:.
\eeq
The phase derivative is now given by
\beq\label{itaphase1}
\frac{d}{d\xi}\,\o_{\eta}(\xi):=-2^{\frac{m}{2}}\,r_j(x,\frac{\xi}{\frac{p_1}{2^m}\,\eta})+
2^{j+\frac{m}{2}}\,x\,+\,l\:.
\eeq
We choose now $\nu_0$, $\nu_1$ a smooth partition of unity on the real line such that
\begin{itemize}
\item $\nu_0+\nu_1=1$;
\item $\nu_0,\,\nu_1\in C^{\infty}(\R)$ with $\nu_0,\,\nu_1\geq 0$;
\item $\textrm{supp}\,\nu_0\subset\{|s|<C(\g)\}$ and $\textrm{supp}\,\nu_1\subset \R\setminus \{|s|<\frac{C(\g)}{2}\}$.
\end{itemize}
Notice that, based on \eqref{fstterm0} and Observation \ref{rinv}, for an appropriate choice of $C(\g)>0$ above we have\footnote{This is the key place where we exploit the fact that, as a result of our discretization at Step 1, our phase oscillates at the linear level.}
 \beq\label{itaphase2}
 \eeq
$$\nu_1(l+2^{j+\frac{m}{2}}\,x\,-2^{\frac{m}{2}}\,r_j(x,\frac{w}{\frac{p_1}{2^m}\,v}))>0\:\:\Rightarrow\:\:
\left|\frac{d}{d\xi}\,\o_{\eta}(\xi)\right| \phi(\xi-w)\,\phi(\eta-v)>\frac{C(\g)}{10}\:.$$

We have now two cases to discuss:

$\newline$
\noindent\textbf{Case 1.} $\nu_0(l+2^{j+\frac{m}{2}}\,x\,-2^{\frac{m}{2}}\,r_j(x,\frac{w}{\frac{p_1}{2^m}\,v}))>0$.

In this situation we might have stationary points of the phase $\o_{\eta}(\xi)$ inside the set $\{(\xi,\eta)\,|\,\phi(\xi-w)\,\phi(\eta-v)>0\}$. As a consequence, we will not perform any operation on $I_{\eta}$ but use instead the smallness of the support of the underlying measure.

$\newline$
\noindent\textbf{Case 2.} $\nu_0(l+2^{j+\frac{m}{2}}\,x\,-2^{\frac{m}{2}}\,r_j(x,\frac{w}{\frac{p_1}{2^m}\,v}))=0$.

In this situation, based on \eqref{itaphase2}, we are ``far" from the set of stationary points and hence it is advantageous to apply integration by parts. Thus
 \beq\label{iet}
 \eeq
$$I_{\eta}=\int_{\R} \frac{1}{i\,\o_{\eta}'(\xi)}\,\frac{d}{d\xi}\left(e^{i\,\o_{\eta}(\xi)}\right)\,\mu^{w,v}(x,p_1,\xi,\eta)\,d\,\xi$$
$$=\int_{\R} e^{i\,\o_{\eta}(\xi)}\,\left(\frac{d}{d\xi}(-\frac{1}{i\,\o_{\eta}'(\xi)}\,\cdot)\right)\,
\mu^{w,v}(x,p_1,\xi,\eta)\,d\,\xi$$
$$=\int_{\R} e^{i\,\o_{\eta}(\xi)}\,\frac{\sigma^{w,v}(x,p_1,\xi,\eta)}{(\o_{\eta}'(\xi))^2}\,d\,\xi\,,$$
where here $\sigma^{w,v}(x,p_1,\xi,\cdot)\in C^{2}([2^{\frac{m}{2}-1},2^{\frac{m}{2}+1}]^2)$ with
 \beq\label{iet1}
 \eeq
\begin{itemize}
\item $\textrm{supp}\,\sigma^{w,v}(x,p_1,\cdot,\cdot)\subseteq  \textrm{supp}\,\phi(\cdot-w)\times \textrm{supp}\,\phi(\cdot-v)$.
\item $\|\partial_{\eta}^{\b}\,\sigma^{w,v}(x,p_1,\xi,\cdot)\|_{L^{\infty}_{\xi,\eta}}\lesssim_{\g} 1$ for any $\b\leq 2$.
\end{itemize}
Now on the support of $\sigma^{w,v}$ we write
 \beq\label{iet2}
 \eeq
$$\frac{1}{(\o_{\eta}'(\xi))^2}=\frac{1}{(\o_{v}'(w))^2\left(1-\frac{\o_{v}'(w)-\o_{\eta}'(\xi)}{\o_{v}'(w)}\right)^2}$$
$$=\frac{1}{(\o_{v}'(w))^2}\sum_{k\geq 1} k\,\left(\frac{\o_{v}'(w)-\o_{\eta}'(\xi)}{\o_{v}'(w)}\right)^{k-1}\,,$$
where here  $\left|\frac{\o_{v}'(w)-\o_{\eta}'(\xi)}{\o_{v}'(w)}\right|\leq \frac{1}{100}$ if $C(\g)$ is chosen large enough.

Based on \eqref{iet} - \eqref{iet2} we can safely rewrite the expression in  \eqref{iet1} as
\beq\label{iet3}
I_{\eta}= \frac{1}{(\o_{v}'(w))^2}\,\int_{\R} e^{i\,\o_{\eta}(\xi)}\,\sigma^{w,v}_0(x,p_1,\xi,\eta)\,d\,\xi\,,
\eeq
where here $\sigma_0^{w,v}$ is a function with the same properties as $\sigma^{w,v}$ - see \eqref{iet1}.

In conclusion, for both Cases 1 and 2 one has that
\beq\label{iet4}
I_{\eta}= \frac{1}{1+(\o_{v}'(w))^2}\,\int_{\R} e^{i\,\o_{\eta}(\xi)}\,\sigma^{w,v}_0(x,p_1,\xi,\eta)\,d\,\xi\,,
\eeq
with $\sigma_0^{w,v}$ obeying \eqref{iet1}.

Proceeding in the same way, one has that
\beq\label{ita5}
\eeq
$$I_{\eta}^1:=\int_{\R}  e^{i\,\frac{p_1}{2^{\frac{m}{2}}}\,\eta\,R_j(x,\frac{\xi_1}{\frac{p_1}{2^m}\,\eta})}
\,e^{-i\,2^{j+\frac{m}{2}}\,x\,\xi_1}\,e^{-i\,l_1\,\xi_1}\,\overline{\mu^{w_1,v}(x,p_1,\xi_1,\eta)}\,d\,\xi_1$$
$$=\frac{1}{1+({\o_{v}^{1}}'(w_1))^2}\,\int_{\R} e^{i\,\o_{\eta}^{1}(\xi_1)}\,\sigma^{w_1,v}_1(x,p_1,\xi_1,\eta)\,d\,\xi_1\,,$$
where here
\beq\label{itaphase21}
\o_{\eta}^{1}(\xi_1):=\frac{p_1}{2^{\frac{m}{2}}}\,\eta\,R_j(x,\frac{\xi_1}{\frac{p_1}{2^m}\,\eta})-
2^{j+\frac{m}{2}}\,x\,\xi_1\,-\,l_1\,\xi_1\:,
\eeq
and $\sigma^{w_1,v}_1$ obeying the same properties as $\sigma_0^{w,v}$.

Inserting now \eqref{iet4} and  \eqref{ita5} in \eqref{keystTtts1} we deduce that
\beq\label{keret}
\eeq
$$|\K^{l_1,s_1,w_1}_{l,s,w}(j,p_1,v,x)|\lesssim 2^{j+\frac{m}{2}}\,\chi_{\Z^{x}}(j)\,
\frac{1}{1+(\o_{v}'(w))^2}\,\frac{1}{1+({\o_{v}^1}'(w_1))^2}\times$$
$$\int_{\R^2}\left|\int_{\R} e^{i\,\o_{\eta}(\xi)}\,e^{i\,\o_{\eta}^{1}(\xi_1)}\,e^{i\,s\,\eta}\,e^{-i\,s_1\,\eta}
\,\sigma^{w,v}_0(x,p_1,\xi,\eta)\,
\sigma^{w_1,v}_1(x,p_1,\xi_1,\eta)\,d\,\eta\right|
\,d\xi_1\,d\,\xi\:.$$
Set now
$$|I(\xi,\xi_1)|:=\left|\int_{\R} e^{i\,\o_{\eta}(\xi)}\,e^{i\,\o_{\eta}^{1}(\xi_1)}\,e^{i\,s\,\eta}\,e^{-i\,s_1\,\eta}
\,\sigma^{w,v}_0(x,p_1,\xi,\eta)\,\sigma^{w_1,v}_1(x,p_1,\xi_1,\eta)\,d\,\eta\right|\,,$$
and notice that $I(\xi,\xi_1)$ can be written as
$$I(\xi,\xi_1)=\int_{\R} e^{i\,\o_{\xi,\xi_1}(\eta)}\,\tilde{\sigma}^{w,w_1,v}(x,p_1,\xi,\xi_1,\eta)\,d\,\eta\,,$$
where
\beq\label{ph12}
\o_{\xi,\xi_1}(\eta):=-\frac{p_1}{2^{\frac{m}{2}}}\,\eta\,R_j(x,\frac{\xi}{\frac{p_1}{2^m}\,\eta})+
\frac{p_1}{2^{\frac{m}{2}}}\,\eta\,R_j(x,\frac{\xi_1}{\frac{p_1}{2^m}\,\eta})\,+\,s\,\eta\,-\,s_1\,\eta\:,
\eeq
and $\tilde{\sigma}^{w,w_1,v}(x,p_1,\xi,\xi_1,\eta):=\sigma^{w,v}_0(x,p_1,\xi,\eta)\,\sigma^{w_1,v}_1(x,p_1,\xi_1,\eta)$.

Based on \eqref{R}, it remains now to notice that
\beq\label{ph121}
\frac{d}{d\eta}\,\o_{\xi,\xi_1}(\eta)=\frac{p_1}{2^{\frac{m}{2}}}\,
\left(Q_j(x,r_j(x,\frac{\xi}{\frac{p_1}{2^m}\,\eta}))-Q_j(x,r_j(x,\frac{\xi_1}{\frac{p_1}{2^m}\,\eta}))\right)\,+\,s\,-\,s_1\:.
\eeq
Repeating now the case discussion from above - that is splitting the domain of integration into two regions close and far from the stationary points - and then, in the latter situation, integrating by parts the resulting components at most twice and using the properties of $\g$ and $\tilde{\sigma}$ we deduce that
\beq\label{ph1eta}
|I(\xi,\xi_1)|\lesssim_{\g}\frac{1}{1+(\o_{w,w_1}'(v))^2}\,\phi(\xi-w)\,\phi(\xi_1-w_1)\,\phi(\frac{\g'_x(2^{-j})}{2^{j-k}}-p_1)\,.
\eeq
Putting now together \eqref{keret} and \eqref{ph1eta}, we deduce that \eqref{Keyestim2} holds, thus ending the proof of our Theorem \ref{l2dec}.
\end{proof}

\noindent\textbf{Proof of Theorem \ref{hmdiag2}}

Relying on \eqref{dualityop}, we have that
\beq\label{ph1etare}
\|\L_{\G,j,k,m} f\|_{L^2(\R^2)}\lesssim_{\g}\,2^{-\ep m}\,\|f*_{x} \check{\phi}_{m+j}*_{y}\check{\phi}_{k}\|_{L^2(\R^2)}\:.
\eeq
Thus, appealing now to duality (as we did several times by now), we deduce that\footnote{For more details on this, one can read the reasonings presented at Stage 1 inside the proof of Theorem \ref{hmdiagp} below, reasonings that are very similar in nature.}
$$\|\M_{\G,m}(f)\|_2^2,\,\|\H_{\G,m}(f)\|_2^2\lesssim_{\g}
\left\|\left(\sum_{{j\in\Z}\atop{k\in\Z}} \left|\L_{\G,j,k,m} f\right|^2 \right)^{\frac{1}{2}}\right\|_{2}^2$$
$$\lesssim_{\g} 2^{-2\,\ep\,m}\,\sum_{{j\in\Z}\atop{k\in\Z}} \|f*_{x} \check{\phi}_{m+j}*_{y}\check{\phi}_{k}\|_{L^2(\R^2)}^2\lesssim_{\g} 2^{-2\,\ep\,m}\,\|f\|_2^2\:,$$
where in the last line we used the standard Littlewood-Paley theory.

\section{The $L^p$ bound}\label{LPS}

Our goal in this section is to provide the proof of Theorem \ref{Diagp}.

As discussed in Section \ref{statph}, the desired $L^p$-bound control in Theorem \ref{Diagp} follows trivially by interpolating our $L^2-$case proved in Theorem \ref{hmdiag2} with the result stated in Theorem \ref{hmdiagp}. Thus, all that remains is to prove Theorem \ref{hmdiagp}, that is, to show that
\beq\label{HMmp0}
\|\M_{\G,m}(f)\|_p,\,\|\H_{\G,m}(f)\|_p\lesssim_{\g,p}
\left\|\left(\sum_{{j\in\Z}\atop{k\in\Z}} \left|\L_{\G,j,k,m} f\right|^2 \right)^{\frac{1}{2}}\right\|_{p}\lesssim_{\g,p} m^2\,\|f\|_p\:.
\eeq

\begin{o0}\label{genpbounds} [\textsf{Fourier $L^2-$methods versus spatial $L^p$-methods}] In any of the themes discussed in the present paper, there is an important, philosophical distinction between the $L^2$-approach and the corresponding $L^p$, $p\not=2$, one, based on the scale-type decay that one needs to extract from the non-zero curvature hypothesis. Indeed, on the one hand, the $L^2$ bound requires a decay in terms of the height of the phase of the multiplier and thus this type of information can only be obtained by appealing to Fourier transform with the aim of exploiting the cancellation captured within the multiplier. Thus, in this situation, one needs to rely on Parseval, frequency discretization of the kernel/operator and ultimately $T T^{*}-$method. In contrast with this, the $L^p-$bound, $p\not=2$, does not require a scale-type decay but only a moderate growth in the corresponding parameter as witnessed by \eqref{HMmp0}. This latter bound is better adapted to \textit{spatial} methods in which the discretization of our operator moves its weight from the oscillation of the multiplier's phase to the location of the input function.

Thus, there should come as no surprise that in our $L^p$-approach(es), $p\not=2$, one starts by performing backwards the frequency decomposition in order to reach to the original spatial localization of the operator.
\end{o0}

$\newline$
\noindent\textit{Proof of Theorem \ref{hmdiagp}.}
$\newline$

We split our proof in two components corresponding to the two inequalities in \eqref{HMmp}. Thus:
$\newline$

\noindent\textbf{Stage 1}. The following holds:
\beq\label{HMmpI}
\|\M_{\G,m}(f)\|_p,\,\|\H_{\G,m}(f)\|_p\lesssim_{\g,p}
\left\|\left(\sum_{{j\in\Z}\atop{k\in\Z}} \left|\L_{\G,j,k,m} f\right|^2 \right)^{\frac{1}{2}}\right\|_{p}\:.
\eeq
$\newline$

We notice that \eqref{HMmpI} is trivially true for $\M_{\G,m}$ since we have the straightforward pointwise estimate
\beq\label{Mmp1}
\M_{\G,m}(f)(x,y)\lesssim
\left(\sum_{j\in\Z} \left|\L_{\G,j,m} f (x,y)\right|^2 \right)^{\frac{1}{2}}\:,
\eeq
while for each $(x,y)\in\R^2$ and $j,\,m$ fixed at most $10$ terms from $\{\phi(\frac{\g'_x(2^{-j})}{2^{m+j-k}})\}_{k\in\Z}$ are non-zero and thus one further has
\beq\label{Mmp2}
 \left|\L_{\G,j,m} f (x,y)\right|\lesssim \left(\sum_{k\in\Z} \left|\L_{\G,j,k,m} f (x,y)\right|^2 \right)^{\frac{1}{2}}\:.
\eeq

Thus, the entire focus falls onto proving this estimate for $\H_{\G,m}$. For this, given $m\in\N$ and $f\in L^p(\R^2)$, $g\in L^{p'}(\R^2)$ with $1<p<\infty$ and $\frac{1}{p}+\frac{1}{p'}=1$, we define
\beq\label{ljm}
\ll_{m}(f,g):=\int_{\R^2}  \H_{\G,m}(f)(x,y)\,g(x,y)\,dx\,dy\:.
\eeq

Recall now, that with the previous notations, we can visualize the following decomposition
\beq\label{hmkj}
\H_{\G,m}(f)(x,y)=\sum_{j\in\Z}\sum_{k\in\Z} \L_{\G,j,k,m}(f)(x,y)\,,
\eeq
where, recalling \eqref{hjmmain1}, one has
\beq\label{hjkmndd}
\L_{\G,j,k,m}(f)(x,y)=\sum_{p_1=2^m}^{2^{m+1}}\L_{\G,j,k,m}^{p_1}(f)(x,y)\,,
\eeq
with
\beq\label{hjkmnddd}
\eeq
$$\L_{\G,j,k,m}^{p_1}(f)(x,y):=2^{-\frac{m}{2}}\,\chi_{\Z^{x}}(j)\,\int_{\R^2}\,\hat{f}(\xi,\eta)$$
$$\times\,
e^{-i\,p_1\,\frac{\eta}{2^k}\,
R_j\left(x,\frac{\frac{\xi}{2^{m+j}}}{\frac{p_1}{2^{m}}\,\frac{\eta}{2^k}}\right)}\,e^{i\,x\,\xi}\,e^{i\,y\,\eta}\,
\Phi\left(x,\frac{\g'_x(2^{-j})}{2^{j-k}}-p_1,\frac{\xi}{2^{m+j}},\,\frac{\eta}{2^k}\right)\,d\xi\,d\eta\,,$$
and $\Phi(x,\cdot,\cdot,\cdot)$ smooth with $\textrm{supp}\,\Phi(x,\frac{\g'_x(2^{-j})}{2^{j-k}}-p_1,\cdot,\cdot)\subset \{\frac{1}{10}<|x|<10\}^2$ and
\beq\label{ph}
\|\Phi(x,\frac{\g'_x(2^{-j})}{2^{j-k}}-p_1,\cdot,\cdot)\|_{C^2((\frac{1}{10},10)^{2})}\lesssim_{\g} 1\,.
\eeq

From \eqref{mdecomp} - \eqref{hjkmnd}, we can assume wlog\footnote{We are simply proceeding backwards in our multiplier decomposition.} that
\beq\label{pph}
\eeq
$$\L_{\G,j,k,m}^{p_1}(f)(x,y)=2^{-\frac{m}{2}}\,\chi_{\Z^{x}}(j)\,\phi(\frac{\g'_x(2^{-j})}{2^{j-k}}-p_1)\,\int_{\R^2}\,\hat{f}(\xi,\eta)\,
e^{i\,\Psi_{x,\frac{\xi}{2^{m+j}},\frac{\eta}{2^k}}(\frac{\g'_x(2^{-j})}{2^{j-k}})}$$
$$e^{i\,x\,\xi}\,e^{i\,y\,\eta}\,\zeta\left(\frac{\xi}{2^{m+j}},\,\frac{\eta\,\g'_x(2^{-j})}{2^{m+j}}\right)\,\phi\left(\frac{\xi}{2^{m+j}}\right)\,
\phi(\frac{\eta}{2^k})\,d\xi\,d\eta\,,$$
and hence by possibly adding similar nature terms we may further assume
$$\L_{\G,j,k,m}^{p_1}(f)(x,y)\approx\chi_{\Z^{x}}(j)\,\phi(\frac{\g'_x(2^{-j})}{2^{j-k}}-p_1)\,
\int_{\R^2}\,\hat{f}(\xi,\eta)\,e^{i\,x\,\xi}\,e^{i\,y\,\eta}\,$$
$$\left(\int_\R  e^{-i\xi 2^{-j}t + i\eta \gamma_x(2^{-j}t)} \rho(t) dt\right)\,\phi\left(\frac{\xi}{2^{m+j}}\right)\,
\phi(\frac{\eta}{2^k})\,d\xi\,d\eta\:.$$
Thus, from now on, we will identify
\beq\label{phys}
\eeq
$$\L_{\G,j,k,m}^{p_1}(f)(x,y)\equiv$$
$$\chi_{\Z^{x}}(j)\,\phi(\frac{\g'_x(2^{-j})}{2^{j-k}}-p_1)\,
\int_{\R} (f*_{x}\phi_{m+j}*_{y}*\phi_{k})\left(x-\frac{t}{2^j},\,y+\gamma_x(\frac{t}{2^j})\right)\,\rho(t) dt\,.$$

We now have

$$\left| \ll_{m}(f,g)\right|=\left|\int_{\R^2} \sum_{j\in\Z}\sum_{k\in\Z} \sum_{p_1=2^m}^{2^{m+1}} \L_{\G,j,k,m}^{p_1}(f)(x,y)\, g(x,y)\,dx\,dy\right|$$
$$\leq \int_{\R^2} \left(\sum_{j\in\Z}\sum_{k\in\Z} \sum_{p_1=2^m}^{2^{m+1}} |\L_{\G,j,k,m}^{p_1}(f)(x,y)|^2\right)^{\frac{1}{2}}$$\,
$$\times\left(\sum_{j\in\Z^{x}}\sum_{k\in\Z} \sum_{p_1=2^m}^{2^{m+1}} \phi(\frac{\g'_x(2^{-j})}{2^{j-k}}-p_1)\, |(g*_{y}\phi_{k})(x,y)|^2\right)^{\frac{1}{2}}\,dx\,dy$$
$$\leq A_{p}(f)\,B_{p'}(g)\,,$$
where here we set
\beq\label{AP}
 A_{p}(f):=\left\|\left(\sum_{j\in\Z}\sum_{k\in\Z} \sum_{p_1=2^m}^{2^{m+1}} |\L_{\G,j,k,m}^{p_1}(f)(x,y)|^2\right)^{\frac{1}{2}} \right\|_{L^p(\R^2)}\,,
\eeq
and
\beq\label{BP'}
B_{p'}(g):=\left\|\left(\sum_{j\in\Z^{x}}\sum_{k\in\Z} \sum_{p_1=2^m}^{2^{m+1}} \phi(\frac{\g'_x(2^{-j})}{2^{j-k}}-p_1)\, |(g*_{y}\phi_{k})(x,y)|^2\right)^{\frac{1}{2}} \right\|_{L^{p'}(\R^2)}\,.
\eeq
Now the last term is easy to estimate based on the previous considerations about our curve $\g$ together with standard Littlewood-Paley theory; indeed, based on property \eqref{variation0}, we know that\footnote{Compare this relation with the closely related \eqref{jcgam} and \eqref{gs}.}
$$\sum_{j\in\Z^{x}}\sum_{p_1=2^m}^{2^{m+1}} \phi(\frac{\g'_x(2^{-j})}{2^{j-k}}-p_1)\lesssim_{\g} 1\,,$$
and thus
$$B_{p'}(g)\lesssim_{\g}\left\|\left(\sum_{k\in\Z} |(g*_{y}\phi_{k})(x,y)|^2\right)^{\frac{1}{2}} \right\|_{L^{p'}(\R^2)}\lesssim_{p'} \|g\|_{L^{p'}(\R^2)}\:.$$
This ends the proof of \eqref{HMmpI}.

$\newline$

\noindent\textbf{Stage 2.} The following holds:
\beq\label{apest}
 A_{p}(f):=\left\|\left(\sum_{j\in\Z}\sum_{k\in\Z} \sum_{p_1=2^m}^{2^{m+1}} |\L_{\G,j,k,m}^{p_1}(f)|^2\right)^{\frac{1}{2}} \right\|_{L^p(\R^2)}
 \lesssim_{\g,p} m^2\, \left\|f\right\|_{L^p(\R^2)}\:.
\eeq
$\newline$

Relation  \eqref{apest} is just a direct consequence of Lemma \ref{translk} for $m=n\in\N$.

\section{Other approaches for the $L^p$-boundedness}\label{otherlp}

In this section we want to provide several other perspectives/angles for approaching the problem of providing $L^p-$bounds for the operators $H_{\Gamma}$ and $M_{\Gamma}$. Each of these other methods comes with it's own ``personality" and is meant to enrich one's understanding of the present topic.

In the first two subsections we discuss different routes for proving\footnote{The logarithmic loss in the RHS may have different exponents, that is instead of $m^2$ we allow $m^{c(p)}$ for some $c(p)>0$.} \eqref{apest}, which is the key estimate needed in order to obtain the desired $L^p-$control for $p\not=2$. The third subsection treats only the case $1<p<2$, while the last one makes the reference to an alternative - but less powerful - route to obtain the fundamental $m-$decay in the case $p=2$.

\subsection{Other approaches (I.1) - the case $1<p<\infty$ with $p\not=2$}

As mentioned in the introduction above, our intention is to provide a different proof for
\beq\label{apestt}
\left\|\left(\sum_{j\in\Z}\sum_{k\in\Z} \sum_{p_1=2^m}^{2^{m+1}} |\L_{\G,j,k,m}^{p_1}(f)|^2\right)^{\frac{1}{2}} \right\|_{L^p(\R^2)}
 \lesssim_{\g,p} m^{c(p)}\, \left\|f\right\|_{L^p(\R^2)}\:,
\eeq
for some $c(p)>0$.

We will show that \eqref{apestt} holds with $c(p)=4$. This is an immediate consequence of  Propositions \ref{maxdom} and \ref{maxestim} below.

\begin{p1}\label{maxdom} Fix $(x,\,y)\in\R^2$ and $j\in\Z$. With the previous notations, recalling \eqref{Hmm1} and \eqref{Hmm10} in Section \ref{Not}, one has the following pointwise estimate:
\beq\label{dommax}
\eeq
$$\left(\sum_{k\in\Z} \sum_{p_1=2^m}^{2^{m+1}} |\L_{\G,j,k,m}^{p_1}(f)|^2(x,y)\right)^{\frac{1}{2}} \lesssim_{\g}\frac{1}{2^m}\sum_{l=2^{m-2}}^{2^{m+2}}\,M_1^{[l]}
 M_2^{[2^m]} (f*_{x}\phi_{m+j})(x,\,y)\,.$$
\end{p1}

\begin{proof}

The proof of this proposition is inspired from the $L^p$-approach ($p\not=2$) developed in \cite{GHLR}. This in turn, was further inspired by the $L^p$-Banach case treatment of the Bilinear Hilbert transform along ``non-flat" curves in \cite{lv10}. For more about this connection, the reader is invited to read the next (sub)section and the follow up study \cite{lvUA3}.

We first notice the following simple fact:
$$\sum_{p_1=2^m}^{2^{m+1}}\phi(\frac{\g'_x(2^{-j})}{2^{j-k}}-p_1)\lesssim_{\g} \phi(\frac{\g'_x(2^{-j})}{2^{m+j-k}})\:.$$

Appealing to the definitions of $\phi_{m+j}$ and $\phi_{k}$, we further deduce\footnote{Strictly speaking we should write $\check{\tilde{\phi}}_{m+j}(s)$ instead of $\check{\phi}_{m+j}(s)$ where $\tilde{\phi}_{m+j}(s)=1$ on the support of $\check{\phi}_{m+j}(s)$. However for notational simplicity we drop the symbol $\tilde{}$. Same applies to $\phi_{k}(u)$. Also here we only treat the case $\textrm{supp}\,\r\subseteq \{\frac{1}{4}<t<1\}$ since the symmetric case (the support of $\r$ is restricted within $\textrm{supp}\, \r\subseteq \{-1<t<-\frac{1}{4}\}$) only requires trivial modifications.}
$$I_{\Gamma, j}f(x,y):=\left(\sum_{k\in\Z} \sum_{p_1=2^m}^{2^{m+1}} |\L_{\G,j,k,m}^{p_1}(f)|^2(x,y)\right)^{\frac{1}{2}}\lesssim_{\g} $$
$$\phi(\frac{\g'_x(2^{-j})}{2^{m+j-k_0}})\,\int_{\R}\left(\int_{\R^2} |f*_{x}\phi_{m+j}|\left(x-\frac{t}{2^j}-s,\,y+\gamma_x(\frac{t}{2^j})-u\right)\,
|\check{\phi}_{m+j}(s)|\,|\check{\phi}_{k_0}(u)|\,ds\,du\right)\,\r(t)\,dt$$
where here $k_0\in\Z$ is - up to at most $C(\g)$ consecutive values - the ``only" integer $k$ depending on $j$ and $x$, for which $\phi(\frac{\g'_x(2^{-j})}{2^{m+j-k}})\not=0$.

From here, we further deduce that
$$I_{\Gamma, j}f(x,y)\lesssim_{\g}\sum_{l=2^{m-2}}^{2^{m+2}}\int_{\frac{l}{2^m}}^{\frac{l+1}{2^m}}\left(\int_{\R^2} |f*_{x}\phi_{m+j}|\left(x-\frac{t}{2^j}-s,\,y+\gamma_x(\frac{t}{2^j})-u\right)\,\right.$$
$$\times\left.\frac{2^{m+j}}{(2^{m+j}\,s)^2+1}\,\frac{2^{k_0}}{(2^{k_0}\,u)^2+1}\,ds\,du\right)\,dt$$
which after making the change of variables $t \rightarrow \frac{t+l}{2^m}$, $s \rightarrow \frac{s-t-l}{2^{m+j}}$ and
$u \rightarrow u + \gamma_x(\frac{t+l}{2^{m+j}})$ becomes
\beq\label{cv1con}
\eeq
$$I_{\Gamma, j}f(x,y)\lesssim_{\g}\frac{1}{2^m}\sum_{l=2^{m-2}}^{2^{m+2}}\int_{0}^{1}\left(\int_{\R^2} |f*_{x}\phi_{m+j}|(x-\frac{s}{2^{m+j}},\,y-u)\,\right.$$
$$\left.\times\frac{1}{(s-t-l)^2+1}\,\frac{2^{k_0}}{(2^{k_0}\,(u + \gamma_x(\frac{t+l}{2^{m+j}})))^2+1}\,ds\,du\right)\,dt\:.$$
From  \eqref{asymptotic0}-\eqref{fstterm0q} and the choice of $k_0$ we notice that
\beq\label{ycon}
|2^{k_0}\,\gamma_x(\frac{t+l}{2^{m+j}})-2^{k_0}\,\gamma_x(\frac{l}{2^{m+j}})|\lesssim_{\g} 1\,,
\eeq
which reduces \eqref{cv1con} to
\beq\label{cv1con1}
\eeq
$$I_{\Gamma, j}f(x,y)\lesssim_{\g}\frac{1}{2^m}\sum_{l=2^{m-2}}^{2^{m+2}}\int_{\R^2} |f*_{x}\phi_{m+j}|(x-\frac{s}{2^{m+j}},\,y-u)\,$$
$$\times\frac{1}{(s-l)^2+1}\,\frac{2^{k_0}}{(2^{k_0}\,(u + \gamma_x(\frac{l}{2^{m+j}})))^2+1}\,ds\,du\:.$$
Based now on \eqref{asymptotic0}-\eqref{fstterm0q} we notice that - except possibly for at most $C(\g)$ values of $j$ for which one can apply the same reasonings as for the set $\Z_{0}^{x}$ - one has
$\phi(\frac{\g'_x(2^{-j})}{2^{m+j-k_0}})\not=0$ iff there exists $c_\g\not=0$ such that for each $l\approx 2^m$ one has
$\phi\left(c_\g\,2^{k_0-m}\,\gamma_x(\frac{l}{2^{m+j}})\right)\not=0$. Given $l$ and assuming for notational simplicity that $c_{\g}>0$, we define $k_l$ such that
$2^{-k_l}:=2^{-m}\,\gamma_x(\frac{l}{2^{m+j}})$. With these, we have
\beq\label{cv1con2}
\eeq
$$I_{\Gamma, j}f(x,y)\lesssim_{\g}\frac{1}{2^m}\sum_{l=2^{m-2}}^{2^{m+2}}\int_{\R^2} |f*_{x}\phi_{m+j}|(x-\frac{s}{2^{m+j}},\,y-u)\,$$
$$\times\frac{1}{(s-l)^2+1}\,\frac{2^{k_l}}{(2^{k_l}\,u + 2^m)^2+1}\,ds\,du\:,$$
which by taking now a supremum over $k_l$ proves the validity of \eqref{dommax}.
\end{proof}

\begin{p1}\label{maxestim} Let $1<p<\infty$. Assume $l\in\N$ with $l\approx 2^m$ is a given, fixed parameter. Then, the following holds:
\beq\label{thpmain1}
\left\|\left(\sum_{j\in\Z}\,\left|M_1^{[l]}
 M_2^{[2^m]} (f*_{x}\phi_{m+j})(x,\,y)\right|^2\right)^{\frac{1}{2}} \right\|_{L^p(\R^2)}\lesssim  m^{4}\,\left\|f\right\|_{L^p(\R^2)} \:.
\eeq
\end{p1}

\begin{proof}
The proof of \eqref{thpmain1} follows via classical Littlewood-Paley theory from applying twice the following one-dimensional vector-valued inequality:
\beq\label{thpmain11}
\left\|\left(\sum_{j\in\Z} \left|M^{(2^m)}(f_j)\right|^2\right)^{\frac{1}{2}} \right\|_{L^p(\R)}\lesssim  m^{2}\,\left\|(\sum_{j\in\Z}|f_j|^2)^{\frac{1}{2}}\right\|_{L^p(\R)}\:,
\eeq
with the latter being a particular case of Theorem 3.1 in \cite{GHLR}.
\end{proof}

Conclude now that \eqref{apestt} follows from Propositions \ref{maxdom} and \ref{maxestim} and an application of Minkowski's inequality for the $l-$summation.

\subsection{Other approaches (I.2) - the case $1<p<\infty$ with $p\not=2$} \label{LPg2}

As already mentioned, the proof of Proposition \ref{maxdom} was inspired by the general $L^p-$approach in \cite{GHLR}, which was further inspired by the Banach triangle case treatment of the ``non-flat" Bilinear Hilbert transform in \cite{lv10} in which we introduced a shifted square function argument. However, in contrast with the approach provided by Proposition \ref{maxdom} which is tailored for the spatial variables,  the approach in \cite{lv10} focuses more on the Fourier variable side and in fact can be adapted to give another alternative approach to the $L^p$ bounds discussed here. Due to space limitations concerns, we will leave this enlightening $L^p$-parallelism between space and Fourier methods for the second part of our study - see \cite{lvUA3}.

\subsection{Other approaches (I.3) - the case $1<p<2$}

In this section we want to provide a different, though \textit{indirect} route of proving that
\beq\label{HMlp}
\|\M_{\G,m}(f)\|_p,\,\|\H_{\G,m}(f)\|_p\lesssim_{p}\left\|\left(\sum_{j\in\Z} \left|\L_{\G,j,m} f\right|^2 \right)^{\frac{1}{2}}\right\|_{p}\lesssim_{\g,p}\|f\|_p\:,
\eeq
for $1<p<2$, where here we set $\L_{\G,j,m}(f):=\sum_{k\in\Z}\L_{\G,j,k,m}(f)$.
$\newline$

\begin{r1}\label{equivmaxsq} This indirect approach has the merit of showing that under suitable conditions on the operators under discussion and on the range of $p$, control on maximal operators is in fact \underline{\emph{equivalent}} with control over the associated square functions.
\end{r1}

In what follows, we will skip the technical considerations that were anyhow presented in great detail in the previous sections, and insist only on the main ideas.

Recalling \eqref{hjkmndd} and \eqref{phys}, we notice that
\beq\label{physs}
\eeq
$$\L_{\G,j,k,m}(f)(x,y)$$
$$\approx\chi_{\Z^{x}}(j)\,\phi(\frac{\g'_x(2^{-j})}{2^{m+j-k}})\,
\int_{\R} (f*_{x}\phi_{m+j}*_{y}\phi_{k})\left(x-\frac{t}{2^j},\,y+\gamma_x(\frac{t}{2^j})\right)\,\rho(t) dt\,.$$

Next, \footnote{For example by possibly slightly modifying the partition of unity in \eqref{firstl4} into $1=\sum_{m,n\in\Z}\phi(\frac{\xi}{2^{m+j}})\,\phi(\frac{\g'_x(2^{-j})\,\eta}{2^{n+j}})$.} we have that
\beq\label{physs1}
\L_{\G,j,m}(f)(x,y)\approx\chi_{\Z^{x}}(j)\,
\int_{\R} (f*_{x}\check{\phi}_{m+j}*_{y}\check{\phi}_{m,j,x})\left(x-\frac{t}{2^j},\,y+\gamma_x(\frac{t}{2^j})\right)\,\rho(t) dt\,,
\eeq
where here we set $\phi_{m,j,x}(\eta):=\phi(\frac{\g'_x(2^{-j})\,\eta}{2^{m+j}})$.

Following similar reasonings with those in Section \ref{LPS}, we have that
\beq\label{HMmps}
\|\M_{\G,m}(f)\|_p,\,\|\H_{\G,m}(f)\|_p\lesssim_{\g,p}
\left\|\left(\sum_{j\in\Z} \left|\L_{\G,j,m} f\right|^2 \right)^{\frac{1}{2}}\right\|_{p}\:.
\eeq

\noindent\textbf{Claim.}
\textit{\underline{Assume} that, given\footnote{Here $L^p(\R^2)$ is the standard $L^p$-space relative to the Lebesgue measure on $\R^2$.} any $\{f_j\}_{j\in \Z}\subset L^r(\R^2)$, the following relation holds:
\beq\label{HMmps1}
\left\|\sup_{j\in\Z} |\L_{\G,j,m} f_j|\right\|_{r}\lesssim_{\g,r} m^2\,\|\sup_{j\in \Z}|f_j|\|_r\,,\:\:\:\:\:\:\forall\:1<r\leq 2\:.
\eeq
Then, for any $1<p\leq 2$, one has
\beq\label{HMmps2}
\left\|\left(\sum_{j\in\Z} \left|\L_{\G,j,m} f\right|^2 \right)^{\frac{1}{2}}\right\|_{p}\lesssim_{\g,p}m^2\,\|f\|_p\:.
\eeq}
Indeed, to see this, we apply the following vector valued interpolation result: assume we are given a sequence of linear operators $\textbf{T}:=\{T_{j}\}_{j\in\Z}$ such that for any $1<p\leq 2$ and any $j\in\Z$ the operator $T_{j}$ is well defined on $L^p(\R^2)$. For $\textbf{f}:=\{f_j\}_{j\in\Z}\subset L^p(\R^2)$ we set
\beq\label{convv}
\textbf{T}(\textbf{f}):=\{T_{j}(f_j)\}_{j\in\Z}\:.
\eeq
With the above conventions, we have that the following holds

\begin{p1}\label{Interp}
Assume that we are given the following:
\begin{itemize}
\item there exists $1<r\leq 2$ such that
\beq\label{HMmps3}
\textbf{T}\,:L^{r}(l^{\infty}(\Z))\,\rightarrow\,L^{r}(l^{\infty}(\Z))\:,
\eeq
with ${\|\textbf{T}\|}_{L^{r}(l^{\infty}(\Z))\,\rightarrow\,L^{r}(l^{\infty}(\Z))}\leq A_r<\infty$.
\item there exists $1<q\leq 2$ such that
\beq\label{HMmps4}
\textbf{T}\,:L^{q}(l^{q}(\Z))\,\rightarrow\,L^{q}(l^{q}(\Z))\:,
\eeq
with ${\|\textbf{T}\|}_{L^{q}(l^{q}(\Z))\,\rightarrow\,L^{q}(l^{q}(\Z))}\leq B_q<\infty$.
\end{itemize}
Then, applying vector valued interpolation, (see \cite{BLo}), one has that
\beq\label{HMmps5}
\textbf{T}\,:L^{p}(l^{2}(\Z))\,\rightarrow\,L^{p}(l^{2}(\Z))\:,
\eeq
with ${\|\textbf{T}\|}_{L^{p}(l^{2}(\Z))\,\rightarrow\,L^{p}(l^{2}(\Z))}\leq C(r,q,A_r, B_q)<\infty$ and $\frac{1}{p}=\frac{1}{2}+(1-\frac{q}{2})\,\frac{1}{r}$.
\end{p1}

Take now $T_j:=\L_{\G,j,m}$. Then, we can place ourselves in the settings offered by our proposition. Indeed, we notice that \eqref{HMmps1} implies  both \eqref{HMmps3} and \eqref{HMmps4} for \textit{any} $1<r\leq2$ and $1<q\leq 2$ respectively,
with $A_r\lesssim_{\g,r} m^2$ and $B_q\lesssim_{\g,q} m^2$.

Conclude thus now based on \eqref{HMmps3}-\eqref{HMmps5} that \eqref{HMmps2} holds. Moreover interpolating \eqref{HMmps2} with the $p=2$ case given by \eqref{HMm2} we have that \eqref{HMlp} follows.

\begin{o0}\label{lpl2}
1) It is worth noticing that the assumption \eqref{HMmps1} holds. Indeed, by following similar steps with those described in either Section \ref{LPS}, Stage 2 or Section \ref{LPg2} one can show that
\beq\label{physss}
\sup_{j\in\Z}|\L_{\G,j,m}(f_j)(x,y)|\lesssim_{\g}  \frac{1}{2^m}\sum_{l=2^{m-2}}^{2^{m+2}}\,M_1^{[l]}
 M_2^{[2^m]}(\sup_{j\in \Z}|f_j|)(x,y)\,.
\eeq
However, \eqref{physss} makes the interpolation argument above superfluous, since \eqref{physss} trivially implies Theorem \ref{hmdiagp} via Proposition \ref{maxestim}.

It would be interesting to identify a distinct method of proving \eqref{HMmps1} that does not need pointwise control over the LHS in \eqref{physss}. In connection with this last comment we propose
\end{o0}

2) \textbf{Open Question.} \textit{Does there exist a (less general) theoretical framework that applies to the case  $T_j:=\L_{\G,j,m}$ and such that Proposition \ref{Interp} above can be strengthened by replacing \eqref{HMmps3} with the milder assumption that there exists some $1<r\leq 2$ such that
\beq\label{HMmps10}
\left\|\sup_{j\in\Z} |T_j(f)|\right\|_{r}\lesssim_{r} A_r\,\|f\|_r\,?
\eeq}

This strengthening is in fact possible, for example, if each $T_j$ is of the form $K_j*f$ with $K_j$ a \textit{positive} kernel. An interesting question is whether this remains true for a sufficiently general class of operators\footnote{It is highly unlikely for such a strengthened analogue to apply in full generality.}  $T_j$ that would include the operators treated in this paper, \textit{i.e.} $T_j:=\L_{\G,j,m}$.

Assume for the moment that our question has an affirmative answer (AA), that is if \eqref{HMmps10} and \eqref{HMmps4} hold then \eqref{HMmps5} holds. Then, as we will see through a bootstrapping argument, taking $T_j:=\L_{\G,j,m}$, condition \eqref{HMmps10} can be further replaced by the estimate\footnote{For a proof of this, see \eqref{physs11} below.}
\beq\label{HMmps11}
\left\|\L_{\G,j,m}(f)\right\|_{r}\lesssim_{r} A_r\,\|f\|_r\,,\:\:\:\:\:\:\forall\:1<r<\infty\:;
\eeq
in other words, if we have (AA) only assuming \eqref{HMmps10}, then the same affirmative answer holds if instead of \eqref{HMmps10} we only have \eqref{HMmps11}.

Indeed, \eqref{HMmps11} trivially implies \eqref{HMmps4} for all $1<q\leq 2$ (and in fact for all $1<q<\infty$) and also
\beq\label{HMmps12}
\|(\sum_{j}|\L_{\G,j,m}(f)|^2)^{\frac{1}{2}}\|_2\lesssim \|f\|_2\;.
\eeq
Now \eqref{HMmps12} further implies \eqref{HMmps10} for $r=2$.

Based on our assumption provided by (AA) for $r=2$, we deduce that \eqref{HMmps5} holds for any $\frac{1}{2}\leq \frac{1}{p}<\frac{3}{4}$. This, however, implies in turn that \eqref{HMmps10} holds for any $2\geq r>\frac{4}{3}$, which by another application of (AA) further implies that \eqref{HMmps5} holds for any $\frac{1}{2}\leq\frac{1}{p}<\frac{7}{8}$. By continuing this bootstrapping algorithm we conclude that \eqref{HMmps5} and hence \eqref{HMmps2} holds in the full range $1<p\leq 2$.

We end this comment by presenting the short proof of \eqref{HMmps11}, as promised above. Indeed, by applying
Jensen's inequality, one immediately has
\beq\label{physs11}
\eeq
$$\left\|\L_{\G,j,m} f\right\|_{r}^r\lesssim $$
$$\int_{\R}\int_{\R^2} \left|(f*_{x}\check{\phi}_{m+j}*_{y}\check{\phi}_{m,j,x})\left(x-\frac{t}{2^j},\,y+\gamma_x(\frac{t}{2^j})\right)\right|^r
\,|\rho(t)| dx\,dy\,dt$$
$$\lesssim \int_{\R}\int_{\R^2} \left|(f*_{x}\check{\phi}_{m+j}*_{y}\check{\phi}_{m,j,x})\left(x,\,y)\right)\right|^r
\,|\rho(t)| dx\,dy\,dt$$
$$\lesssim \int_{\R^2} |M_1 M_2f(x,y)|^r\,dx\,dy\lesssim_{r} \int_{\R^2}|f(x,y)|^r\,dx\,dy\:.$$

\subsection{Other approaches (I.4) - the case $p=2$}

In this section, we want to briefly mention a different approach to the $L^2-$case that is explored by the author in \cite{lvBLC}.  The starting point for both our present approach and the one in \cite{lvBLC} is the formula \eqref{Bjmintdec}. From this on, in order to better understand the main differences and subtleties, we proceed in an antithetical manner:
\begin{itemize}
\item the approach developed in Section \ref{ldecsec} here focuses on \emph{linearization}, that is the discretization of our phase is performed such that the localized pieces oscillate at the linear level. This is revealed by 1) the choice of the phase-localization \eqref{xloc} reducing \eqref{Bjmintdec} to \eqref{hjkmnd} and 2) the ``cuttings" \eqref{xiloc} and \eqref{etaloc} that further triggers the Gabor frame decomposition \eqref{Gs}.

\item in contrast with this, the approach in \cite{lvBLC} highlights a \emph{non-linear} analysis, that is the discretization of the phase allows the resulting pieces to oscillate up to the second order. This is reflected by 1) the specific choice of the Gabor frame decomposition whose wave-length oscillation\footnote{Relative to the $m-$parameter.} is $2^m$ as opposed to $2^{\frac{m}{2}}$ in the first case above; 2) the two required consecutive applications of \emph{stationary phase} principle as opposed to time-frequency correlation at linear level based on integration by parts and application of \emph{non-stationary} phase principle in the first case.
\end{itemize}

Due to space constrains, we defer a more extended discussion to \cite{lvBLC}. Here, we will only limit ourselves to say that this second approach developed in \cite{lvBLC} was verified only for the \emph{tensor-product} situation $\g(x,t)=u(x)\,\g(t)$ where $u(\cdot)$ is measurable and $\g\in\n\f$ is a suitable ``non-flat" curve as introduced in \cite{lv4}. It seems however that this second approach has more limitations than the first one detailed in the present paper.

\section{Proof of the Main Theorem, Part (II)}

In this section, based on the reasonings involved for proving Part (I) and on the heuristic presented in Section \ref{HCa}, we will provide the proof of Part (II) in our Main Theorem.

Assume thus throughout this section that $\g\in \mathbf{M}_x\mathbf{NF}_{t}$. Recalling \eqref{defph} and taking now in \eqref{hgj}
\beq\label{bdHVTpm0}
f(x,y)=h(x,y):= f(x)\, (g*\check{\phi})(y)
\eeq
one deduces that
\beq\label{hgj1}
H_{\G} h(x,y) = \int_{\R}\hat{g}(\eta)\,\hat{\phi}(\eta)\,e^{-i\,\eta\,y}\, C_{\g}^{\eta} f(x)\,d\eta\,,
\eeq
where
\beq\label{carlgj1}
C_{\g}^{\eta} f(x):=\sum_{j\in\Z} C_{\g, j}^{\eta} f(x):= \sum_{j\in\Z}\int_{\R}f(x-t)\, e^{i\,\eta\,\g(x,t)}\,\r_j(t)\,dt\,.
\eeq

Using now the fact that the $\eta$ variable stays within the support of the function $\phi$, we rewrite \eqref{firstl6} in the form
\beq\label{firstl6s}
  \phi(\eta)\,\mathfrak{m}_{j}(x,\xi,\eta)\approx\sum_{m,n\in\Z}  \mathfrak{m}_{j,m,n,0}(x,\xi,\eta)\:,
\eeq
and thus we immediately deduce that
\beq\label{carlgj11}
\eeq
$$\phi(\eta)\,C_{\g, j,m,n}^{\eta} f(x):= \int_{\R} \hat{f}(\xi)\,\mathfrak{m}_{j,m,n,0}(x,\xi,\eta)\,e^{i x \xi}\,d\xi$$
$$=\phi(\eta)\,\phi(\frac{\g'_x(2^{-j})}{2^{n+j}})\,\int_{\R}(f*\check{\phi}_{m+j})(x-t)\, e^{i\,\eta\,\g(x,t)}\,\r_j(t)\,dt\,.$$

As a consequence, with the obvious correspondences, we are following the same line of thought as in Section \ref{anmult}, see \eqref{LF} - \eqref{firstl8sss}, and decompose
\beq\label{decC}
C_{\g}^{\eta} = C_{\g}^{L,\eta} + C_{\g}^{H\not\Delta,\eta}+ C_{\g}^{H\Delta,\eta}\,.
\eeq

\subsection{The low-frequency case.} In this section we will prove the analogue of Theorem \ref{lowf}:

\begin{p1}\label{lowfpcc} With the above notations we have
\beq\label{mlc1}
| C_{\g}^{L,\eta}f(x)|\lesssim_{\g} M H f(x)\,+\, M f(x)\:,
\eeq
uniformly in the parameter $|\eta|\lesssim 1$ and thus, for any $1<p<\infty$, one has
\beq\label{mhl1p0}
\|\sup_{|\eta|\lesssim 1} |C_{\g}^{L,\eta}f| \|_{p}\lesssim_{\g,p} \|f\|_p\:.
\eeq
\end{p1}

\begin{proof}
Applying \eqref{mlre}, one deduces that
\beq\label{mlrew}
\eeq
$$|C_{\g}^{L,\eta} f(x)|\lesssim$$
$$\left|\sum_{{p+\ell>0}\atop{p,\,\ell\in\N}} \sum_{j\in\Z}\chi_{\Z^{x}}(j) \sum_{(m,n)\in (\Z_-)^2} \frac{i^{p+\ell}(-1)^\ell}{p!\ell !}\,2^{m\,l}\,2^{n\,p}\,C_{p,\ell,j,x} \,(f*\check{\tilde{\phi}}_{l,m+j})(x)\,\tilde{\phi}_{p}(\frac{\g'_x(2^{-j})}{2^{n+j}})\,\tilde{\phi}_{p,0}(\eta)\right|$$
As for $\H_{\G}^{L}$, due to the fast coefficient decay, it is enough to discuss only the two extreme cases

$\newline$
\noindent \textbf{Case 1} $l=0$, $p=1$,  $m\in \Z_{-}$ and $n=-1$
$\newline$

In this situation, we have

$$|C_{\g}^{L,1,\eta} f(x)|:=
\left|\sum_{j\in\Z} \chi_{\Z^{x}}(j)\sum_{m\in \Z_-} C_{1,0,j,x} \,(f*\check{\tilde{\phi}}_{0,m+j})(x)\,\tilde{\phi}_{1}(\frac{\g'_x(2^{-j})}{2^{j-1}})\,\tilde{\phi}_{1,0}(\eta)\right|$$
$$\lesssim M f(x)\,\sum_{j\in\Z} \left| C_{1,0,j,x} \,\tilde{\phi}_{1}(\frac{\g'_x(2^{-j})}{2^{j-1}})\,\tilde{\phi}_{1,0}(\eta)\right|\lesssim Mf(x)\:.$$

$\newline$
\noindent \textbf{Case 2} $l=1$, $p=0$,  $m=-1$ and $n\in \Z_{-}$
$\newline$

In this second case, we recall that $C_{0,1,j,x}=C$ and thus we focus on
$$C_{\g}^{L,2,\eta} f(x):=
\sum_{j\in\Z} \chi_{\Z^{x}}(j)\sum_{n\in \Z_{-}} \,(f*\check{\tilde{\phi}}_{1,j-1})(x)\,\phi(\frac{\g'_x(2^{-j})}{2^{n+j}})\,\phi(\eta)$$
$$=\phi(\eta)\,\sum_{j\in\Z} \chi_{\Z^{x}}(j)\,(f*\check{\tilde{\phi}}_{1,j-1})(x)\,\psi(\frac{\g'_x(2^{-j})}{2^{j}})\:.$$

Setting now  $\S(\tilde{\phi})_{j}:=\sum_{l=-\infty}^{j} \check{\tilde{\phi}}_{1,l-1}$ and $\psi_{\g,j}(x):=  \psi(\frac{\g'_x(2^{-j})}{2^{j}})$ and letting\footnote{Here we let $a^x,\,b^x\in \Z\cup\{\pm\infty\}$ with $a^x\leq b^x$.} $\Z^{x}=[a^x,\,b^x]\cap \Z$ we proceed as in the proof of Lemma \ref{trm22p} and apply an Abel summation argument

$$\sum_{j=a^x}^{b^x}\psi_{\g,j}(x)\, (f*\check{\tilde{\phi}}_{1,j-1})(x)
=\sum_{j=a^x}^{b^x}\psi_{\g,j}(x)\, (f*(S(\tilde{\phi})_{j}-S(\tilde{\phi})_{j-1}))(x)$$
$$= \psi_{\g,b^x}(x)\, (f*(S(\tilde{\phi})_{b^x}))(x)\,- \,\psi_{\g,a^x}(x)\, (f*(S(\tilde{\phi})_{a^x-1}))(x)$$
$$+\,\sum_{j=a^x}^{b^x-1}(\psi_{\g,j}(x)-\psi_{\g,j+1}(x))\, (f*(S(\tilde{\phi})_{j})(x)$$
From \eqref{gamprop} and the above chain of equalities, we immediately have that
$$|C_{\g}^{L,2,\eta} f(x)|\lesssim MHf(x)+Mf(x)\:,$$
thus finishing our proof.
\end{proof}

\subsection{The high-frequency far from diagonal case.} In this section we will prove the analogue of Theorem \ref{hofdiag}:

\begin{p1}\label{lowfpccs} Let $1<p<\infty$. Then, with the above notations and conventions, we have
\beq\label{mhl1hndp}
\|\sup_{|\eta|\lesssim 1} |C_{\g}^{H\not\Delta,\eta}f| \|_{p}\lesssim_{\g,p} \|f\|_p\:.
\eeq
\end{p1}

\begin{proof}
We start by reminding that from now on we will identify the two multipliers $\mathfrak{m}_{j,m,n,0}$ and $\underline{\mathfrak{m}}_{j,m,n,0}$ since the mean zero condition of $\r_j$ will no longer play any role.

Setting (recalling) now
$$\mathfrak{m}_{m,n}(x,\xi,\eta)=\sum_{j\in\Z} \mathfrak{m}_{j,m,n,0}(x,\xi,\eta)\,,$$
$$\mathfrak{m}^{j,H\not\Delta}(x,\xi,\eta)=\sum_{(m,n)\in \Z^2\setminus ((\Z_-)^2 \cup \Delta)} \mathfrak{m}_{j,m,n,0}(x,\xi,\eta)\,,$$
and
$$\mathfrak{m}^{H\not\Delta}(x,\xi,\eta)=\sum_{j\in\Z} \mathfrak{m}^{j,H\not\Delta}(x,\xi,\eta)\,,$$
we proceed exactly as in the proof of Theorem \ref{hofdiag} and split the discussion in two cases: (I) $m>|n|+C(\g)$ and (II) $n>|m|+C(\g)$. Focusing on the first case - the second one has a similar treating - one applies \eqref{mfdmex} in order to deduce
\beq\label{mfdmexx}
\eeq
 $$\mathfrak{m}_{m,n}(x,\xi,\eta)= $$
$$\frac{1}{2^m}\,\sum_{j\in\Z} \chi_{\Z^{x}}(j)\,\left(\int_{\R} e^{-i\, \frac{\xi}{2^j}\, t}\: e^{i \,\eta \g_x(\frac{t}{2^j})}\,\r_j(x,t)\,dt\right)\,
\phi\left(\frac{\xi}{2^{m+j}}\right)\,\phi(\frac{\g'_x(2^{-j})}{2^{n+j}})\,\phi(\eta)\,.$$
From this, with the obvious correspondences, we further have
\beq\label{mlreww}
|C_{\g}^{H\not\Delta,\eta} f(x)|\approx\left|\phi(\eta)\,\sum_{j\in\Z}\chi_{\Z^{x}}(j) \sum_{{(m,n)\in \Z^2\setminus ((\Z_-)^2 \cup \Delta)}\atop{m>|n|+C(\g)}}\frac{1}{2^m}\,C_{\g, j,m,n}^{\eta} f(x)\right|\,.
\eeq
Since taking a supremum in $\eta$ would have as an effect - after the linearization - the modification of $\eta$ into $\eta(x)$ with the latter function being absorbed into the expression $\g(x,t)$ it is enough thus to estimate our expressions for say $\eta=1$. Taking now $1<p<\infty$ and applying a standard square argument, we have thus that
\beq\label{mclreww0}
\eeq
$$\|C_{\g}^{H\not\Delta,\eta} f\|_p$$
$$\lesssim_{p, \g}
\sum_{(m,n)\in \Z^2\setminus ((\Z_-)^2 \cup \Delta)}\,\frac{1}{2^{\max\{m,n\}}}\,
\left\|\left(\sum_{j\in\Z}\chi_{\Z^{x}}(j) |C_{\g, j,m,n}^{1} f(x)|^2\right)^{\frac{1}{2}}\right\|_p$$
$$\lesssim \sum_{(m,n)\in \Z^2\setminus ((\Z_-)^2 \cup \Delta)}\,\frac{1}{2^{\max\{m,n\}}}\,
\left\|\left(\sum_{j\in\Z} |M (f*\check{\phi}_{m+j})|^2\right)^{\frac{1}{2}}\right\|_p\lesssim_{p} \|f\|_p\:.$$

\end{proof}

\subsection{The high-frequency close to diagonal case.} Finally, in this section we are proving the analogue of Theorem \ref{Diagp}:

\begin{p1}\label{diagpp} Let $1<p<\infty$. Then, we have
\beq\label{mhl1hndp}
\|\sup_{|\eta|\lesssim 1} |C_{\g}^{H\Delta,\eta}f| \|_{p}\lesssim_{\g,p} \|f\|_p\:.
\eeq
\end{p1}

\begin{proof}
In our setting, we have
$$C_{\g}^{H\Delta,\eta}= \sum_{m\in\N} C_{m,\g}^{H\Delta,\eta}\,,$$
with
$$C_{m,\g}^{H\Delta,\eta}\approx\sum_{j\in\Z} C_{\g,j,m,m}^{\eta}\,.$$

Notice now that for any $m\in\N$ and $j\in\Z$ we trivially have
\beq\label{keyq}
|C_{\g,j,m,m}^{\eta}f|(x)\lesssim
M (f*\check{\phi}_{m+j})(x)\,.
\eeq
Moreover, following the same argument as in \eqref{mclreww0}, we now get - uniformly in $m\in\N$ - that
\beq\label{mclreww1}
\|\sup_{|\eta|\lesssim 1}\,C_{m,\g}^{H\Delta,\eta} f\|_p
\lesssim_{\g,p}\left\|\left(\sum_{j\in\Z}|M (f*\check{\phi}_{m+j})|^2\right)^{\frac{1}{2}}\right\|_p\lesssim_{p} \|f\|_p\:.
\eeq

Now, our Proposition \ref{diagpp} follows from \eqref{mclreww1}, the result below and a standard interpolation argument

\begin{p1}\label{diagpp2} There exists $\bar{\ep}>0$, depending only on $\g$, such that, for any $m\in\N$ and $j\in\Z$, the following holds:
\beq\label{mhl1hndp2}
\|\sup_{|\eta|\lesssim 1} |C_{\g,j,m,m}^{\eta}f| \|_{2}\lesssim_{\g} 2^{- \bar{\ep}\,m}\, \|f\|_2\:.
\eeq
\end{p1}

This last proposition however is a direct application of Theorem \ref{hmdiag2} together with the following observation

\beq\label{hgj11}
\H_{\G,m}^{H\Delta}(h)(x,y)=\H_{\G,m}(h)(x,y)= \int_{\R}\hat{g}(\eta)\,\hat{\phi}(\eta)\,e^{-i\,\eta\,y}\, C_{m,\g}^{H\Delta,\eta} f(x)\,d\eta\,,
\eeq
from which we deduce after an application of Parseval that
\beq\label{hgj112}
\|\sup_{|\eta|\lesssim 1} |C_{m,\g}^{H\Delta,\eta}|\|_{L^2(\R)\rightarrow L^2(\R)}\leq \|\H_{\G,m}\|_{L^2(\R^2)\rightarrow L^2(\R^2)} \,.
\eeq
\end{proof}

\section{The richness of the class $\mathbf{M}_x\mathbf{NF}_{t}$ - Proof of Theorem \ref{Gencurvpolyn}.}\label{richclass}

\noindent\textbf{Theorem \ref{Gencurvpolyn} - Restatement.}\label{polynmc} \textit{Let $d\in\N$, $d\geq 1$ and assume we are given $\{\a_k\}_{k=1}^{d}\subset\R\setminus\{0,\,1\}$. Define
\beq\label{py0}
\p_d:=\left\{\g(x,t)\,|\,\g(x,t)=\sum_{k=1}^d a_k(x)\,t^{\a_k}\:\:\textrm{with}\:\:\{a_k(\cdot)\}_{k=1}^d\:\textrm{real measurable}\right\}\:.
\eeq
Then
\beq\label{py1}
\p_d\subset \mathbf{M}_x\mathbf{NF}_{t}\:.
\eeq}

\begin{proof}

$\newline$
\noindent\textbf{Step 1.} Verifying \eqref{u}-\eqref{variation0}.
$\newline$

Throughout our discussion we fix $\g(x,t)\in \p_d$ and $x\in\R$ and assume wlog that $\{\a_k\}_{k=1}^d\subset\R\setminus\{0,\,1\}$ arranged in a strictly increasing order.

We start by defining
$$\R_0:=\{x\in\R\,|\,a_1(x)=a_2(x)=\ldots=a_d(x)=0\}\,,$$
and set
$\R_1:=\R\setminus \R_0\,.$

Next, based on some elementary reasonings, we can always reduce our discussion to the situation $a_k(x)\not=0$ for all $x\in\R_1$ and $1\leq k\leq d$. Given this, we can
simply set $A=1$. Also, in order to avoid the discussion about the definition of $t^{\a_k}$ for $t<0$, wlog we can assume throughout this section that $t>0$.

Now, in order to be able to define the sets $\{\Z^{x}_{\b}\}_{\b}$, we proceed as follows:

Firstly, we notice that
\beq\label{defjkd0}
\eeq
$$\g(x,t)=\sum_{k=1}^d a_k(x)\,t^{\a_k}\,,$$
$$\frac{d}{dt}\,\g(x,t):=\g'(x,t)=\sum_{k=1}^{d} \a_k a_k(x)\,t^{\a_k-1}\:,$$
and
$$\frac{d^2}{dt^2}\,\g(x,t):=\g''(x,t)=\sum_{k=1}^{d} \a_k\,(\a_k-1)\,a_k(x)\,t^{\a_k-2}\:.$$

Secondly, we fix $x\in\R_1$ and for $\mathfrak{K}>0$ properly chosen\footnote{For example $\mathfrak{K}= C\,d$ with $C>>1$ absolute constant works.} and each pair $(n,k)\in\N^2$ with $d\geq n>k\geq 1$, define\footnote{Throughout the paper $R_{+}^{*}$ stands for $(0,\infty)$.}
\beq\label{defjk}
\eeq
$$V_{nk}(x):=\left\{t\in\R_{+}^{*}\,|\,
\frac{1}{\mathfrak{K}}\,\min\left\{1,\,\frac{|\a_k|}{|\a_n|},
\,\frac{|\a_k\,(\a_k-1)|}{|\a_n\,(\a_n-1)|}\right\}\,\frac{|a_k(x)|}{|a_n(x)|}\leq |t|^{\a_n-\a_k}\right\}$$
$$\bigcap\left\{t\in\R_{+}^{*}\,|\,|t|^{\a_n-\a_k}\leq
\mathfrak{K}\,\max\left\{1,\,\frac{|\a_k|}{|\a_n|},
\,\frac{|\a_k\,(\a_k-1)|}{|\a_n\,(\a_n-1)|}\right\}\,\frac{|a_k(x)|}{|a_n(x)|}\right\}\:.$$

Thirdly, we notice that if $t\in V_{nk}^{c}(x):=\R_{+}^{*}\setminus V_{nk}(x)$, then we are in one of the two situations below:
\beq\label{defjk1}
\eeq
\noindent - either
$|t|^{\a_n}\,|a_n(x)|>\mathfrak{K}\,|t|^{\a_k}\,|a_k(x)|$ and $|\a_n|\,|t|^{\a_n}\,|a_n(x)|>\mathfrak{K}\,|\a_k|\,|t|^{\a_k}\,|a_k(x)|$ and
 $|\a_n|\,|\a_n-1|\,|t|^{\a_n}\,|a_n(x)|>\mathfrak{K}\,|\a_k|\,|\a_k-1|\,|t|^{\a_k}\,|a_k(x)|$;

\noindent - or $|t|^{\a_n}\,|a_n(x)|<\mathfrak{K}^{-1}\,|t|^{\a_k}\,|a_k(x)|$ and $|\a_n|\,|t|^{\a_n}\,|a_n(x)|<\mathfrak{K}^{-1}\,|\a_k|\,|t|^{\a_k}\,|a_k(x)|$ and
 $|\a_n|\,|\a_n-1|\,|t|^{\a_n}\,|a_n(x)|<\mathfrak{K}^{-1}\,|\a_k|\,|\a_k-1|\,|t|^{\a_k}\,|a_k(x)|$.

Deduce from this that
\beq\label{defjk2}
\bigcap_{1\leq k<n\leq d} V_{nk}^{c}(x)= \bigcup_{\b=1}^{O(d^2)} J_{\b}^{x}\,,
\eeq
with each $J_\b^{x}$ interval so that there exist $k_\b\in\{1,\ldots,\,d\}$ and such that for every $t\in J_{\b}^{x}$ one has
\beq\label{defjk3}
\eeq
\begin{itemize}
\item $|t|^{\a_{k_\b}}\,|a_{k_\b}(x)|>\mathfrak{K}\,\max_{k\not=k_\b}\{|t|^{\a_k}\,|a_k(x)|\}$ and
\item $|\a_{k_\b}|\,|t|^{\a_{k_\b}}\,|a_{k_\b}(x)|>\mathfrak{K}\,\max_{k\not=k_\b}\{|\a_k|\,|t|^{\a_k}\,|a_k(x)|\}$ and
\item $|\a_{k_\b}|\,|\a_{k_\b}-1|\,|t|^{\a_{k_\b}}\,|a_{k_\b}(x)|>\mathfrak{K}\,\max_{k\not=k_\b}\{|\a_k|\,|\a_k-1|\,|t|^{\a_k}\,|a_k(x)|\}$.
\end{itemize}
Thus, taking now $J_{0}^{x}:=\R_{+}^{*}\setminus \bigcap_{1\leq k<n\leq d} V_{nk}^{c}(x)$ and setting $\Z_{\b}^{x}:=\log_2 (J_{\b}^{x})\cap \Z$, we deduce that
\beq\label{defjk31}
\#\Z_{0}^{x}\leq \sum_{1\leq k<n\leq d}\frac{100}{\a_n-\a_k}\,\left(\log \mathfrak{K} + \left|\log \frac{|\a_k|}{|\a_n|}\right|\, + \, \left|\log \frac{|\a_k-1|}{|\a_n-1|}\right| \right)\,,
\eeq
and hence \eqref{contr0} holds.

Next, for $1\leq \b\leq B=O(d^2)$, we obviously have $\g(x,\cdot)\in C^{\infty}(J_{\b}^{x})$ and for $t\in J_{\b}^{x}$ one has
\beq\label{defjk4}
\g''(x,t)= \a_{k_\b}\,(\a_{k_\b}-1)\,t^{\a_{k_\b}-2}\,a_{k_\b}(x)\,(1\,+\, O(d\,\mathfrak{K}^{-1}))\,,
\eeq
and hence \eqref{deriv0} holds as well.

Finally, recalling that $\Z^{x}=\bigcup_{\b=1}^{B} \Z_{\b}^{x}$, we have
\beq\label{defjkd1}
\eeq
$$\sup_{x\in\R_{1}}\sup_{c\in\R_{+}}\#\{j\in\Z^{x}\,|\,|2^{-j}\,\g'_x(2^{-j})|\in[c,2c]\}\leq$$
$$\sum_{\b=1}^{B}\sup_{x\in\R_{1}}\sup_{c\in\R_{+}}\#\{j\in\Z_{\b}^{x}\,|\,
|\a_{k_\b}\,(2^{-j})^{\a_{k_\b}}\,a_{k_\b}(x)|\,(1\,+\, O(d\,\mathfrak{K}^{-1}))\in[c,2c]\}$$
$$\lesssim d +\sum_{k=1}^d\frac{1}{|\a_k|}=N(\g)\:.$$

$\newline$
\noindent\textbf{Step 2.} Verifying \eqref{asymptotic0}-\eqref{fstterm0}.
$\newline$

Fix $x\in\R_{1}$ and $1\leq \b\leq B$. For $\frac{1}{4}\leq|t|\leq 4$ and $j\in \Z^{x}_{\b}$ we have
\beq\label{asymptotic01}
Q_{j}(x,t)=\frac{\g_x(2^{-j}\,t)}{2^{-j}\,\g'_x(2^{-j})}
=\frac{t^{\a_{k_\b}}}{\a_{k_\b}}\,(1\,+\, O(d\,\mathfrak{K}^{-1}))\,,
\eeq
and thus both \eqref{asymptotic0} and \eqref{fstterma0} are satisfied.

Similarly, we have
\beq\label{asymptotic011}
Q^{''}_{j}(x,t)=\frac{2^{-j}\,\g_x^{''}(2^{-j}\,t)}{\g'_x(2^{-j})}
=(\a_{k_\b}-1)\,t^{\a_{k_\b}-2}\,(1\,+\, O(d\,\mathfrak{K}^{-1}))\,,
\eeq
which implies that \eqref{fstterm0} holds.

$\newline$
\noindent\textbf{Step 3.} Verifying \eqref{ndeg0}.

In what follows, without loss of generality we fix $\a=1$, and $1=\b=O(d^2)$. Next, following the same reasonings and notations as in the proof of Lemma \ref{nondegcont}, we see that our task is to show that there exists $\bar{\ep}\in (0,1)$ such that \eqref{ndeg00} holds.  Now, linearizing the supremum in $w$, we would like to show that, for any given map $l\,\rightarrow\,w_l$, the following holds\footnote{Notice that \eqref{ndeg00} is equivalent with \eqref{ndeg00e}.}
\beq\label{ndeg00e}
I_{m}:=2^{-\frac{m}{2}}\,\sum_{l\in\mathcal{R}_m}\,\left(\int_{-\frac{1}{2}}^{\frac{1}{2}}\frac{1}
{\left( \lfloor 2^{\frac{m}{2}}\,\left(\tilde{q}(x,\,l\,+\,x)-w_l\right)\rfloor\right)^2}\,dx\right)\lesssim_{\g} 2^{-2\,\bar{\ep}\,m}\:.
\eeq
Fix from now on the map $l\,\rightarrow\,w_l$. Next, for a suitable parameter $\d\in (0,1)$ chosen later, recalling the defined concepts in \eqref{A} - \eqref{H}, we have
\beq\label{ndeg00e1}
I_{m}:=2^{-\frac{m}{2}}\,\sum_{{l\in\mathcal{R}_m}\atop{(l,w_l)\in \L_{\tilde{q},\d}}}\,\left(\ldots\right)
\,+\,2^{-\frac{m}{2}}\,\sum_{{l\in\mathcal{R}_m}\atop{(l,w_l)\in \H_{\tilde{q},\d}}}\,\left(\ldots\right)
=:I_{m}^{L,\d}\,+\,I_{m}^{H,\d}\:.
\eeq
We now have
\beq\label{ndeg00e2}
\eeq
$$I_{m}^{L,\d}=2^{-\frac{m}{2}}\,\sum_{{l\in\mathcal{R}_m}\atop{(l,w_l)\in \L_{\tilde{q},\d}}}\,
\left(\int_{A_{\tilde{q},\d}(l,w_l)} (\ldots)\,+\,\int_{[-\frac{1}{2},\,\frac{1}{2}]\setminus A_{\tilde{q},\d}(l,w_l)} (\ldots)\right)$$
$$\lesssim  2^{-\frac{m}{2}}\,\sum_{{l\in\mathcal{R}_m}\atop{(l,w_l)\in \L_{\tilde{q},\d}}} |A_{\tilde{q},\d}(l,w_l)|
+2^{-4\d m}\lesssim 2^{-2\,\d\,m}\:,$$
while
\beq\label{ndeg00e3}
I_{m}^{H,\d}=2^{-\frac{m}{2}}\,\sum_{{l\in\mathcal{R}_m}\atop{(l,w_l)\in \H_{\tilde{q},\d}}}\, 1\lesssim
 2^{-\frac{m}{2}} \, (\# \H_{\tilde{q},\d})\:.
\eeq

For $X\subseteq\R$ measurable, we now define
\beq\label{PY1}
\P_d(X\times I):=\left\{q(x,t)=\sum_{k=1}^d a_k(x)\,t^{\a_k}\,\big|\, \begin{array}{ll}
      \{a_k(\cdot)\}_{k=1}^d\:\textrm{measurable},\\\|a_k\|_{L^{\infty}(X)}\lesssim 1\:\:\forall\:1\leq k\leq d\,,\\
      \inf_{x\in X}|\frac{\partial}{\partial t}q(x,t)|\gtrsim_{\{\a_j\}_{j=1}^d,\, d} 1\,.
    \end{array}
\right\}\,,
\eeq
where here $d\in\N$, $\a_1<\ldots <\a_d$ with $\{\a_j\}_{j=1}^d\subset\R\setminus\{0\}$.

Our goal is to show the following

\begin{p1}\label{fewhe} Fix $\tilde{q}\in \P_d([-\frac{1}{2},\,\frac{1}{2}]\times I)$. Then, with the previous notations and for a proper choice of $\d$, there exists $\nu=\nu(\d, d)\in (0,\frac{1}{2})$ such that
\beq\label{fhe0s}
\# \H_{\tilde{q},\d}\lesssim 2^{(\frac{1}{2}-\nu) m}\:.
\eeq
\end{p1}

Notice that if we believe for a moment Proposition \ref{fewhe} above, then from \eqref{ndeg00e1} - \eqref{ndeg00e3}, we deduce that
\beq\label{ndeg00ee}
I_{m}\lesssim 2^{-2\,\d\,m}\,+\,2^{-\nu\, m}\:,
\eeq
which proves \eqref{ndeg00e} for $\bar{\ep}=\min \{\d,\,\frac{\nu}{2}\}$.

$\newline$
\noindent\textsf{The proof of Proposition \ref{fewhe}.}
$\newline$

In what follows we will need the following

\begin{l1}\label{largeint} Let $n,\,N,\,M\in\N$ with $n,\,M\leq N$. Assume we are given  $\{I_l\}_{l=1}^{N}$ sets such that for any $1\leq l\leq N$ the following properties hold:
\beq\label{propset}
\eeq
\begin{itemize}
\item $I_l\subseteq [-\frac{1}{2},\,\frac{1}{2}]$;

\item $|I_l|\geq M^{-1}$.
\end{itemize}
Then, if
\beq\label{n1}
N\geq\,2\,M^{n}\,n^n\,,
\eeq
there exists a subset $S\subset\{1,\ldots, N\}$ such that
\beq\label{propS}
\#S=n\,,
\eeq
and
\beq\label{propS1}
|\bigcap_{l\in S} I_{l}|\geq \frac{1}{2\,M^{n}}\,.
\eeq
\end{l1}

\begin{proof}
Wlog we can assume $n\geq 2$ since the case $n=1$ is trivial.

On the one hand, from H\"older's inequality, one has
\beq\label{Ho}
\I:=\int_{-\frac{1}{2}}^{\frac{1}{2}}\left(\sum_{l=1}^{N}\chi_{I_l}(x)\right)^n\,dx\geq \left(\int_{-\frac{1}{2}}^{\frac{1}{2}}\sum_{l=1}^{N}\chi_{I_l}(x)\,dx\right)^n\geq_{\eqref{propset}}\left(\frac{N}{M}\right)^n\,.
\eeq
On the other hand
\beq\label{Ho1}
\I=\sum_{S\subset \bigotimes_{n}\{1,\ldots, N\}} \left|\bigcap_{l\in S} I_l\right|=\sum_{r=1}^n\sum_{{S\subset \bigotimes_{n}\{1,\ldots, N\}}\atop{\#\tilde{S}=r}} \left|\bigcap_{l\in S} I_l\right|\,,
\eeq
where here if $S=(i_1, i_2,\ldots, i_n)$ we let $\tilde{S}:=\bigcup_{l=1}^n \{i_l\}$.

It is now simple to notice that
\beq\label{Ho2}
\sum_{{S\subset \bigotimes_{n}\{1,\ldots, N\}}\atop{\#\tilde{S}=r}} \left|\bigcap_{l\in S} I_l\right|\leq n!\,N^{r}\:.
\eeq
Thus, based on \eqref{Ho} - \eqref{Ho2}, we must have
\beq\label{Ho3}
\sum_{{S\subset \bigotimes_{n}\{1,\ldots, N\}}\atop{\#\tilde{S}=n}} \left|\bigcap_{l\in S} I_l\right|\geq
\left(\frac{N}{M}\right)^n-n!\sum_{r=1}^{n-1} N^r\geq \left(\frac{N}{M}\right)^n-n^n N^{n-1}
\eeq
Using now \eqref{n1} and pigeonhole principle we immediately conclude that \eqref{Ho3} implies \eqref{propS1}.

\end{proof}

\begin{l1}\label{coeffcontr} Let $d\in\N$ and let $\S,\,\A$ be two collections of strictly increasing real numbers with $\S=\{s_i\}_{i=1}^{2^{d}} \subset [1,4]$ and $\A=\{\a_j\}_{j=1}^d\subset\R$. Set\footnote{By convention, if $d=1$ we set $\Pi\A=1$.} $\|\S\|:=\min_{{1\leq i<j\leq 2^d}}|s_j-s_i|$ and $\Pi_j\A:=\prod_{{i=1}\atop{i\not=j}}^d |\a_j-\a_i|$.

Assume we are given $A>0$ and a sequence $\{a_j\}_{j=1}^d$ of real numbers such that
\beq\label{conts}
\left|\sum_{j=1}^d a_j s^{\a_j}\right|\leq A\:\:\:\:\:\textrm{for any}\:s\in \S\:.
\eeq
Then, for any $1\leq j\leq d$, we have
\beq\label{contcoef}
|a_j|\leq \frac{4^{d((\a_d-\a_1)+\max_{1\leq j\leq d} |\a_j|+2)}\,A}{\|\S\|^{d-1}\,\Pi_j\A}\:.
\eeq
\end{l1}

\begin{proof}
We will prove in fact the slightly improved version
\beq\label{contcoefw}
|a_j|\leq \frac{4^{d((\a_d-\a_1)+\min_{1\leq j\leq d} |\a_j|+2)}\,A}{\|\S\|^{d-1}\,\Pi_j\A}\:.
\eeq
This last statement will be verified by induction over the values of $d$.

The case $d=1$ is trivial. Assume now that \eqref{contcoefw} holds for $1,\,2,\,\ldots,\,d-1$.

We now multiply in \eqref{conts} by $s^{-\a_1}$ in order to deduce that\footnote{Here $\a_{-}=-\a$ if $\a<0$ and $0$ otherwise.}
\beq\label{conts1}
\left|\sum_{j=1}^d a_j s^{\a_j-\a_1}\right|\leq 4^{\a_{1-}}\,A\:\:\:\:\:\textrm{for any}\:s\in \S\:.
\eeq
Applying now the mean value theorem on each of the intervals $(s_{2i-1},\,s_{2i})$ with $i\in\{1,\ldots,\,2^{d-1}\}$ we obtain a collection of intermediate points $\S^{1}:=\{s^1_i\}_{i=1}^{2^{d-1}}$ such that
\beq\label{conts2}
\left|\sum_{j=2}^d a_j\,(\a_j-\a_1) s^{\a_j-\a_1}\right|\leq \frac{4^{\a_{1-}+2}\,A}{\|\S\|}\:\:\:\:\:\textrm{for any}\:s\in \S^{1}\:.
\eeq
Applying now the induction hypothesis, for any $2\leq j\leq d$, we get that
\beq\label{conts3}
|a_j|\,(\a_j-\a_1) \leq \frac{4^{(d-1)((\a_d-\a_2)+\min_{2\leq j\leq d}|\a_j-\a_1|+2))}}{\|\S^1\|^{d-2}\,\Pi_j\A^1}\,\frac{4^{\a_{1-}+2}\, A}{\|\S\|}\:,
\eeq
where here we set $\A^1=\{\a_j-\a_1\}_{j=2}^{d}$.

Notice now that from our construction of the sets $\A^1$ and $\S^{1}$ one has that $\|\S^1\|\geq \|\S\|$ and that
$(\a_j-\a_1)\,\Pi_j\A^1=\Pi_j\A$. With these one immediately verifies our induction hypothesis \eqref{contcoefw} for $2\leq j\leq d$.

Finally, in order to get the similar relation for the term $a_1$ it is enough to repeat the above reasoning by replacing  the role played by $\a_1$ with $\a_d$.
\end{proof}

\begin{o0}\label{xeropolyn} [\textsf{Number of real distinct roots for a generalized polynomial}] Assume we are given a generalized polynomial
\beq\label{genpolyna1}
P(t)=\sum_{k=1}^d a_k\,t^{\a_k}\,,
\eeq
with $d\in\N$, $\{a_k\}_{k=1}^{d}$ real numbers and $\{\a_k\}_{k=1}^d$ real positive numbers. Then, the number of real distinct roots of $P$ obeys the following relation
\beq\label{genpolyna2}
\#\{t\in\R_{+}\,|\, P(t)=0\}\leq 2 d\,.
\eeq
Notice that even in the case when $\{\a_k\}_{k=1}^d$ are (strictly increasing) positive integers, \eqref{genpolyna2} is far from obvious (assuming of course we are in the nontrivial case $\a_d>2d$), offering a refinement over the information one could get from the fundamental theorem of algebra: indeed, in this latter instance, applying the fundamental theorem of algebra one could only say that the LHS of \eqref{genpolyna2} is bounded by $\a_d$, while \eqref{genpolyna2} states that the relevant information stays in the \textit{number} of non-zero monomials as opposed to the degree of $P$.

The proof of this result is done via induction over $d$, and is inspired from \cite{swmul}. Moreover, this observation  served as a toy model in our approach for the proof of the more involved Lemma \ref{coeffcontr}.
\end{o0}

We pass now to the actual proof of our proposition. This will be based on reductio ad absurdum. Thus, we assume that \eqref{fhe0s} fails, and hence
\beq\label{fhe00}
\# \H_{\tilde{q},\d}\geq 2^{(\frac{1}{2}-\nu) m}\:,
\eeq
for a suitable small $\nu>0$.

We now define
\beq\label{lproj}
\H_{\tilde{q},\d}^{L}:=\{l\,|\,(l,w)\in \H_{\tilde{q},\d}\}\,,
\eeq
and let
\beq\label{lproj1}
\H_{\tilde{q},\d}^{L,S}:=\textrm{a maximal separated subset of}\:\:\H_{\tilde{q},\d}^{L}\:\:\,,
\eeq
such that
\beq\label{lproja1}
\#\H_{\tilde{q},\d}^{L,S}= 2^{2^{3d}\,2\d\,m}\,
\eeq
which is possible if $\frac{1}{2}-\nu-4\d> 2^{3d}\,2\d$, since, based on \eqref{A} - \eqref{H}, one has that for any $l\in \H_{\tilde{q},\d}^{L}$
\beq\label{lproja11}
\#\{w\,|\,(l,w)\in \H_{\tilde{q},\d}\}\lesssim 2^{4\,\d\, m}\,.
\eeq

Next, for each $l\in\H_{\tilde{q},\d}^{L,S}$, we let
\beq\label{lprojA}
A_l:=A_{\tilde{q}}(l,w_l)\:\:\textrm{with}\:\:w_l\:\:\textrm{defined by}\:|A_{\tilde{q}}(l,w_l)|=\max_{w:\:(l,w)\in\H_{\tilde{q},\d}}|A_{\tilde{q},\d}(l,w)|\:.
\eeq

Now, for\footnote{Recall here that $d$ stands for the number of ``generalized" monomials (\textit{i.e.} we allow the power/degree to be an arbitrary real number) of $q$.} $n= 2^{2d}$, we apply Lemma \ref{largeint} for $I_l:=A_l$, $M=2^{2\d m}$ and $N=\#\H_{\tilde{q},\d}^{L,S}=2^{2^{3d}\,2\d\,m}$ to deduce that there exist $\bar{\H}_{\tilde{q}}\subseteq \H_{\tilde{q},\d}^{L,S}$ with $\#\bar{\H}_{\tilde{q}}=n$ and a set $X:=\bigcap_{l\in\bar{\H}_{\tilde{q}}}A_l$ such that
\beq\label{lproja110}
\eeq
\begin{itemize}
\item $|\tilde{q}(x,x+l)-w_l|\leq 2^{-(\frac{1}{2}-2\d)\,m}$ for any $l\in \bar{\H}_{\tilde{q}}$ and $x\in X$;

\item $|X|\gtrsim 2^{-2^{2d}\,2\d\,m}$.
\end{itemize}

First thing to notice is that, from \eqref{lproja110} and Lemma \ref{coeffcontr}, one immediately has
\beq\label{contcoefcont}
\max_{1\leq j\leq d}|a_j(x)|\leq \frac{4^{d((\a_d-\a_1)+\max_{1\leq j\leq d} |\a_j|+2)}\,2}{\|\bar{\H}_{\tilde{q}}\|^{d-1}\,\inf_{j}\prod_{{i=1}\atop{i\not=j}}^d |\a_j-\a_i|}\:,
\eeq
for any $x\in X$.

Exploiting the first item in \eqref{lproja110}, we notice that for any $x,\,x_0\in X$ one has
\beq\label{erw}
|\sum_{j=1}^d a_j(x)\,(x+l)^{\a_j}-\sum_{j=1}^d a_j(x_0)\,(x_0+l)^{\a_j}|\lesssim 2^{-(\frac{1}{2}-2\d)\,m}\:\:\:\forall\: l\in\bar{\H}_{\tilde{q}}\,.
\eeq
Fix $x,\,x_0\in X$. We now apply Taylor's formula in order to deduce that
\beq\label{Tay}
(x_0+s)^{\a_j}=(x+s)^{\a_j}\,+\,\a_j\,(x_0-x)\,(x+s)^{\a_j-1}\,+\,O(|x-x_0|^2)\,.
\eeq
Putting together \eqref{contcoefcont}, \eqref{erw} and \eqref{Tay}, we further deduce that\footnote{Here $O_{\A}(B)$ designates an expression that is bounded from above by a quantity of the form $C_{\A}\times B$ where $C_{\A}>0$ depends on the properties of the set $\A=\{\a_j\}_{j=1}^d$.}
\beq\label{erw1}
\eeq
$$\left|\sum_{j=1}^d \big\{(a_j(x)- a_j(x_0))\,(x+l)^{\a_j}-\a_j\,(x_0-x)\,a_j(x_0)\,(x+l)^{\a_j-1}\big\}\right|$$
$$\lesssim 2^{-(\frac{1}{2}-2\d)\,m}\,+\,\frac{1}{\|\bar{\H}_{\tilde{q}}\|^{d-1}}\,O_{\A}(|x-x_0|^2)\:\:\textrm{for any}\:\: l\in\bar{\H}_{\tilde{q}}\,.$$
Applying now Lemma \ref{coeffcontr} for the coefficient of the ``lowest degree generalized monomial", we deduce that
\beq\label{a1}
|a_1(x)|,\,|a_1(x_0)|\lesssim  \frac{4^{8d(\max_{1\leq j\leq d} |\a_j|+1)}  (2^{-(\frac{1}{2}-2\d)\,m}\,+\,\frac{1}{\|\bar{\H}_{\tilde{q}}\|^{d-1}}\,|x-x_0|^2)}{|x-x_0|\,|\a_1|\,\
|\bar{\H}_{\tilde{q}}\|^{2d-1}\,\prod_{i=2}^d |\a_i-\a_1|}\:.
\eeq
Once at this point we set $\I_1$ be the collection of same length $=2^{-\frac{m}{4}}$ intervals partitioning the interval $[1,4]$. Set now $\I_1^{0}$ be the collection of intervals $J\in \I_1$ for which
\beq\label{a2}
|X\cap J|\geq 2^{-5}\,2^{-2^{2d}\,2\d\,m}\,|J|\,,
\eeq
and set
\beq\label{a3}
X_1:=\bigcup_{J\in \I_1^{0}} X\cap J\;.
\eeq
From \eqref{lproja110}, \eqref{a2} and \eqref{a3} one deduces immediately that
\beq\label{a4}
|X_1|\geq \frac{|X|}{2}\geq 2^{-2}\,2^{-2^{2d}\,2\d\,m}\;.
\eeq
Applying now \eqref{a1} within each set $X\cap J$ with $J\in \I_1^{0}$ we notice for any $x\in X\cap J$ we can choose $x_0\in X\cap J$ such that $|x-x_0|\geq 2^{-10}\,2^{-2^{2d}\,2\d\,m}\,|J|$ and thus we have that
\beq\label{a5}
|a_1(x)|\lesssim  \frac{4^{8d(\max_{1\leq j\leq d} |\a_j|+1)}  (2^{-(\frac{1}{4}-2\d)\,m}\,2^{2^{2d}\,2\d\,m})}{|\a_1|\,\|\bar{\H}_{\tilde{q}}\|^{3d-2}\,\prod_{i=2}^d |\a_i-\a_1|}\:\;\:\:\textrm{for any}\:x\in X_1\:.
\eeq
Inserting now \eqref{a5} in \eqref{erw} we have that for any $l\in\bar{\H}_{\tilde{q}}$ and $x,\,x_0 \in X_1$ one has
\beq\label{erw10}
|\sum_{j=2}^d a_j(x)\,(x+l)^{\a_j}-\sum_{j=2}^d a_j(x_0)\,(x_0+l)^{\a_j}|\lesssim \frac{4^{8d(\max_{1\leq j\leq d} |\a_j|+1)}  (2^{-(\frac{1}{4}-2\d)\,m}\,2^{2^{2d}\,2\d\,m})}{|\a_1|\,\|\bar{\H}_{\tilde{q}}\|^{3d-2}\,\prod_{i=2}^d |\a_i-\a_1|}\,.
\eeq
At this point we can repeat the above algorithm to deduce an upper bound for the size of the terms $a_2(x)$ where here $x\in X_2$ with $X_2\subset X_1$ such that $|X_2|\geq \frac{1}{2}\,|X_1|$. Iterating this argument $d$ times we conclude that for any $1\leq k\leq d$ and $x\in X_d$ with $X_d\subset X$ and $|X_d|\geq \frac{1}{2^{d+1}}\,|X|$, one must have
\beq\label{a6}
|a_k(x)|\lesssim  \frac{4^{8d\,k(\max_{1\leq j\leq d} |\a_j|+1)}  (2^{-(\frac{1}{2^{k+1}}-2\d)\,m}\,2^{2^{2d}\,k\,2\d\,m})}{|\a_1|\cdots|\a_k|\,\|\bar{\H}_{\tilde{q}}\|^{(3d-2)k}\,
\prod_{j=1}^{k}\prod_{{i=1}\atop{i\not=j}}^d |\a_i-\a_j|}\:.
\eeq
Notice now that $\|\bar{\H}_{\tilde{q}}\|\geq 2^{-(2^{3d}\,2\d+\nu)\,m}$. Thus, from \eqref{a6}, we deduce that for any $x\in X_d$ and $1\leq k\leq d$ one has
\beq\label{a7}
|a_k(x)|\lesssim  \frac{4^{8d^2\,(\max_{1\leq j\leq d} |\a_j|+1)}}{|\a_1|\cdots|\a_k|\,
\prod_{j=1}^{k}\prod_{{i=1}\atop{i\not=j}}^d |\a_i-\a_j|} \,2^{-(\frac{1}{2^{d+1}}-2\d)\,m}\,2^{2^{2d}\,d\,2\d\,m}\,2^{3d^2(2^{3d}\,2\d+\nu)\,m}\:.
\eeq
Choosing now $\d,\,\nu>0$ small enough such that
\beq\label{a8}
\frac{1}{2^{d+1}}-2\d-2^{2d}\,d\,2\d-3d^2(2^{3d}\,2\d+\nu)>0\,,
\eeq
we deduce that for any $x\in X_d$ and $l\in \bar{\H}_{\tilde{q}}$ one has
\beq\label{a9}
\left|\frac{\partial}{\partial t}q(x,t)\big|_{t=x+l}\right|<<1\,,
\eeq
which violates the definition in \eqref{PY1}.

Thus  we reach our desired contradiction in \eqref{fhe00} ending the proof of Proposition \ref{fewhe}.
\end{proof}

\begin{o0}\label{compconst}
A simple inspection of the proof of Proposition \ref{fewhe} shows that one can choose the following values for our parameters: $\d=\frac{1}{2^{5d+5}}$ and $\nu=\frac{1}{2^{2d+5}}$. This choice implies the desired exponential decay in $m$ in \eqref{ndeg00e} with $\bar{\ep}=\d=\frac{1}{2^{5d+5}}$. In the case in which $q(x,t)$ is a polynomial in $t$ with no constant coefficient - \textrm{i.e}  $\a_j\in\N$ in  \eqref{PY1} - the optimal dependence in \eqref{ndeg00e} as well as in \eqref{dualityop} is in fact linear; that is, instead of the above exponential dependence, with supplementary ideas that we choose not to present here, one can obtain in the polynomial case that in fact $\ep,\,\bar{\ep}\approx\frac{1}{d}$.
\end{o0}

\section{Proofs of the remaining corollaries}\label{Remcorol}

$\newline$
\noindent\textbf{Proof of Corollary \ref{Polyncasegen}.}
$\newline$

This follows from Main Theorem and the proof of Theorem \ref{Gencurvpolyn} by inspecting the bound dependencies on the exponents $\{\a_j\}_j$. More specifically, the form of \eqref{lpbddep} follows from tracking the constant/bounds dependencies in \eqref{defjkd0} - \eqref{asymptotic011}, \eqref{contcoef} and \eqref{a7} together with the $L^p$ bound interpolation argument provided in Lemma \ref{translk}.

$\newline$
\noindent\textbf{Proof of Corollary \ref{Polyncaseap}.}
$\newline$

Corollary \ref{Polyncaseap} is a straightforward application of Corollary \ref{Polyncasegen}. For completeness, we provide here the argumentation: from Corollary \ref{Polyncasegen} we know that for $d\in\N$, $\{\a_j\}_{j=1}^d\subset \mathbb{R}\setminus\{0,\,1\}$ and $\{a_j\}_{j=1}^d$ measurable functions defining
$\g(x,t):=\sum_{j=1}^{d} a_j(x)\,t^{\a_j}$ we have that $C_{\g}$ is $L^p$ bounded for $1<p<\infty$.

Applying now a standard linearization argument we notice that
 \beq\label{SWgenlin}
C_{\vec{\a},d}f(x)\approx \int_{\R} f(x-t)\, e^{i\,\sum_{l=1}^d a_l(x)\,t^{\a_l}}\,\frac{1}{t}\,dt=C_{\g}\,,
\eeq
finishing our proof.

$\newline$
\noindent\textbf{Proof of Corollary \ref{Crossprod}.}
$\newline$

It is straightforward to see that the curve $\g(x,t)=u(x)\,\tilde{\g}(t)$ defined in the statement of Corollary \ref{Crossprod}, see in particular relations \eqref{crp1} - \eqref{crfstterm0}, verifies relations \eqref{u} - \eqref{fstterm0}. The only part left is to verify the nondegeneracy condition \eqref{ndeg0} - which in our present context becomes

\beq\label{ndeg0co}
\sup_{{k,\,n\in\Z}\atop{|j|\geq |j_0|}} \int_{1<|s|<4}\,\sup_{t\in\R}\,\left(\frac{1}{2^j}\int_{(k-\frac{1}{2})\,2^{-j}}^{(k+\frac{1}{2})2^{-j}}
\frac{\chi_{\R_\a}(x)\,\chi_{\Z^{x}}(j)\,\phi(\frac{\g'_x(2^{-j})}{2^n})}{\left\lfloor2^{\frac{m}{2}}
\left(\tilde{q}_j(s+2^{j}\,x-k)-t\,\frac{2^n}{\g'_x(2^{-j})}\right)\right\rfloor^2}\,dx\right)\,ds\lesssim_{\g} 2^{-2\,\bar{\ep}\,m}\:,
\eeq
where here we used that our hypothesis implies $q_j(x,t)=\tilde{q}_j(t):=\frac{\tilde{\g}'(2^{-j}\,t)}{\tilde{\g}'(2^{-j})}$.

Combining the reasoning from the proof of Lemma \ref{nondegcont} with those from the proof of Proposition \ref{fewhe} we see that we can reduce \eqref{ndeg0co} to proving the analogue of \eqref{fhe0} with $\H_{\tilde{q},\ep}$ there replaced by $\H_{\tilde{q}_j,\g,\d}$, where here
\beq\label{Hgam}
\H_{\tilde{q}_j,\g,\d}:=\{(l,w)\in \mathcal{R}_m^2\,|\,|A_{\tilde{q}_j,\g,\d}(l,w)|> 2^{-2\d m}\}\:,
\eeq
with
\beq\label{Agam}
A_{\tilde{q}_j,\g,\d}(l,w):=\{x\in [-\frac{1}{2},\,\frac{1}{2}]\,|\,|\tilde{q}_j(l\,+\,x)-w\,a(x)|\leq 2^{-(\frac{1}{2}-2\d)\,m}\}\:,
\eeq
and $a(\cdot)$ a suitable measurable function (depending on $u$) with $\frac{1}{10}\leq |a(x)|\leq 10$ for any $x\in [-\frac{1}{2},\,\frac{1}{2}]$.

Following now the same steps as in the proof by contradiction of Proposition \ref{fewhe}, specifically  \eqref{fhe00} - \eqref{lproja110}, by taking there $d=1$, one gets the modified form of \eqref{lproja11} as
\beq\label{lproja11m}
\eeq
\begin{itemize}
\item $|\tilde{q}_j(x+l)-w_l a(x)|\leq 2^{-(\frac{1}{2}-2\d)\,m}$ for any $l\in \bar{\H}_{\tilde{q}_j}$ and $x\in X$;

\item
$\# \bar{\H}_{\tilde{q}_j}\geq 2$ and $\min_{l_1\not=l_2\in \bar{\H}_{\tilde{q}_j}}|l_1-l_2|\geq 2^{-(16\,\d+\nu)\,m} $;

\item $|X|\geq 2^{-8\,\d\,m}$;
\end{itemize}
where here $\nu,\, \d\in (0,1)$ suitable chosen obeying \eqref{a8} for $d=1$.

This immediately implies that that there exists $x,\,y\in X$ with $|x-y|\geq 2^{-8\,\d\,m}$ and $l_1,\,l_2\in \bar{\H}_{\tilde{q}_j}$ with $|l_1-l_2|\geq 2^{-(16\,\d+\nu)\,m}$ such that
\beq\label{lproja11m1}
|\frac{\tilde{q}_j(x+l_k)}{a(x)}-\frac{\tilde{q}_j(y+l_k)}{a(y)}|\lesssim 2^{-(\frac{1}{2}-2\d)\,m}\:\:\:\textrm{for}\:\:\:k\in\{1,\,2\}\:.
\eeq
Now according to the two cases under discussion we have:

\noindent i) if $u=u_0$ (almost everywhere) constant then $a=a_0$ constant and then from the mean value theorem and condition \eqref{crfstterm0} we have that \eqref{lproja11m1} can't hold.

\noindent ii) for general measurable $u$ we deduce from \eqref{lproja11m1} that
\beq\label{lproja11m2}
|\frac{\tilde{q}_j(x+l_1)}{\tilde{q}_j(y+l_1)}-\frac{\tilde{q}_j(x+l_2)}{\tilde{q}_j(y+l_2)}|\lesssim 2^{-(\frac{1}{2}-2\d)\,m}\:.
\eeq
However, an application of the mean value theorem together with hypothesis \eqref{nondeg} gives
\beq\label{lproja11m3}
|\frac{\tilde{q}_j(x+l_1)}{\tilde{q}_j(y+l_1)}-\frac{\tilde{q}_j(x+l_2)}{\tilde{q}_j(y+l_2)}|\gtrsim M |l_1-l_2|\,|x-y|\:,
\eeq
which contradicts \eqref{lproja11m2} for a proper choice of $\nu,\, \d$.

$\newline$
\noindent\textbf{Proof of Corollary \ref{Crossprod1}.}
$\newline$

This is straightforward based on Corollary \ref{Crossprod} above an the fact that $\n\f\subset \mathbf{NF}$. The last containment is immediate by just inspecting the definitions of the classes $\mathbf{NF}$ (in the present paper) and  $\n\f$
(in \cite{lv4}); indeed, one has
\begin{itemize}
\item condition \eqref{crp1} here is the same as conditions (2) (for $j\in \Z_{+}$) and respectively (8) (for $j\in \Z_{-}$) in \cite{lv4};

\item conditions \eqref{crfstterma0} and \eqref{crfstterm0} here are more general than (but analogues to) conditions
(3), (4), (5) and respectively (9), (10), (11) in  \cite{lv4};

\item condition \eqref{nondeg} is morally equivalent with (6) and respectively (12)  in  \cite{lv4}.
 \end{itemize}

\section{Final Remarks}

In this section we start by clarifying some aspects related to the newly introduced class of curves $\mathbf{M}_x\mathbf{NF}_{t}$ and then we end with some open problems that have arisen naturally in the course of the present study.

\subsection{Analyzing the class $\mathbf{M}_x\mathbf{NF}_{t}$}

As always when introducing a new concept or definition around which an entire paper revolves, an honest discussion is advisable along at least three directions:
\begin{itemize}
\item How natural/intuitive is the new definition?

\item What are its merits?

\item What are its limitations?
\end{itemize}

\subsubsection{The intuition} We first stress that Definition \ref{defgam} is tailored around the concept of curvature which is the keystone in each of the three themes approached in our study. The prototype that served as an initial model for our class of curves (and hence must end up residing within it) is the set of all $t-$polynomials of a given degree $d\in\N$, $d\geq 2$
\beq\label{poy}
P(x,t)=\sum_{k=1}^d a_k(x)\,t^{k+1}\,,
\eeq
where here $\{a_k(\cdot)\}_{k=1}^d$ are only assumed to be real measurable functions.

Two features are quintessential in the above expression: i) one allows only measurability in $x$; ii) one does not allow linear terms, and thus a suitable non-zero curvature must be present.

With this in mind, we can now easily notice that \eqref{mpp} expresses in more specific terms i) while preparing the ground for ii). Conditions \eqref{u} - \eqref{deriv0} provide a natural splitting of the $(x,t)-$plane into regions where we gain some more structure and control on the properties of $\g$ in analogy - it is instructive here to consult Step 1 in Section \ref{richclass} - with a partition guided by the roots or coefficient localization of a polynomial with variable coefficients given by \eqref{poy}. Within each such region one is able to extract a first, qualitative curvature condition reflected in \eqref{deriv0} thus partially fulfilling ii) above.

Condition \eqref{variation0} is necessary in order to impose the standard almost disjointness of the Littlewwod - Paley projections in the $t-$variable for a fixed $x$ without which no square-function type argument would be applicable. \footnote{This condition is trivially satisfied for any polynomial with variable coefficients - even if allowing constant and linear terms.}

We pass now to the two fundamental requirements for our class of curves:
\begin{itemize}
\item \textsf{Doubling and uniform non-zero curvature (non-flatness)}:  The underlying motivation for the doubling condition \eqref{asymptotic0} stems from the natural desire to exploit the dilation symmetry of the kernel $\frac{1}{t}$ evoked by the Whitney decomposition \eqref{wit}. This in turn invites the multiplier discretization \eqref{firstlt} and thus the analysis of its corresponding phase \eqref{firstl1}, with the preeminent role played by its derivative \eqref{firstl2}. Indeed, the behavior of \eqref{firstl2} is essential for the location of the stationary points, which in turn guides our entire subsequent approach. At this point, we insert a simplification of our analysis by assuming that \eqref{firstl3} is a good approximation for \eqref{firstl2}, or equivalently, by morally rephrasing our doubling condition \eqref{asymptotic0}. From this point on, the uniform  upper-boundedness condition \eqref{fstterma0} is the doubling condition extended uniformly across all scales while \eqref{fstterm0} is the completion of ii) in terms of a uniform - in scale - non-zero curvature condition.

\item \textsf{Non-degeneracy}: Finally, condition \eqref{ndeg0} is a dilation invariant way of measuring the ``twisted non-zero curvature" of $\g(x,t)$. More precisely, in our context it is not enough simply to control the behavior and
  curvature of $\g$ in t. Rather, we need to consider a relative interplay between $x$, $t$, and how (the derivative of) $\g$ varies around a given point.\footnote{For more on this, please see the discussion within the limitations in the definition of our class of curves.}

\end{itemize}

\subsubsection{The merits}\label{merits}   The primary advantage of the $\mathbf{M}_x\mathbf{NF}_{t}-$definition lies within its malleability, demonstrated by the large classes of examples that it subsumes. Indeed, following Corollary \ref{Crossprod} and the definition therein, Observation \ref{W}, Theorem \ref{Gencurvpolyn} and its restatement \eqref{py0} - \eqref{py1}, we do have that
\begin{eqnarray}\label{keyincl}
& \mathbf{NF}\subset \mathbf{M}_x\mathbf{NF}_{t}\,,\\
& \mathbf{M}_x\cdot\n\f:=\{u(x)\,\g(t)\,|\,u\:\textrm{measurable}\:\g\in \n\f\} \subset \mathbf{M}_x\mathbf{NF}_{t}\,,\\
& \p_d\subset \mathbf{M}_x\mathbf{NF}_{t}\:.
\end{eqnarray}
Moreover, by repeating step by step the proof in Section \ref{richclass} - with the obvious analogue of Lemma \ref{coeffcontr} - one can extend the last inclusion above, by actually showing that any expression of the form
\beq\label{plog}
P(x,t)=\sum_{k=1}^d a_k(x)\,t^{\a_k}\,\log^{\b_k} |t|\,\in\mathbf{M}_x\mathbf{NF}_{t}\,,
\eeq
with $d\in\N$, $\{a_k(\cdot)\}_{k=1}^d$ measurable, $\{\a_k\}_{k=1}^d \subset\mathbb{R}\setminus \{0,\,1\}$ and $\{\b_k\}_{k=1}^d\subset\mathbb{R}$.

Thus, indeed, many of the most relevant and/or standard examples that one might consider are proven to be part of our newly defined class or curves.

\subsubsection{The limitations}\label{lim} The class $\mathbf{M}_x\mathbf{NF}_{t}$ has several downsides that will be discussed below:

\begin{itemize}

\item \textsf{Exclusion of $t-$non-doubling curves}: The doubling condition \eqref{asymptotic0} is clearly restrictive and - unlike the non-zero curvature condition - should not play an important role in the boundedness of the operators that we consider. In fact, it is known that both $H_{\G}$ and $M_{\G}$ are $L^2(\R^2)$-bounded operators for $\G$ corresponding to $\g(x,t)=x\,e^{-\frac{1}{t^2}}$, (\cite{CWWhf}), or more generally $\g(x,t)=P(x)\,e^{-\frac{1}{t^2}}$ with $P$ polynomial (\cite{Benh}). Although we believe that our results could in principle be extended to include these latter situations, this would require a non-trivial amount of extra technicalities that would not be justified in the context of our main focus derived from the inclusions in \eqref{keyincl} - \eqref{plog}.

\item \textsf{Absence of non-translation invariance}: A far more consequential restriction is the lack of translation invariance for $\mathbf{M}_x\mathbf{NF}_{t}$. While unappealing, this should not be surprising since our definition for the set of curves relies fundamentally on the concept of \textit{non-zero curvature near zero and infinity}, a property that is compatible with dilation but not with translation symmetry. Consequently, we remark that
    \begin{eqnarray}\label{nottrinv}
& \g(x,t)\in \mathbf{M}_x\mathbf{NF}_{t}\:\:\nRightarrow\:\: \g(x,t-a)\in \mathbf{M}_x\mathbf{NF}_{t}\:\:\:\textrm{for}\:\:\:a\in\R\setminus\{0\}\,,\\
& \label{trx}\g(x,t)\in \mathbf{M}_x\mathbf{NF}_{t}\:\:\nRightarrow\:\: \g(x,t-x)\in \mathbf{M}_x\mathbf{NF}_{t}\:.
\end{eqnarray}

\item \textsf{Sufficiency versus necessity}: We return to the two key conditions i) doubling, uniform boundedness and non-zero curvature (non-flatness) and ii) non-degeneracy. While sufficient if taken together, neither of these two are necessary conditions. With regard to i) we already clarified this in discussing the first limitation above. Passing now to ii), we exemplify the lack of necessity by appealing to \eqref{trx}: on the one hand, while $\g_0(x,t)=t^2\in \mathbf{M}_x\mathbf{NF}_{t}$ trivially, one can easily check that $\g(x,t):=\g_0(x,x-t)=(t-x)^2$ does not satisfy \eqref{ndeg0} and hence $\g\notin\mathbf{M}_x\mathbf{NF}_{t}$; on the other hand, taking as usual $\Gamma\equiv\{(t,\,-\g(x,t))\}_{\{x\in\R\}}$ we have that\footnote{In what follows we ignore the principal value symbol.}
    $$H_{\G}(f)(x,y)=\int_{\R} f(x-t,\,y+\g(x-t))\,\frac{dt}{t}=_{c.v.\:t\,\rightarrow\,x-t}\int_{\R} f(t,\,y+\g(t))\,\frac{dt}{x-t}$$
    $$=:\int_{R}f_y(t)\,\frac{dt}{x-t}=H(f_y)(x)\:,$$
    where here $H$ stands for the standard one dimensional Hilbert transform.

    Thus, for any $1<p<\infty$, one deduces that
    \beq\label{hbd}
    \int_{\R^2}|H_{\G}(f)(x,y)|^p\,dx\,dy=\int_{\R^2}|H(f_y)(x)|^p\,dx\,dy\lesssim_{p}\int_{\R^2}|(f_y)(x)|^p\,dx\,dy=\|f\|_p^p\,.
    \eeq

\end{itemize}

\subsection{More about the non-degeneracy condition \eqref{ndeg0}}

As mentioned earlier, the non-degeneracy condition \eqref{ndeg0} is a dilation invariant way of measuring the ``twisted non-zero curvature" of $\g(x,t)$. Of course, following the proof, one might find various other - essentially equivalent or slightly more general - conditions that could have replaced \eqref{ndeg0} with the same consequences; indeed one could have asked as an alternative non-degeneracy hypothesis relation \eqref{fhe0} or,  even more generally, condition \eqref{ndegeq0}. However, all of these fundamentally ask for the same type of behavior - a suitable decay condition based on the interplay between the $x$ and $t$ parameters.

Passing now to more concrete aspects related to our non-degeneracy condition, we first comment on \eqref{ctrw1}. Although it is tempting to hope that one could replace $1$ by a factor of the form $2^{-\ep m}$ in the right-hand side of \eqref{ctrw1}, this is in fact impossible in general, as one can simply check for example in the case $q(x,t)=a(x)\,t$. This is our motivation for introducing the definitions of light and heavy pairs in \eqref{L} and \eqref{H}.

Finally, we briefly touch on the specific choice of where to place the supremum in $t$ in \eqref{ndeg0}. Indeed, it would have been natural to try to relax this condition more and move the $\sup_{t}$ expression in front of the first integral sign; however, it turns out that such a condition would hold for any $\g$ obeying \eqref{asymptotic0} - \eqref{fstterm0} and even for zero-curvature curves (\textit{e.g.} $\g(x,t)=a(x) t$ with $a(\cdot)$ measurable), thus being of no help. In this context, one could have tried to move the $\sup_{t}$ inside the second integral, but in this situation no decay in $m$ would be possible. Accordingly, the chosen location of the supremum in $t$ is the only one possible.

\subsection{Some open questions}

We end our paper with several open questions starting with the more pedestrian ones and slightly moving towards those with more philosophical content.

\subsubsection{An interpolation question - Generalized Vandermonde matrices}

A positive answer to the question below provides a more elegant alternative and also more efficient bounds to Lemma \ref{coeffcontr}:
$\newline$

\noindent \textbf{Open problem} [\textsf{Control over generalized Vandermonde determinants}]

\textit{Let $d\in\N$ and let $\S,\,\A$ be two collections of strictly increasing positive real numbers with $\S=\{s_i\}_{i=1}^{d}$ and $\A=\{\a_i\}_{i=1}^d$. Assume there exists $C_0>0$ such that the following holds:
\beq\label{bdbl}
\|\A\|,\,\|\S_0\|\geq C_0\:,
\eeq
where here $\S_0=\S\cup \{0\}$. $\newline$
\indent Define the Generalized Vandermonde matrix relative to the sets $\S,\,\A$ as
\beq\label{gvd}
\V_{\S,\A}=
\begin{bmatrix}
    s_1^{\a_1} & s_2^{\a_1}&\dots  & s_n^{\a_1}\\
    s_1^{\a_2} & s_2^{\a_2}& \dots & s_n^{\a_2} \\
    \vdots & \vdots &  \ddots & \vdots \\
    s_1^{\a_n}&s_2^{\a_n} &\dots & s_n^{\a_n}
\end{bmatrix}\:.
\eeq
\indent Is is true that there exists a constant $C=C(d,\,C_0)>0$ depending solely on $d$ and $C_0$ but not on the other characteristics of the sets $\S$ and $\A$ such that the following holds
\beq\label{bdblq}
|\det \V_{\S,\A}|\geq C\,?
\eeq}

\subsubsection{In relation with $\mathbf{M}_x\mathbf{NF}_{t}$.}
i) \textsf{Relaxing the requirements (A):} It would be interesting to investigate whether the non-degeneracy condition \eqref{ndeg0} is in fact needed for our Main Theorem to hold. As we have already seen in Subsection \ref{lim} we know that \eqref{ndeg0} is not a necessary condition; accordingly, conditions \eqref{mpp} - \eqref{fstterm0} in Definition \ref{defgam} alone might in fact be sufficient for our main results. A similar question can be addressed towards the - already mentioned as unnecessary - doubling requirement.

ii) \textsf{Relaxing the requirements (B):}  As mentioned in Definition \ref{defgam}, it is enough to require $C^{2+}$-regularity in $t$. While our proof was performed assuming $C^{4}-$regularity this was just for exposition purposes - as the only key fact required is the $L^1$-summability of each of the components in \eqref{Keyestim2} relative to the parameters $l,\,l_1$ and $s$ (or $s_1$) which is satisfied if the power exponent for each of the fractions is any number strictly greater than one (as opposed to two in \eqref{Keyestim2}). In this context, one possible direction of investigation would be to see if one could lower the $t-$regularity requirement for our curves to the limiting case 2.

iii) \textsf{End-point bounds for $C_\g$ and $M_\g$:} The presence of non-zero curvature and the highly oscillatory nature of the phase creates a good context for the following natural question: what is the behavior of $C_{\g}$ and $M_{\g}$ near $L^1$? For both operators it is expected for one to see an improvement over the known estimates from the flat (zero-curvature) case.

\subsubsection{Some final questions}

1) Returning to the historical evolution to of our theme (I), can our main results regarding the $L^p$-boundedness of $H_{\G}$ and $M_{\G}$ be given any interesting interpretation in terms of parabolic differential operators with variable coefficients?

2) Following the detailed discussion in Section \ref{Dich}, it would be of real interest if one could remove the non-zero curvature condition (as well as the non-degeneracy one) in order for our Main Theorem, Part (I) to cover the situation when $\g$ is given by a polynomial with linear term included and having $x-$measurable coefficients, \textit{i.e.} $\g(x,t)=\sum_{j=1}^{n} a_j(x)\,t^j$ with $\{a_j(\cdot)\}_j$ arbitrary real measurable functions.

3) Finally, one could be tempted to ask even more generally for similar boundedness results for analogous operators $H_{\G}$ and $M_{\G}$ in the situation in which we allow an extra $y-$dependence, \textit{i.e.} $\g=\g(x,y,t)$, and, as before, impose minimal regularity in $x,\,y$ and some suitable non-zero curvature and non-degeneracy conditions in $t$.

\end{document}